\documentclass[a4paper, english]{smfart}

\usepackage[utf8]{inputenc}	

\usepackage{url}

\usepackage[small]{caption}

\usepackage[all]{xy}
\usepackage{pb-diagram, pb-xy}

\usepackage{amssymb}

\usepackage{smfthm}

\usepackage{mathrsfs}

\numberwithin{equation}{section}

\renewcommand{\thesubsection}{\arabic{section}.\alph{subsection}}

\usepackage{hyperref}

\theoremstyle{plain}
\newtheorem{thm}[equation]{Theorem}
\newtheorem{lem}[equation]{Lemma}
\newtheorem{cor}[equation]{Corollary}
\newtheorem{propo}[equation]{Proposition} 

\theoremstyle{definition}
\newtheorem{rem}[equation]{Remark}
\newtheorem{ex}[equation]{Example}

\let\tilde = \widetilde

\newcommand{\bbA}{{\mathbb A}}
\newcommand{\bbC}{{\mathbb C}}

\newcommand{\bbG}{{\mathbb G}}

\newcommand{\bbN}{{\mathbb N}}

\newcommand{\bbR}{{\mathbb R}}

\newcommand{\bbZ}{{\mathbb Z}}

\newcommand{\rmD}{{\mathrm D}}
\newcommand{\rmM}{{\mathrm M}}
\newcommand{\rmS}{{\mathrm S}}
\newcommand{\rmT}{{\mathrm T}}

\newcommand{\Ad}{{\mathscr A\!d}}

\newcommand{\scrC}{{\mathscr C}}
\newcommand{\scrD}{{\mathscr D}}
\newcommand{\scrE}{{\mathscr E}}
\newcommand{\scrF}{{\mathscr F}}
\newcommand{\scrG}{{\mathscr G}}
\newcommand{\scrH}{{\mathscr H}}
\newcommand{\scrI}{{\mathscr I}}

\newcommand{\scrL}{{\mathscr L}}
\newcommand{\scrM}{{\mathscr M}}
\newcommand{\scrN}{{\mathscr N}}
\newcommand{\scrO}{{\mathscr O}}
\newcommand{\scrP}{{\mathscr P}}
\newcommand{\scrR}{{\mathscr R}}

\newcommand{\scrT}{{\mathscr T}}

\newcommand{\frakS}{{\mathfrak S}}

\newcommand{\Dbcoh}{\mathrm{D}^b_\mathrm{coh}}
\newcommand{\Dbhol}{\mathrm{D}^b_\mathrm{hol}}

\newcommand{\Mcoh}{\mathrm{Coh}}
\newcommand{\Mhol}{\mathrm{Hol}}
\newcommand{\Vect}{\mathrm{Vect}}

\newcommand{\MTM}{\mathrm{MTM}}
\newcommand{\MT}{\mathrm{MT}}

\newcommand{\Char}{\mathrm{Char}}
\newcommand{\Supp}{\mathrm{Supp}}
\newcommand{\Stab}{\mathrm{Stab}}

\newcommand{\one}{{\mathbf{1}}}

\newcommand{\an}{{\mathrm{an}}}
\newcommand{\semisimple}{{\mathrm{ss}}}

\DeclareMathOperator{\Gl}{GL}
\DeclareMathOperator{\Sp}{Sp}
\DeclareMathOperator{\SO}{SO}

\DeclareMathOperator{\Sl}{SL}
\DeclareMathOperator{\Pic}{Pic}
\newcommand{\Hom}{{\mathrm{Hom}}}
\newcommand{\End}{{\mathrm{End}}}

\newcommand{\Aut}{{\mathit{Aut}}}
\newcommand{\id}{{\mathit{id}}}

\newcommand{\p}{{(\mathrm{p})}}

\newcommand{\omegaFM}[1]{{\omega_{\hspace*{0.1em}\FM, #1}}}
\newcommand{\omegaMicro}[1]{{\omega_{\hspace*{0.1em} F_{#1}}}}

\newcommand{\CartierDual}[1]{{{\mathrm{Hom}}(#1, \bbG_m)}}

\newcommand{\RHom}{{R\mathscr{H}\!\mathit{om}}}
\newcommand{\scrHom}{{\mathscr{H}\kern -.9pt om}}
\newcommand{\scrExt}{{\mathscr{E}\kern -.9pt xt}}
\newcommand{\scrTor}{{\mathscr{T}\kern -.9pt or}}

\DeclareMathOperator{\Rep}{{Rep}}
\DeclareMathOperator{\Mod}{{Mod}}
\DeclareMathOperator{\LS}{LS}
\DeclareMathOperator{\VB}{VB}
\DeclareMathOperator{\Gal}{Gal}
\DeclareMathOperator{\D}{{D}}
\DeclareMathOperator{\DM}{{DM}}
\DeclareMathOperator{\M}{{M}}

\DeclareMathOperator{\MM}{{MM}}
\DeclareMathOperator{\Coh}{{Coh}}

\DeclareMathOperator{\CC}{{CC}}
\DeclareMathOperator{\FM}{{FM}}

\newcommand{\DR}{\mathrm{DR}}
\newcommand{\Dol}{\mathrm{Dol}}
\DeclareMathOperator{\Sym}{{Sym}}

\newcommand\circled[1]{{\footnotesize \textcircled{{\tiny #1}}}}

\renewcommand{\theenumi}{\alph{enumi}}
\renewcommand\labelenumi{(\theenumi)}

\author{Thomas Kr\"amer}
\address{Institut f\"ur Mathematik, Humboldt-Universit\"at zu Berlin\\ Unter den Linden 6, 10099 Berlin (Germany)}
\email{thomas.kraemer@math.hu-berlin.de}
\title[Characteristic cycles and the Gauss map]{Characteristic cycles and the microlocal geometry of the Gauss map, I} 

\begin{document}

\frontmatter



\begin{abstract}
We propose two new approaches to the Tannakian Galois groups of holonomic $\scrD$-modules on abelian varieties. The first is an interpretation in terms of principal bundles given by the Fourier-Mukai transform, which shows that they are almost connected. The second constructs a microlocalization functor relating characteristic cycles to Weyl group orbits of weights. This explains the ubiquity of minuscule representations, and we illustrate it with a Torelli theorem and with a bound for decompositions of a given subvariety as a sum of subvarieties. The appendix sketches a twistor variant that may be useful for $\scrD$-modules not coming from Hodge theory.
\end{abstract}

\begin{altabstract}
Nous proposons deux nouvelles approches aux groupes tannakiens des $\scrD$-modules holonomes sur les vari\'et\'es ab\'eliennes. La premi\`ere est une interpr\'etation en termes de fibrés principaux d\'efinis par la transformation de Fourier-Mukai, ce qui implique qu'ils sont essentiellement connexes. La deuxi\`eme fournit un foncteur de microlocalisation qui relie les cycles caract\'eristiques aux orbites des groupes de Weyl sur les poids. Cela explique l'ubiquit\'e des repr\'esentations minuscules, et nous l'illustrons par un th\'eor\`eme de Torelli et par une borne pour les d\'ecompositions d'une sous-vari\'et\'e donn\'ee comme somme d'autres sous-vari\'et\'es. L'appendice donne une variante twistorielle qui peut \^etre utile pour les $\scrD$-modules ne provenant pas de la th\'eorie de Hodge.
\end{altabstract}

\maketitle
\tableofcontents

\mainmatter

\section{Introduction} \label{sec:results}

A common thread in algebraic geometry is the study of various realizations of the six functor formalism like holonomic $\scrD$-modules, Hodge modules or~$\ell$-adic perverse sheaves. On semiabelian varieties one has a Tannakian correspondence between such objects and representations of certain algebraic groups~\cite{KrWVanishing, KraemerSemiabelian, GaL}. These groups are interesting not only over finite fields where they are related to equidistribution~\cite{KatzSatoTate} and generic vanishing~\cite{WeissauerVanishing}, but also over the complex numbers where they can be used for instance to study singular subvarieties~\cite{KrWSchottky, KraemerThreefolds}.
Our goal in the present paper and its sequel~\cite{KraemerMicrolocalII} is to understand the geometric meaning of the Tannakian formalism in the case of holonomic $\scrD$-modules on complex abelian varieties.

\subsection{Motivation and overview} 
 
Let $A$ be a complex abelian variety. As we will recall in the next section, the Tannakian formalism naturally attaches to any holonomic $\scrD_A$-module $\scrM$ an affine complex algebraic group $G(\scrM)$; this paper is motivated by the following

\begin{ex}
For a proper closed subvariety $Z\subset A$, let $\scrM = \delta_Z$ be the regular holonomic $\scrD_A$-module whose associated perverse sheaf is the perverse intersection complex of the subvariety. There are two extreme cases:
\smallskip
\begin{enumerate} 
\item 
If $Z$ is a finite set of points, one easily sees that $G(\delta_Z)=\Hom(\Gamma, \bbG_m)$ is the Cartier dual of the group $\Gamma \subset A(\bbC)$ generated by these points.\smallskip
\item If $Z$ is a divisor on $A$, then the group $G(\delta_Z)$ is related to the Gauss map sending a smooth point of the divisor to its normal direction. An interesting case is the theta divisor on intermediate Jacobians of cubic threefolds, for which $G(\delta_Z)$ is an exceptional group of type $E_6$ and the monodromy of the Gauss map of the theta divisor can be identified with the Weyl group $W(E_6)$. See~\cite{KraemerThreefolds}.
\end{enumerate}
\end{ex} 

\noindent We generalize both observations to the groups $G=G(\scrM)$ for any $\scrM$.
Our first approach identifies them as structure groups of principal bundles by Schnell's work on the Fourier-Mukai transform~\cite{SchnellHolonomic}. If $G^\circ \subseteq G$ denotes the connected component of the identity, we deduce that
\[
 G/G^\circ \;=\; \Hom(X, \bbG_m)
 \quad  \textnormal{for a finite subgroup} \quad
 X \;\subset\; A(\bbC)
\] 
as in the first example above (theorem~\ref{thm:componentgroup}). This complements the analogous theorem by Weissauer for $\ell$-adic perverse sheaves on abelian varieties over finitely generated fields~\cite[th.~2 and the paragraph after lemma~11]{WeConnected}. The latter works over base fields of arbitrary characteristic, but since it is a result about perverse sheaves, it can in characteristic zero only deal with {\em regular} holonomic $\scrD$-modules. Our proof is very different and works directly with $\scrD$-modules --- as such it is restricted to base fields of characteristic zero but includes the case of {\em irregular} holonomic $\scrD$-modules.

\medskip

Given the above result, we have a fairly good understanding of the group $G/G^\circ$ of connected components, so that the main task will be to understand the connected component $G^\circ \subseteq G$ of the identity for the groups $G=G(\scrM)$. If $\scrM$ is semisimple, then this is a connected reductive group; our second approach studies its root system via subgroups of multiplicative type that are related to Gauss maps. For this we consider the characteristic variety $\Char(\scrM)\subset T^* A$, which is a conic Lagrangian subvariety of the cotangent bundle. For abelian varieties the cotangent bundle is the trivial bundle $T^*A = A \times V$ with fiber $V=H^0(A, \Omega_A^1)$. We define the {\em Gauss map} to be the projection
\[
 \gamma: \quad \Char(\scrM) \;\subset\; T^*A \;=\; A\times V \;\twoheadrightarrow\; V
\]
to the cotangent space. This definition generalizes the classical notion of Gauss maps for divisors on abelian varieties. In the cases relevant to this paper, it is dominant and generically finite. In sections~\ref{sec:germs} and \ref{sec:microlocalization} we construct an embedding 
$
 \Hom(\Gamma, \bbG_m) \hookrightarrow G
$
for the subgroup
\[
 \Gamma \;=\; \langle \, a\in A(\bbC) \mid (a, u) \in \Char(\scrM) \, \rangle
 \;\subset\;
 A(\bbC)
\]
that is generated by the finitely many points in the fiber of the Gauss map over a fixed, very general $u\in V(\bbC)$ (by a {\em very general} point on a variety we mean a point outside countably many proper closed subvarieties). The resulting link between Weyl groups and Gauss maps generalizes the second example from above~(corollary~\ref{cor:weyl}), explains the ubiquity of minuscule representations, and provides a way to compute Tannaka groups that covers all previously known cases in a uniform way (theorem~\ref{thm:smoothsubvarieties}).

\medskip

For geometric applications we refer the reader to~\cite{KraemerMicrolocalII} and only illustrate our ideas with two easy examples: A lower bound for summands of subvarieties (proposition~\ref{prop:sumofsubvarieties}) that shows why for intermediate Jacobians of cubic threefolds the theta divisor is not a sum of curves, a special case of Schreieder's beautiful result~\cite{SchreiederTheta}; and a Torelli theorem which generalizes those for curves and Fano surfaces of cubic threefolds by recovering a smooth subvariety from Tannakian data (corollary~\ref{cor:torelli}). The appendix gives an alternative approach via the theory of twistor modules by Mochizuki and Sabbah that may be useful for holonomic $\scrD$-modules not coming from Hodge modules.

\medskip 

\subsection{The Tannakian framework} \label{sec:tannakiansetting}

Let $\scrD_A$ be the sheaf of algebraic differential operators and 
\[
 \Mhol(\scrD_A) \;\subset\; \Dbhol(\scrD_A)
\]
the abelian category of holonomic right $\scrD_A$-modules resp.~the derived category of bounded algebraic $\scrD_A$-module complexes with holonomic cohomology sheaves. On the latter the addition morphism $a: A\times A \rightarrow A, (x,y)\mapsto x+y$ defines a convolution product
\[
\scrM_1 * \scrM_2 \;=\; a_\dag(\scrM_1 \boxtimes \scrM_2)
\]
which makes the derived category into a {\em tensor category}, i.e.~a $\bbC$-linear symmetric monoidal category. The abelian subcategory $\Mhol(\scrD_A) \subset \Dbhol(\scrD_A)$ is not stable under convolution, but this can be overcome by passing to a quotient category as described in a more general axiomatic framework in~\cite[th.~13.2]{KrWVanishing}. For convenience we briefly recall how the construction works in the case relevant for this paper: For any~$\scrM\in \Mhol(\scrD_A)$ the analytic de~Rham complex \bigskip 
\begin{align*}
 \DR(\scrM) \;=\; \Bigl[ \; \cdots 
 \, \longrightarrow \, \scrM\otimes_{\scrO_A} \wedge^2\scrT_A^\an 
 \, \longrightarrow \,  \scrM \otimes_{\scrO_A} \scrT_A^\an 
 \, \longrightarrow \,  \scrM \otimes_{\scrO_A} \scrO_A^\an & \;   \Bigr]  \\
 %
 \textnormal{degree $-2$} \hspace*{4em} 
 \textnormal{degree $-1$} \hspace*{4em} \textnormal{degree $0$} \hspace*{1em} & \bigskip
\end{align*}
is a perverse sheaf whose hypercohomology has non-negative Euler characteristic by Kashiwara's index formula and the fact that the characteristic cycle of any perverse sheaf is effective~\cite[cor.~1.4]{FKGauss}. In what follows we say that $\scrM$ is~{\em negligible} if the Euler characteristic of its hypercohomology vanishes. Since the Euler characteristic is additive on short exact sequences, the negligible modules form a thick abelian subcategory $\rmS(A)\subset \Mhol(\scrD_A)$. Furthermore, if $\rmT(A)\subset \Dbhol(\scrD_A)$ denotes the thick triangulated subcategory of complexes whose cohomology sheaves are negligible, then by~\cite[prop.~3.6.1]{GaL} the abelian quotient category  \bigskip 
\[ \rmM(A) \;=\; \Mhol(\scrD_A)/\rmS(A) \medskip \]
is naturally equivalent to the heart of the induced $t$-structure on the triangulated Verdier quotient $\rmD(A)=\Dbhol(\scrD_A)/\rmT(A)$ via the following diagram: \bigskip
\[
\xymatrix@M=0.5em{
 \Mhol(\scrD_A) \ar@{->>}[d] \ar@{^{(}->}[r] & \Dbhol(\scrD_A) \ar@{->>}[d]  \\
 \rmM(A) \ar@{^{(}->}[r] & \rmD(A) \medskip
}
\]
From the generic vanishing theorem in~\cite[th.~1.1]{KrWVanishing} and~\cite[th.~4.1]{SchnellHolonomic} one now deduces as in~\cite[\S 5]{KraemerSemiabelian} the following result:

\begin{thm} 
The convolution product descends to the quotient category $\rmD(A)$.~It preserves the full abelian subcategory $\rmM(A)\subset \rmD(A)$ and makes the latter into an abelian tensor category which is neutral Tannakian.
\end{thm}

Recall that $\rmM(A)$ being a {\em neutral Tannakian category} means that it is a rigid abelian tensor category whose unit object $\one$ satisfies $\End(\one)=\bbC$ and which admits a {\em fiber functor}, i.e.~a faithful  exact $\bbC$-linear tensor functor to the category $\Vect(\bbC)$ of finite dimensional complex vector spaces. The existence of such a fiber functor is guaranteed by Deligne's internal characterization of Tannakian categories, since we are working over an algebraically closed field~\cite[\S 6.4]{CoulembierTannakian}. We do not know how to write down explicitly a fiber functor on the entire tensor category  $\rmM(A)$, but on any finitely generated tensor subcategory there are many geometric ways for doing so as we will see below. Here by a {\em finitely generated tensor subcategory} we mean the smallest rigid abelian full tensor subcategory which contains a given $\scrM \in \M(A)$ and is stable under subquotients. We denote this subcategory by~$\langle \scrM \rangle \subset \M(A)$. Any fiber functor~$\omega$ on it induces an equivalence \medskip 
$$
  \langle \scrM \rangle \; \stackrel{\sim}{\longrightarrow} \; \Rep(G) \medskip
$$
with the tensor category of finite dimensional algebraic representations of the affine algebraic group $G=\Aut^\otimes(\omega)$ of its tensor automorphisms. Up to isomorphism this group is determined uniquely by~$\scrM$ because on any Tannakian category over an algebraically closed field any two fiber functors are isomorphic~\cite[th.~3.2(b)]{DM}. So by abuse of notation we often simply write
$G = G(\scrM)$.

\subsection{The Fourier-Mukai transform} \label{sec:fouriermukai-intro}

In~\cite{KrWVanishing} fiber functors have been obtained as the hypercohomology of the twist by generic local systems of rank one. The latter correspond to pairs $(\scrL, \nabla)$ where $\scrL\in \Pic^\circ(A)$ and $\nabla: \scrL \rightarrow \Omega^1_A\otimes \scrL$ is a flat connection. The moduli space of such pairs is a smooth quasiprojective variety~$A^\natural$ that can be viewed as the universal vector extension of the dual abelian variety via the forgetful map~$A^\natural \rightarrow \hat{A} = \Pic^\circ(A)$, see~\cite{MazurMessing} and~\cite[sect.~9]{SchnellHolonomic}. On~$A\times A^\natural$ the pullback of the Poincar\'e bundle has a universal relative flat connection which can be used to define a Fourier-Mukai transform, see section~\ref{subsec:fouriermukai}. It has been shown by Laumon~\cite{LaumonTransformation} and Rothstein~\cite{RothsteinSheaves} that this Fourier-Mukai transform gives an equivalence \bigskip 
\[
\FM: \quad \Dbcoh(\scrD_A) \;\stackrel{\sim}{\longrightarrow} \; \Dbcoh(\scrO_{A^\natural}) \medskip
\]
between the bounded derived categories of complexes of algebraic $\scrD_A$-modules with coherent cohomology sheaves, and complexes of algebraic $\scrO_{A^\natural}$-modules with coherent cohomology sheaves. While $\FM$ is not an exact functor with respect to the standard $t$-structures, it has been shown by Schnell~\cite[th.~18.1]{SchnellHolonomic} that its restriction to the subcategory $\Dbhol(\scrD_A) \subset \Dbcoh(\scrD_A)$ becomes an exact functor when $\Dbcoh(\scrO_{A^\natural})$ is equipped with a certain perverse coherent $t$-structure~\cite{ABPerverse, KashiwaraTStructures}. In particular, for any~$\scrM\in \Mhol(\scrD_A)$ there is an open subset $U\subseteq A^\natural$ with complement of codimension at least two such that \bigskip
\begin{equation} \label{eq:gvt} 
 \scrH^i(\FM(\scrM))|_U
 \;\; \textnormal{is} \;\;
 \begin{cases}
 \; \textnormal{locally free} & \textnormal{for $i=0$}, \\
 \; \textnormal{zero} & \textnormal{for $i\neq 0$}.
 \end{cases} 
\end{equation}
In fact $U \subseteq A^\natural$ can be taken to be the complement of a finite union of translates of proper linear subvarieties, where by a {\em linear subvariety} we mean the image of an embedding $B^\natural \hookrightarrow A^\natural$ given by an epimorphism $A\twoheadrightarrow B$ of abelian varieties.~After removing finitely many further translates of linear subvarieties, we may assume all the subquotients of~$\scrM$ also satisfy~\eqref{eq:gvt}. 
Under this assumption we may lift our Tannakian quotient category to an equivalent subcategory $\langle \scrM \rangle \subset \Mhol(\scrD_A)$, and the fiber functors in~\cite{KrWVanishing} arise as fibers of the Fourier-Mukai transform: Any~$u\in U(\bbC)$ corresponds by definition to a flat connection $(\scrL, \nabla)$ of rank one, and then the fiber of the Fourier-Mukai transform at $u$ is the de Rham cohomology twisted by the local system $L = \ker(\nabla)$: \medskip
\begin{align*}
 \omegaFM{u}: \quad \langle \scrM \rangle \;\longrightarrow \; & \; \Vect(\bbC), \\[0.5em]
 \scrN \;\mapsto\; & \; \scrH^0(\FM(\scrN))(u) \;=\; H^0(A, \DR(\scrM)\otimes L). \medskip
\end{align*}
This leads to the following interpretation:

\begin{thm}[= theorem~\ref{structure_group}] \label{thm:structure_group}
The vector bundle $\scrE = \scrH^0(\FM(\scrM))|_U$ is induced by an algebraic principal bundle whose structure group is isomorphic to the Tannaka group $G(\scrM)=\Aut^\otimes(\omegaFM{u} \!\mid\! \langle \scrM \rangle)$ for any $u\in U(\bbC)$. If $\scrM$ is semisimple, then $G(\scrM)$ is the unique minimal reductive structure group for this vector bundle. 
\end{thm}

Here {\em minimal} means that any other reduction of $\scrE$ to an algebraic principal bundle with a reductive structure group is induced from the given one.~If such a minimal reductive reduction exists, it is unique up to isomorphism; the existence has been shown in~\cite[th.~2.1]{BogomolovStable} on any variety $U$ with~$\dim H^0(U, \scrO)=1$. This last condition is valid in our case, indeed \medskip 
\[
 \dim H^0(U, \scrO_U) \;=\; \dim H^0(A^\natural, \scrO_{A^\natural}) \;=\; 1 \medskip
\]
where the first identity holds since we assumed that the complement of $U\subseteq A^\natural$ is of codimension at least two, and the second is shown in~\cite[th.~2.4.1]{LaumonTransformation}. So for semisimple holonomic $\scrD_A$-modules $\scrM$ the Tannaka group is determined uniquely by the Fourier-Mukai transform on the open dense subset $U\subseteq A^\natural$ as expected from the reconstruction result in~\cite[cor.~21.3]{SchnellHolonomic}. We will use this in section~\ref{sec:almostconnected} for a proof of the following fact that extends a result of Weissauer~\cite{WeConnected} from perverse sheaves to arbitrary holonomic $\scrD_A$-modules:

\begin{thm} [= theorem~\ref{componentgroup}]\label{thm:componentgroup}
In the above setting, the group of connected components of the Tannaka group $G=G(\scrM)$ is naturally a quotient
$
 \pi_1(\hat{A}, 0) \twoheadrightarrow  G/G^\circ
$
of the fundamental group of the dual abelian variety $\hat{A}$.
\end{thm}

In particular, it follows that $G/G^\circ$ is a finite abelian group, so any object of the full subcategory $\Rep(G/G^\circ) \subseteq \Rep(G)$ is a sum of characters $\chi_a: \pi_1(\hat{A}, 0) \rightarrow G/G^\circ \rightarrow \bbC^*$ of finite order. The corresponding flat line bundles are parametrized by certain torsion points $a \in A=\Pic^\circ(\hat{A})$, and if $\delta_a$ denotes the Dirac module supported on such a point, the associated line bundle is the Fourier-Mukai transform $\FM(\delta_a)$. So the group of components is the Cartier dual 
\[
 G/G^\circ \;=\; \CartierDual{K}
\]
of the finite abelian group
\[
 K \;=\; 
\{ \, a\in A(\bbC) \;\textnormal{torsion point} \; | \; \delta_a \in \langle \scrM \rangle
\,  \}.
\]
For reference we include in section~\ref{sec:almostconnected} a proof of the following result:

\begin{cor} \label{cor:isogeny} \label{cor:restriction-irreducible}
Let $G=G(\scrM)$ and define the finite group $K\subset A(\bbC)$ as above. 
\begin{enumerate} 
\item If $p: A\twoheadrightarrow B$ is an isogeny of abelian varieties, then $\scrN = p_\dag(\scrM)$ has as its Tannaka group a subgroup 
$G(\scrN) \subseteq G$
of finite index. This subgroup is given by
$$ \quad \qquad G(\scrN)/G^\circ \;=\; \CartierDual{K/K\cap \ker(p)} \;\subseteq\; G/G^\circ \;=\; \CartierDual{K}. $$
\item Let $\scrM$ be simple, corresponding to an irreducible representation $W\in \Rep(G)$ of the Tannaka group.  Then the restriction $W|_{G^\circ}\in \Rep(G^\circ)$ is irreducible if and only if
$$ 
 \quad \qquad t_x^*(\scrM) \; \not\simeq \; \scrM 
 \quad \textnormal{\em for the translations $t_x: A\rightarrow A$ by all $x\in K\setminus \{0\}$}.
$$
\end{enumerate}
\end{cor}

\noindent
It remains to understand the connected component of the identity of the Tannaka groups. This leads us to the main goal of this paper: To explain the relation between Weyl groups and the monodromy of Gauss maps observed in~\cite{KraemerThreefolds}.

\subsection{Characteristic cycles and Gauss maps}

For the relation with Gauss maps we will use results from microlocal analysis as developed by Sato, Kashiwara and Kawai~\cite{SKK, KashiwaraSystems}. Here {\em microlocal} means that we study differential operators on a complex manifold by working locally on its cotangent bundle. To any differential operator one may attach its principal symbol, which is a function on the cotangent bundle; motivated by the theory of elliptic differential operators in real analysis, we want a ``local inverse'' for a differential operator at each point of the cotangent bundle where its principal symbol is non-zero. In the theory by Sato, Kashiwara and Kawai this is achieved by replacing the sheaf of differential operators by a larger sheaf of microdifferential operators on the cotangent bundle, see section~\ref{sec:microdifferential}. Generalizing the vanishing locus of principal symbols, one may attach to any coherent $\scrD$-module its characteristic variety. The latter can be defined algebraically as the support of the associated graded for any local good filtration; this allows to attach multiplicities to its components, thus defining the {\em characteristic cycle}~\cite[sect.~2.2]{HTT}. For $\scrM \in \Mhol(\scrD_A)$ this cycle has the form \medskip
\[
 \CC(\scrM) \;=\; \sum_{\Lambda} \, m_\Lambda(\scrM) \cdot \Lambda
 \quad \textnormal{with} \quad
 m_\Lambda(\scrM) \;\in\; \bbN_0 \;=\; \bbN \cup \{0\},
\]
where $\Lambda \subset T^*A$ runs over the irreducible closed conic Lagrangian subvarieties of the cotangent bundle. Any such $\Lambda$ is the closure of the conormal bundle to the smooth locus of some closed subvariety of $A$. Since $T^* A = A\times V$ is the trivial bundle with fiber $V=H^0(A, \Omega_A^1)$, we may generalize the usual Gauss map for a smooth divisor by considering the projection
\[
 \gamma_\Lambda: \quad  \Lambda \;\hookrightarrow\; A\times V \;\twoheadrightarrow\;  V
\]
to the cotangent space. This {\em Gauss map} is generically finite and we denote by $d_\Lambda \in \bbN_0$ its generic degree. If $\gamma_\Lambda$ is not dominant, we put $d_\Lambda = 0$ and say $\Lambda$ is {\em degenerate}. In both cases the degree $d_\Lambda$ is the cardinality of the fiber of the Gauss map over a general point $u\in V(\bbC)$; to avoid selecting such a point, we will later also view~$V$ as a scheme and work over its generic point. By Kashiwara's index formula~\cite[th.~1.3 and prop.~2.2]{FKGauss} we have

\[
 \sum_{i\in \bbZ} \; (-1)^i \dim_\bbC H^i(A, \DR(\scrM)) \;=\; \sum_{\Lambda} \, m_\Lambda(\scrM) \cdot d_\Lambda.
\]
So $\scrM$ is negligible if and only if all components of $\CC(\scrM)$ are degenerate. To explain the link between Tannaka groups and conormal geometry, we will categorify the index formula by naturally attaching to any holonomic $\scrD_A$-module $\scrM$ a local system of rank $m_\Lambda(\scrM)$ on an open dense subset of each non-degenerate irreducible component $\Lambda$ of the characteristic variety $\Char(\scrM) = \Supp(\CC(\scrM)) \subset T^*A$. 

\medskip 

There is a classical such construction, the {\em second microlocalization}~\cite[sect.~6]{GinzburgCharacteristic}, but it results in {\em twisted local systems}: For a holonomic $\scrD_A$-module whose characteristic variety $\Lambda = \Lambda_Y \subset T^* A$ is the conormal variety to a smooth subvariety $Y\subset A$, its second microlocalization will be a holonomic $\scrD_{\Lambda^{1/2}}$-module. Here we denote by~$\scrD_{\Lambda^{1/2}}$ the sheaf of {\em twisted differential operators} acting on a local square root of the relative canonical bundle $\omega_{\Lambda/Y}$ \cite[rem.~2.6.5]{KashiwaraRepresentation}; by loc.~cit.~this sheaf of twisted differential operators is globally well-defined even if the square root of the line bundle exists only locally. In order to see a monodromy operation we must remove the twist. For this it seems natural to construct a square root of the line bundle $\omega_{\Lambda/Y}$ away from the branch locus of the Gauss map $\gamma: \Lambda \to V$ by pulling back a square root of the trivial line bundle $\omega_V \simeq \scrO_V$. However, for the Tannakian formalism we need to make sure that our chosen way of removing the twists is compatible with the convolution product. To achieve this we follow a slightly different route: We pass to microdifferential modules on the cotangent bundle (sections~\ref{sec:microdifferential} and~\ref{sec:microconvolution}) and compare these with a carefully chosen class of simple microdifferential modules on abelian varieties that behaves well under convolution (sections~\ref{sec:simplemodules} and~\ref{sec:secondmicro}). Each of our simple modules will only be defined on the locus where the corresponding Gauss map is finite \'etale. Since this locus cannot be fixed in advance, we introduce  in section~\ref{sec:germs} a tensor category $\LS(A, \eta)$ of germs of local systems near the generic point $\eta \in V$ to get a tensor functor 
$F: \M(A) \longrightarrow \LS(A, \eta)$.
For $\scrM\in \Mhol(\scrD_A)$, consider then the Gauss map 
\[
 \gamma: \quad \Char(\scrM) \; \subset \; A\times V
 \; \twoheadrightarrow \; V
\]
on the characteristic variety. Over some Zariski open dense subset this is a finite \'etale cover; we define its monodromy group to be the automorphism group $\Gal(\gamma)$ of its Galois closure. For very general $u\in V(\bbC)$, i.e.~all $u$ outside countably many proper closed subvarieties, we will show that the above tensor functor $F$ factors on the subcategory $\langle \scrM \rangle \subseteq \M(A)$ over a tensor functor $F_u: \langle \scrM \rangle \longrightarrow \LS(A, u)$ to an analogous category of Zariski germs of local systems over the closed point $u$. This will then induce a new fiber functor $\omegaMicro{u}$ as in the following diagram: \bigskip
\[
\xymatrix@C=1.5em@M=0.5em{
 \langle \scrM \rangle \ar@{..>}[d]_-{\exists F_u} \ar[rr]^-{\omegaMicro{u}} && \Vect(\bbC) \\
 \LS(A, u) \ar[rr] && \Vect_{X_u}(\bbC) \ar[u] 
}
\]
Here $\Vect_{X_u}(\bbC)$ is the category of finite dimensional vector spaces with a grading by the group 
\[
 X_u \;=\; X_u(\scrM) \;=\; \bigl\langle a \in A(\bbC) \mid (a,u) \in \Char(\scrM) \bigr\rangle  \;\subset\; A(\bbC)
\]
generated by the finitely many points in the fiber of the Gauss map. Notice that even though any two fiber functors on a Tannakian category over an algebraically closed fields are non-canonically isomorphic, there is no obvious geometric way to relate the new fiber functors coming from microlocal analysis to those that we found earlier via the Fourier-Mukai transform: The fiber functors $\omegaFM{u}$ are defined for a general point~$u\in A^\natural(\bbC)$ and do not seem directly related to Gauss maps or characteristic cycles, whereas the fiber functors $\omegaMicro{u}$ are defined only for {\em very} general $u\in V(\bbC)$ but make the relation with Gauss maps transparent. Since all fiber functors are non-canonically isomorphic and we mostly care about the Tannaka groups only up to isomorphism, we now take the group \bigskip
\[
 G_u \;=\; G_u(\scrM) \;=\; \Aut^\otimes(\omegaMicro{u} \!\mid \! \langle \scrM \rangle) \medskip
\] 
of tensor automorphisms of the above fiber functor as a specific realization of the Tannaka groups. Then the factorization of $\omegaMicro{u}$ over a fiber functor to $\Vect_{X_u}(\bbC)$ leads to the following relation between Tannaka groups and Gauss maps; here for a subgroup $T\hookrightarrow G$ of multiplicative type (see section~\ref{sec:multiplicativeI}) in an algebraic group $G$ we denote by 
\[ W(G, T) = N_G(T)/Z_G(T) \medskip \]
the quotient of its normalizer by its centralizer, generalizing the notion of Weyl groups: \pagebreak

\begin{thm}
\label{thm:gausstorus} 
Fix $\scrM\in \Mhol(\scrD_A)$ and a very general point $u\in V(\bbC)$.
\begin{enumerate} 
\item 
With notation as above, we have embeddings
\begin{eqnarray*}
	 T_u \; :=\; \CartierDual{X_u} &\;\hookrightarrow \;& G_u \;=\; G_u(\scrM), 
	 \\
	 \Gal(\gamma) &\;\hookrightarrow \;& W_u \;=\; W(G_u, T_u).
\end{eqnarray*} 
\item Each point $(a,u)\in \gamma^{-1}(u) \subset \Char(\scrM)$ defines a character $\chi_a\in \CartierDual{T_u}$, and 
\[
 \omega_u(\scrM)|_{T_u} \;\in\; \Rep(T_u)
\]
splits as a direct sum of such characters. In this decomposition, the multiplicity of~$\chi_a$ coincides with the multiplicity $m_\Lambda(\scrM)$ of the characteristic cycle $\CC(\scrM)$ along the unique irreducible component $\Lambda \subseteq \Char(\scrM)$ with $(a, u)\in \Lambda$. 
\end{enumerate} 
\end{thm}

\noindent
Part (a) is an direct consequence of theorem~\ref{thm:microlocal-fiberfunctor}, theorem~\ref{thm:gradedfiberfunctor} and lemma~\ref{lem:specialization}.
In fact the monodromy operation gives rise to a homomorphism $\Gal(\gamma) \to \Aut(T_u)$, and the proof of theorem~\ref{thm:gradedfiberfunctor} will show that this monodromy action is induced by the conjugation action~$W_u \rightarrow \Aut(T_u)$. For part (b), note that any representation of a group of multiplicative type splits as a direct sum of characters; so it only remains to observe that the dimensions of the graded pieces of $\omega_u(\scrM) \in \Vect_{X_u}(\bbC)$ match the multiplicities $m_\Lambda(\scrM)$ of the characteristic cycle, which will be clear from our construction of the fiber functor in theorem~\ref{thm:microlocal-fiberfunctor}.

\subsection{Weight spaces and minuscule representations} \label{sec:minuscule}

In many   applications the subgroups of multiplicative type that we obtained in theorem~\ref{thm:gausstorus} are quite big. Let us say that $\scrM \in \Mhol(\scrD_A)$ is {\em generically multiplicity-free} if any non-degenerate component of its characteristic cycle has multiplicity one and if the same condition is also satisfied for the pushforward under any isogeny of abelian varieties --- in other words, if all non-degenerate components $\Lambda_{Z_1}, \Lambda_{Z_2}\subset \Char(\scrM)$ enter with multiplicity one and satisfy 
\[
 Z_1 \;\neq \; Z_2 + x \quad \textnormal{for all torsion points} \quad x\in A(\bbC) \setminus \{0\}.
\]
If $Y\subset A$ is a smooth irreducible subvariety which is not stable under any non-trivial translation, then $\delta_Y$ is generically multiplicity-free; this in particular applies to any smooth summand of a divisor defining a principal polarization, such as a smooth curve in its Jacobian, the Fano surface of lines in the intermediate Jacobian of a smooth cubic threefold, or any smooth theta divisor. If $G$ is a reductive group, let
$W(G) = W(G^\circ, T)$
denote the Weyl group for a given maximal torus~$T\subseteq G^\circ$ inside the connected component of the identity. For any representation $U\in \Rep(G)$ the restriction $U|_T \in \Rep(T)$ splits into a sum of one-dimensional representations which we call the {\em weights} of the representation, and these weights are permuted by the Weyl group. We then obtain the following answer to the question raised in~\cite[conj.~8]{KraemerThreefolds}:

\begin{cor} \label{cor:weyl}
If $\scrM$ is semisimple and generically multiplicity-free, theorem~\ref{thm:gausstorus} gives  
\begin{enumerate} 
\item an embedding
$\Gal(\gamma) \;\hookrightarrow \; W(G_u)$, and \smallskip
%
\item a bijection between $\gamma^{-1}(u)\subset \Char(\scrM)$ and the weights in $\omegaMicro{u}(\scrM) \in \Rep(G_u)$ so that the monodromy and Weyl group actions match.
\end{enumerate}
\end{cor}

{\em Proof.} By definition the direct image of a generically multiplicity free module under an isogeny remains generically multiplicity free. So by corollary~\ref{cor:isogeny} we may assume that the group $G=G_u$ is connected and the subgroup $X=X_u \subset A(\bbC)$ is torsion-free, which means that its Cartier dual $T=T_u$ is a torus. The generic multiplicity-freeness ensures that even after pushforward under an isogeny all non-degenerate components of $\CC(\scrM)$ will be reduced; so theorem~\ref{thm:gausstorus}(b) says that the restriction of the faithful representation $W=\omegaMicro{u}(\scrM)\in \Rep(G)$ to the torus~$T$ splits as a sum of pairwise distinct characters. The centralizer $Z_G(T)$ preserves this decomposition and hence it is contained in a subgroup of diagonal matrices in~$\Gl(W)$. Its connected component of the identity $Z_G(T)^\circ$ is therefore a subtorus of $G$ and hence a maximal torus, since any other subtorus containing $T$ must also lie in the centralizer. The image of the fundamental group under the monodromy representation normalizes this maximal torus since
$N_G(T) \subseteq N_G(Z_G(T)) \subseteq  N_G(Z_G(T)^\circ )$. \qed

\begin{rem} 
The bijection in part (b) of the corollary is only for the weights of one specific representation; in general the subgroup of multiplicative type $T_u \subset G_u$ is {\em not} a maximal torus: If $A=E$ is an elliptic curve, then by minimal extension of suitable branched local systems of rank $r\geq 2$ one finds cases where $G_u(\scrM)\simeq \Sl_3(\bbC)$ and $\CC(\scrM)=\Lambda_{\{0\}} + \Lambda_{\{p\}} + \Lambda_{\{-p\}} + r\cdot \Lambda_E$ for some non-torsion point $p\in E(\bbC)\setminus \{0\}$; then
\[
 T_u \;:=\; \CartierDual{X_u} \;\simeq\; \bbG_m 
\]
is a torus of rank one and hence not a maximal torus inside the group $\Sl_3(\bbC)$. Note that in this case the restriction of the $3$-dimensional standard representation still splits into three pairwise distinct characters of the subtorus $T_u \subset G_u$, so the latter still distinguishes the weights of the standard representation from each other as claimed in (b). In this example the monodromy action is trivial because the Gauss map is a trivial disconnected cover; we will see in theorem~\ref{thm:smoothsubvarieties} that when the monodromy is big, then in certain cases we do get a maximal torus from theorem~\ref{thm:gausstorus}.
\end{rem} 

\medskip

Recall that an irreducible representation is {\em minuscule} if its weights form a single Weyl group orbit; minuscule representations have very special properties, a detailed discussion of their combinatorics can be found in~\cite[chapter~5]{GreenMinuscule}. For the simple Lie algebras the only non-trivial minuscule representations are the fundamental ones of the following dimensions:  \medskip
{\small
\[
\def\arraystretch{1.4}
\begin{array}{|r|r|r|r|r|r|} \hline
  A_n & B_n & C_n & D_n & E_6 & E_7 \\ \hline 
 {n+1 \choose k} \; \textnormal{for} \; 1\leq k \leq n & 2^n & 2n  & 2n, \; 2^{n-1} & 27 & 56 \\ \hline 
\end{array}
\]
}

\noindent
In particular there are no minuscule representations for the types $E_8$, $F_4$, $G_2$, so we may exclude these types in many cases:

\begin{cor} \label{cor:minuscule}
Let $\Lambda \subset T^*A$ be an irreducible Lagrangian subvariety not stable under any translation by a point of $A(\bbC)$.
If $\scrM \in \Mhol(A)$ has the characteristic cycle $\CC(\scrM) \equiv \Lambda$ modulo a linear combination of degenerate subvarieties, then the corresponding representation of $G_u(\scrM)^\circ \subset G_u(\scrM)$ is minuscule.
\end{cor}

{\em Proof.} The irreducibility of $\Lambda$ says that the monodromy group acts transitively on the fiber of the Gauss map, hence the Weyl group of the connected component of the identity acts transitively on the weights in the corresponding representation. \qed

\begin{rem}
Our proof of theorem~\ref{thm:gausstorus} uses the analytic microlocalization by Sato, Kashiwara and Kawai to categorify characteristic cycles of $\scrD$-modules. More recently Beilinson has introduced a notion of characteristic varieties also for $\ell$-adic perverse sheaves over arbitrary base fields~\cite{BeilinsonHolonomic}, and Saito has extended it to a theory of characteristic cycles and an analog of Kashiwara's index formula~\cite{SaitoCC}. We do not know if our results have an~$\ell$-adic analog, using some categorification parallel to the analytic theory of microdifferential operators. If such a theory exists, then in positive characteristic the link with conormal geometry will be weaker due to Deligne's examples of non-Lagrangian characteristic cycles~\cite[sect.~1.3]{BeilinsonHolonomic}.
\end{rem}

\section{Geometric applications} \label{sec:applications}

Before we come to the proof of our main results, let us illustrate them with a few simple applications. For further developments we refer the reader to~\cite{KraemerMicrolocalII}.

\subsection{Examples with big monodromy} 

The results described above provide an effective way to determine Tannaka groups if we have upper bounds on these groups and lower bounds on the monodromy of Gauss maps. Upper bounds come from invariants in tensor powers. For $\scrM \in \Mhol(\scrD_A)$, let~$\gamma: \Char(\scrM) \to V$ be the Gauss map and $d=\deg(\gamma)$ its generic degree, which coincides with the dimension of the representation corresponding to $\scrM$ under the Tannakian formalism. Translating back from multilinear algebra, we say that \smallskip 
\begin{enumerate} 
\item $\scrM$ has {\em trivial determinant} if its top exterior convolution power $\mathrm{Alt}^{*d}(\scrM)$ is isomorphic to the unit object $\one = \delta_0$ modulo negligible summands. This always holds up to a translation by a point, see lemma~\ref{lem:character}. \smallskip
\item $\scrM$ is~{\em self-dual} if
it is isomorphic to its dual $\scrM^\vee = (-\id_A)^* \, \rmD(\scrM)$, where we denote by
\[
\qquad \rmD: \;\, \Mhol(\scrD_A) \,\stackrel{\sim}{\longrightarrow}\, \Mhol(\scrD_A), \;\, \scrM \;\mapsto\; \scrExt^{\hspace{0.1em}g}_{\scrD_A}(\scrM, \omega_A\otimes \scrD_A) 
\;\; \textnormal{for} \;\; g = \dim A
\]
the usual duality functor. In this case the characteristic cycle $\CC(\scrM)$ is stable under the involution $(-\id_A) \times \id_V$. If this characteristic cycle is moreover reduced and has no component supported over a point $p\in A(\bbC)$ with $2p = 0$, then the Gauss map $\gamma: \Char(\scrM)\to V$ will have as its general fiber a collection of $d/2$ distinct pairs of mutually inverse points, in which case the monodromy operation is given by a subgroup $\Gal(\gamma) \hookrightarrow H_d = (\pm 1)^{d/2} \rtimes \frakS_{d/2}$.
\item \smallskip $\scrM$ admits a {\em cubic form} if there is a non-trivial morphism $\scrM^{*3} \rightarrow \one$. \smallskip
\end{enumerate}
Concerning lower bounds on monodromy, we say that a subgroup of the symmetric group $\frakS_d$ is {\em irreducible} if it acts irreducibly on the natural $(d-1)$-dimensional complex representation which is the quotient of the standard permutation representation by its one-dimensional subspace of invariants. If $d$ is even, we say that a subgroup of~$H_d = \{\pm 1\}^{d/2} \rtimes \frakS_{d/2}$ is {\em irreducible} if it acts irreducibly on the natural complex representation of dimension $d/2$ in which  $\{\pm 1\}^{d/2}$ acts via diagonal matrices and~$\frakS_{d/2}$ acts by permuting the factors.

\medskip

\begin{thm} \label{thm:smoothsubvarieties}
Let $G=G_u(\scrM)$ where $\scrM\in \Mhol(\scrD_A)$ has trivial determinant, and assume that modulo degenerate components the characteristic cycle $\CC(\scrM)$ is irreducible, reduced and not stable under translation by any $a\in A(\bbC) \setminus \{0\}$. Denote by $\gamma: \Char(\scrM) \longrightarrow V$ its Gauss map and put $d=\deg(\gamma)$. 
\smallskip%
\begin{enumerate}
\item If the monodromy group $\Gal(\gamma)$ is irreducible in $\frakS_d$, then $G\simeq \Sl_d(\bbC)$. \medskip 
\item If $\scrM$ is self-dual with $\Gal(\gamma)$ irreducible in $H_d$, then $G\in\{ \SO_d(\bbC),\Sp_d(\bbC)\}$. 
\item \medskip If $\scrM$ admits a cubic form, $\Gal(\gamma) \simeq W(E_6)$ and $d=27$, then $G\simeq E_6(\bbC)$.
\end{enumerate}
\smallskip In all three cases the subgroup $T_u \subset G$ is a maximal torus for very general $u$.
\end{thm}

{\em Proof.} The triviality of the determinant implies $G\hookrightarrow \Sl_d(\bbC)$, and in the self-dual case $G\hookrightarrow \SO_d(\bbC)$ or $G\hookrightarrow\Sp_d(\bbC)$ since every irreducible self-dual representation is orthogonal or symplectic. By assumption~$\CC(\scrM)$ is irreducible modulo degenerate subvarieties, hence $\scrM$ corresponds to a minuscule representation $W\in \Rep(G)$ by corollary~\ref{cor:minuscule}. In particular, the representation $W$ is irreducible and therefore the group $G$ is reductive. By Schur's lemma its center acts on $W$ by scalars, and since by assumption the determinant $\det(W)$ is the trivial representation, it follows that the center is finite. Hence $G$ is a semisimple group, and its Lie algebra $\mathrm{Lie}(G)$ splits as a product of simple Lie algebras such that $W$ splits as a tensor product of minuscule representations from the table in section~\ref{sec:minuscule}. Since $\dim(W) = d$, one deduces that in case 3 one of the following four options must occur:

{\small
\[
	\def\arraystretch{1.4}
	\begin{array}{|r||r|r|r|r|r|} \hline
	\mathrm{Lie}(G) & \mathfrak{sl}_{27}(\bbC) &  \mathfrak{sl}_3(\bbC) \times \mathfrak{sl}_9(\bbC) &  \mathfrak{sl}_3(\bbC) \times \mathfrak{sl}_3(\bbC) \times \mathfrak{sl}_3(\bbC) & \mathfrak{e}_6(\bbC) \\ \hline 
	W & \bbC^{27} & \bbC^3\otimes \bbC^9  & \bbC^3\otimes \bbC^3\otimes \bbC^3 & \bbC^{27} \\ \hline 
	\end{array}
\]
}

\noindent Here the entries of the lower row refer to the natural representations of the respective dimensions. Since $\mathfrak{sl}_n(\bbC)$ does not preserve a cubic form on its natural $n$-dimensional representation, the first three options from the above list are excluded; the only remaining option is $\mathrm{Lie}(G)\simeq \mathfrak{e}_6(\bbC)$, acting on $W$ via one of its two $27$-dimensional fundamental representations. Each of the latter admits a unique invariant cubic form, and since the simply connected group $E_6(\bbC)\subset \Gl_{27}(\bbC)$ is the full stabilizer of this cubic form~\cite[thm.~7.3.2]{SpringerVeldkamp}, it follows that in case 3 the Tannaka group is $G\simeq E_6(\bbC)$.

\medskip 

 It remains to show that for very general~$u$ the rank of the finitely generated abelian group $X=X_u$ is at least $d-1$ in case~1, at least $d/2$ in case 2 and at least six in case~3. This rank is the dimension of the representation
$X_\bbC = X \otimes_\bbZ \bbC \in \Rep(\Gal(\gamma))$,
and the latter is faithful because by generic multiplicity-freeness no two points in a general fiber of $\gamma$ differ by a torsion point. Hence in cases~1 and~2 the claim follows since this representation is a quotient of the natural representation of dimension $d-1$ resp.~$d/2$, on which $\Gal(\gamma)$ acts irreducibly by assumption. For case 3 recall that the smallest dimension of a non-trivial representation of $W(E_6)$ is six, see~\cite[\S 13.2, p.~415]{CarterFinite} or the character table in~\cite[p.~27]{ATLAS}. \qed

\medskip

\begin{ex} \label{ex:jacobians}
Theorem~\ref{thm:smoothsubvarieties} provides a uniform approach to all the examples of intersection cohomology sheaves discussed in~\cite[th.~9]{KraemerThreefolds}:

\smallskip 

(a) If $C\subset A$ is a smooth projective curve of genus $g > 1$, embedded in its Jacobian variety via a suitable translate of the Abel-Jacobi map, then it is shown in~\cite[th.~14]{WeBN} and~\cite[th.~6.1]{KrWSmall} that 
\[
 G(\delta_C) \;\simeq\; 
 \begin{cases} 
 \Sp_{2g-2}(\bbC) & \textnormal{if $C$ is hyperelliptic}, \\
 \, \Sl_{2g-2}(\bbC) & \textnormal{if $C$ is not hyperelliptic},
 \end{cases}
\]
and $\delta_C$ corresponds to the natural representation of dimension $2g-2$. 
Note that for hyperelliptic curves we can embed $C=-C\subset A$ as a symmetric subvariety in its Jacobian, with the inversion map inducing the hyperelliptic involution; this gives a symplectic pairing on the corresponding representation. But for non-hyperelliptic curves we always have $C\neq -C$, and the choice between the curve and its negative amounts to the choice between the natural representation and its dual. There is no preferred such choice: We defined $G(\delta_C)$ only up to non-canonical isomorphism and the two representations are conjugate under the outer automorphism of $\Sl_{2g-2}(\bbC)$.

\smallskip 

(b) For the theta divisor $\Theta \subset A$ of a generic principally polarized abelian variety of dimension $g$, it has been shown in~\cite[th.~1.4]{KrWSchottky} that for a suitable translate of the divisor one has
\[
 G(\delta_\Theta) \;\simeq\;
 \begin{cases}
 \; \Sp_{g!}(\bbC) & \textnormal{if $2\mid g$}, \\
 \SO_{g!}(\bbC) & \textnormal{if $2\nmid g$},
 \end{cases}
\]
Again the existence of a bilinear form in both cases comes from the fact that one can choose the translate of the theta divisor to be symmetric. 

\smallskip

(c)  For the Fano surface $S\subset A$ of lines on a smooth cubic threefold, embedded in its intermediate Jacobian via the Albanese morphism, it has been shown in~\cite[th.~1]{KraemerThreefolds} by an ad hoc argument that
\[
 G(\delta_S) \;\simeq\; E_6(\bbC).
\]
We can now explain this geometrically: Here the monodromy group of the Gauss map is known to be the Weyl group $W(E_6)$, acting in its permutation representation on the $27$ lines on a smooth hyperplane section of the cubic threefold, so to apply theorem~\ref{thm:smoothsubvarieties} we only need to show  there is a non-trivial morphism
\[
 \delta_S * \delta_S * \delta_S \;\longrightarrow\; \one.
\]
But this easily follows from the fact that the addition morphism $a: S\times S\times S \rightarrow A$ collapses the incidence variety of coplanar triples of lines on the threefold to a single point~\cite[eq.~(11.9)]{CGIntermediateJacobian}, which we may assume to be the origin. 
We can now also see why for the theta divisor $\Theta = S - S \subset A$ the characteristic cycle~$\CC(\delta_\Theta)$ contains the component $\Lambda_{\{0\}}$ with multiplicity six as noted in~\cite[rem.~7b]{KraemerThreefolds}: We know from loc.~cit.~that in this case the object $\delta_\Theta$ corresponds to the adjoint representation of the Tannaka group $G(\delta_S)$ on its Lie algebra, so the multiplicity of the zero weight is the rank of a maximal torus inside this Tannaka group.
\end{ex}

\subsection{Summands of subvarieties} \label{sec:sumofsubvarieties}

Given a subvariety $W\subset A$, fix an arbitrary fiber functor 
\[
 \omega: \quad \langle \delta_W \rangle \; \stackrel{\sim}{\longrightarrow} \; \Rep(G(\delta_W))
\] 
and let $\Ad_W\in \Mhol(\scrD_A)$ be the unique semisimple holonomic $\scrD_A$-module without negligible summands such that~$\omega(\Ad_W)$ is the adjoint representation of the reductive group $G(\delta_W)$ on its Lie algebra. Note that $\Ad_W$ does not depend on the chosen fiber functor; there is an intrinsic way to define the adjoint object in any finitely generated Tannakian category~\cite[def.~2.4]{CommelinMTC}. In order to have $\dim \omega(\delta_W) > 0$ we want the Gauss map to be dominant, so we will always assume $W\subset A$ to be irreducible with finite stabilizer
\[ \Stab(W) \;=\; \{a\in A(\bbC) \mid W + a = W\} \] 
as in~\cite{WeissauerDegenerate}. We then have the following obstruction for the existence of non-trivial summands of the given subvariety; see also~\cite[th.~2]{KraemerMicrolocalII} where this is  strengthened to an estimate which holds not just for one but for all the summands in a hypothetical sum decomposition:

\begin{propo} \label{prop:sumofsubvarieties}
Suppose that $W\subset A$ is a closed subvariety which decomposes as a sum of subvarieties
\[
 W=Y_1+\cdots + Y_n 
 \quad \textnormal{\em with} \quad 
 \dim(W)=\dim(Y_1)+\cdots+\dim(Y_n).
\]
If $\Stab(W)$ is finite and the Lie algebra of the Tannaka group~$G(\delta_W)$ is simple, then we have\smallskip
\[
 \dim(Y_i) \; \geq \; \tfrac{1}{2} \dim(\Supp(\Ad_W)) 
 \quad \textnormal{\em for some} \quad
 i\;\in\; \{1,2,\dots, n\}.
\] 
\end{propo}

{\em Proof.} The addition morphism $f: Z= Y_1\times \cdots \times Y_n \twoheadrightarrow W$ is a surjective morphism between varieties of the same dimension, so it restricts over some smooth open dense subset $W_0\subseteq W$ of the target to a finite \'etale cover $f_0: Z_0 = f^{-1}(W_0) \twoheadrightarrow W_0$ of smooth varieties. The trace morphism for this finite \'etale cover gives an embedding as a direct summand
\[ 
\delta_{W_0} 
\; \hookrightarrow \; 
f_{0\dag}(\delta_{Z_0})
\; = \;
f_\dag(\delta_Z)|_{W_0}
\;=\;
 (\delta_{Y_1}*\cdots * \delta_{Y_n})|_{W_0}.
\]
By the decomposition theorem for the corresponding perverse sheaves this summand extends to a direct summand $\delta_W \hookrightarrow \delta_{Y_1} * \cdots * \delta_{Y_n}$ in the derived category of regular holonomic $\scrD_A$-modules. In particular, passing to Tannakian quotient categories we have
\[
 \delta_W 
 \; \in \; \langle \delta_{Y} \rangle \;=\; \langle \delta_{Y_1} \oplus \cdots \oplus \delta_{Y_n} \rangle
 \quad \textnormal{where} \quad 
 Y \;=\; Y_1 \cup \cdots \cup Y_n.
\]
By~\cite[prop.~2.21]{DM} this defines an epimorphism
$p: G(\delta_{Y}) \twoheadrightarrow  G(\delta_W)$. Taking direct images under an isogeny of abelian varieties we may assume by corollary~\ref{cor:isogeny} that the reductive groups $G(\delta_Y)$ and $G(\delta_W)$ are connected. In fact we can moreover assume that these groups are semisimple: Indeed, replacing the subvarieties $Y_i \subset A$ by translates we can assume that $\delta_{Y_i}$ has trivial determinant for $i=1,\dots, n$, see lemma~\ref{lem:character}. Since by Schur's lemma the center of $G(\delta_{Y_i})$ acts on the corresponding irreducible representation by a scalar, the triviality of the determinant forces $G(\delta_{Y_i})$ to be semisimple. Of course the chosen translations will also result in a translation of the sum $W=Y_1 + \cdots + Y_n \subseteq A$, but by~\cite[lemma~4.3.2]{KraemerMicrolocalII} such a translation does not affect the derived group of the connected component of the group $G(\delta_W)$, hence it does not change the dimension $\dim(\Supp(\Ad_W))$. The Tannakian formalism gives epimorphisms 
\[ G(\delta_{Y_1}) \times \cdots \times G(\delta_{Y_n})
\; \twoheadrightarrow \;
 G(\delta_Y)
 \; \twoheadrightarrow \;
 G(\delta_W), \]
hence if the groups on the left hand side are semisimple, then so are the groups $G(\delta_Y)$ and $G(\delta_W)$. In what follows we will always make this assumption.

\medskip

Let $\tilde{G}(\delta_Y) \twoheadrightarrow G(\delta_Y)$ and $\tilde{G}(\delta_W) \twoheadrightarrow G(\delta_W)$ be the universal cover of the above connected semisimple groups. Then by the correspondence between simply connected groups and their Lie algebras the epimorphism $p$ lifts to an epimorphism $\tilde{p}$ with a section $\iota$ as indicated in the diagram  
\[
\xymatrix@M=0.5em@R=4em@C=4em{
 \tilde{G}(\delta_Y)  \ar@{->>}[d]^-{\,\tilde{p}} \ar@{->>}[r] & G(\delta_Y) \ar@{->>}[d]^-p  \\
 \tilde{G}(\delta_W) \ar@{->>}[r] \ar@/^2pc/@{..>}[u]^-{\exists\, \iota} & G(\delta_W) 
}
\]
such that the image of the section $\iota$ commutes with the kernel~$K=\ker(\tilde{p})$. So we get an isomorphism
\[ 
 G \;=\; K\times \tilde{G}(\delta_W)\; \stackrel{\sim}{\longrightarrow} \; \tilde{G}(\delta_Y)
\]
by means of which any irreducible representation becomes isomorphic to a tensor product of representations of the two factors.
Since the Lie algebra of~$G(\delta_W)$ is simple, any non-trivial representation of $G(\delta_W)$ is almost faithful in the sense that its kernel is finite. Since the tensor product of any almost faithful representation with its dual contains the adjoint representation, we therefore obtain an inclusion
\[
 \omega(\Ad_W) \;=\; \one \boxtimes Ad \; \hookrightarrow \; \omega(\delta_Y) \otimes \omega(\delta_Y)^\vee \;=\; \omega(\delta_Y * \delta_{-Y})
\]
in $\Rep(G)$. But then $\Ad_W \hookrightarrow \delta_Y * \delta_{-Y}$ and hence $\Supp(\Ad_W)\hookrightarrow Y-Y \subset A$. It follows that 
\[
 \dim \Supp(\Ad_W) \;\leq\; \dim(Y-Y) \;\leq\; 2\dim(Y) \;=\; 2\max \{\, \dim(Y_i) \mid i=1,\dots, n\, \} 
\]
and therefore $\dim(Y_i) \geq \tfrac{1}{2} \dim \Supp(\Ad_W)$ for at least one $i$. \qed

\medskip

Schreieder has shown that Jacobians of smooth projective curves are the only indecomposable principally polarized abelian varieties whose theta divisor is a sum of curves~\cite[cor.~3]{SchreiederTheta}. We do not reprove his result, but we easily recover:

\begin{cor} \label{cor:sumofcurves}
The theta divisor on the intermediate Jacobian of a smooth cubic threefold is not a sum of curves.
\end{cor}

{\em Proof.} For some translate of the theta divisor $\Theta \subset A$ on the intermediate Jacobian we know that $\omega(\delta_\Theta)\in \Rep(G(\delta_\Theta))$ is the adjoint representation~\cite[th.~2]{KraemerThreefolds}, hence we have $\Ad_\Theta = \delta_\Theta$. Since the intermediate Jacobian has dimension $g=5$, proposition~\ref{prop:sumofsubvarieties} shows that any decomposition of the theta divisor as a sum involves at least one summand of dimension at least $\dim \Supp(\Ad_\Theta)/2 = (g-1)/2 = 2$. \qed

\subsection{A Torelli theorem for subvarieties} \label{sec:torelli}

To rigidify the situation let us now fix a very general~$u\in V(\bbC)$ and work in the setting of theorem~\ref{thm:gausstorus}.
The following result generalizes the Torelli theorem for curves~\cite{AndreottiTorelli} and for Fano surfaces of cubic threefolds~\cite[th.~13.11]{CGIntermediateJacobian}:

\begin{cor} \label{cor:torelli}
Fix a simple module $\scrM \in \Mhol(\scrD_A)$. For $i=1,2$, let~$Y_i \subset A$ be subvarieties such that the following properties hold: \smallskip
\begin{enumerate}
\item $\CC(\delta_{Y_i})$ is reduced, irreducible and not stable under any translation, \smallskip
\item $\scrM \in \langle \delta_{Y_i} \rangle$ and the induced map
$p_i:  G_u(\delta_{Y_i}) \twoheadrightarrow G_u(\scrM)$ is an isogeny, \smallskip
\item there is an isomorphism $\varphi: G_u(\delta_{Y_1}) \stackrel{\sim}{\longrightarrow} G_u(\delta_{Y_2})$ with $\varphi^*(\omegaMicro{u}(\delta_{Y_2})) \simeq \omegaMicro{u}(\delta_{Y_1})$ making the following diagram commute:
\[
\xymatrix@M=0.5em@C=1em@R=2em{
 G_u(\delta_{Y_1}) \ar[dr]_-{p_1} \ar[rr]^-\sim_-{\varphi} && G_u(\delta_{Y_2}) \ar[dl]^-{p_2}\\
 & G_u(\scrM) &
}
\]
\end{enumerate}
Then $Y_1 = Y_2 + a$ for some $n$-torsion point $a\in A(\bbC)$, with $n=\deg(p_1)=\deg(p_2)$.
\end{cor}

{\em Proof.} Put $G=G_u(\scrM)$, $G_i = G_u(\delta_{Y_i})$, and let $T \subseteq G$ and $T_i \subseteq G_i$ be the subtori that are the connected components of the subgroups of multiplicative type in theorem~\ref{thm:gausstorus}. By assumption $\scrM \in \langle \delta_{Y_i} \rangle$, so for some $n\in \bbN$ we have an embedding as a direct summand
\[
 \scrM \;\hookrightarrow\; \scrN^{*n} \;=\; \scrN * \cdots * \scrN 
 \quad \textnormal{for} \quad 
 \scrN = \delta_{Y_i} \oplus  \delta_{-Y_i}.
\]
Kashiwara's estimate for the characteristic variety of direct images~\cite[th.~4.27]{KashiwaraDModules} then says 
\[
 \Char(\scrM) \;\subseteq\; \Char(\scrN) * \cdots * \Char(\scrN),
\]
where on the right hand side we use the natural convolution product for subvarieties of the cotangent bundle: For $\Lambda_1, \Lambda_2 \subseteq T^* A = A\times V$ we put
$
 \Lambda_1 * \Lambda_2 = 
 \varpi(\Lambda_1 \times_V \Lambda_2)$
where
\[
 \varpi: \quad T^*A \times_V T^*A \;=\; A\times A \times V \;\longrightarrow\; T^* A \;=\; A\times V, \quad 
 (z_1, z_2, v) \;\mapsto\; (z_1+z_2, v)
\]
is induced by the addition morphism, see section~\ref{sec:localsystems}. It follows from this description that any $(a,u) \in \Char(\scrM)$ satisfies $a = a_1 + \cdots + a_n$ for certain~$(a_i, u)\in \Char(\scrN)$; for the groups generated by these points on the abelian variety, we therefore obtain an inclusion
\[
  \left \langle \, a\in A(\bbC) \mid (a, u) \in \Char(\scrM) \, \right \rangle
 \;\subseteq\;
 \left \langle \, a\in A(\bbC) \mid (a, u) \in \Char(\delta_{Y_i}) \, \right \rangle.
\]
The Cartier dual of an inclusion is an epimorphism, hence the isogeny~$p_i: G_i \twoheadrightarrow G$ restricts to an isogeny $T_i \twoheadrightarrow T$ of tori. Then by connectedness~$T_i$ coincides with the connected component of the preimage $p_i^{-1}(T) \subseteq G_i$, and we have $\varphi(T_1)=T_2$ since~$p_2\circ \varphi = p_1$. The same centralizer argument as in the proof of corollary~\ref{cor:weyl} shows that each $T_i$ lies in a unique maximal torus $Z_i\subseteq G_i$. The uniqueness of these maximal tori implies that $\varphi(Z_1) = Z_2$ and that the maximal torus $Z=p_1(Z_1)=p_2(Z_2)\subseteq G$ is uniquely determined as well. So we have found a distinguished maximal torus in each of our groups such that $p_1$, $p_2$ and $\varphi$ map these tori onto each other. Looking at weight spaces for these tori we want to identify suitable irreducible components in the characteristic varieties of convolution powers of $\delta_{Y_i}$ which suffice to recover $Y_i$.

\medskip

The kernel $\ker(p_i) \subseteq G_i$ is a central subgroup of order $n$, so by Schur's lemma it acts on the  irreducible representation $U_i = \omegaMicro{u}(\delta_{Y_i}) \in \Rep(G_i)$ as multiplication by $n$-th roots of unity. This action induces the trivial action on the $n$-th tensor power, so via the fully faithful embedding
$p_i^*: \Rep(G) \hookrightarrow \Rep(G_i)$
we may view the symmetric power 
\[
 W_i \;=\; S^n(U_i)
\]
as an object of $\Rep(G)$, corresponding to a module $\scrN_i \in \langle \scrM \rangle$. Note that $\Char(\scrN_i)$ has a unique irreducible component 
$$
 \Delta_i \;\subseteq\;  \Char(\scrN_i) \;\subseteq\; \Char(\delta_{Y_i} * \cdots * \delta_{Y_i})
$$
that is the image of the diagonal under the addition map 
\[
 \Char(\delta_{Y_i}) \times_V \cdots \times_V \Char(\delta_{Y_i})
 \;\longrightarrow \; \Char(\delta_{Y_i}*\cdots * \delta_{Y_i}) 
\]
and hence surjects onto 
$
 n \cdot Y_i = \{ n y \mid y\in Y_i \} \subset A
$.
In fact the weights of $Z$ on~$W_i$ contain a unique Weyl group orbit of weights which pull back to $n$ times a weight of $Z_i$ on $U_i$ via the epimorphism $Z_i \twoheadrightarrow Z$, and since $\CC(\delta_{Y_i})$ is irreducible, this Weyl group orbit corresponds to a single orbit of the monodromy group $\Gal(\gamma)$ by corollary~\ref{cor:weyl}, hence to an irreducible component of $\Char(\scrN_i)$.

\medskip

Now by assumption we have an isomorphism $\varphi^*(U_2) \simeq U_1$ in $\Rep(G_1)$, hence an isomorphism $W_2\simeq W_1$ in $\Rep(G)$, and under this latter isomorphism our Weyl group orbits correspond to each other because our maximal tori are compatible with $p_1$, $p_2$ and $\varphi$. Geometrically this gives an isomorphism $\scrN_2 \simeq \scrN_1$ in $\langle \scrM \rangle$ and shows that the diagonal components
\[
 \Delta_1, \Delta_2 \;\subseteq\; \Char(\scrN_1) \;=\; \Char(\scrN_2)
\]
coincide. Projecting to $A$ we obtain $n\cdot Y_1 = n\cdot Y_2$, and by taking preimages under the isogeny $[n]: A\to A$ we obtain that $Y_2 = Y_1 + a$ for some $n$-torsion point $a\in A(\bbC)$. \qed

\medskip

The classical Torelli theorem for curves can be recovered from the above result as follows. Let $C_1, C_2$ be smooth projective curves for whose Jacobians there exists an isomorphism
\[
 f: \quad \mathrm{Jac}(C_1) \;\stackrel{\sim}{\longrightarrow}\;
 \mathrm{Jac}(C_2)
\]
of principally polarized abelian varieties. Using any translate of the Abel-Jacobi map,~we regard  $C_i \subset \mathrm{Jac}(C_i)$ as a subvariety of its Jacobian. Then a theta divisor defining the polarization can be written as a sum $W_{g-1}(C_i) = C_i + \cdots + C_i \subset \mathrm{Jac}(C_i)$ of $g-1$ copies of the curve. Since $f$ is an isomorphism of polarized abelian varieties and the principal polarization determines the theta divisor up to a translation, it follows that
\[
 f(W_{g-1}(C_1)) \;=\; W_{g-1}(C_2) + a 
 \quad \textnormal{for some} \quad a \in \mathrm{Jac}(C_2)(\bbC).
\] 
Replacing the two curves by suitable translates inside their Jacobian, we can assume that $a=0$. Then corollary~\ref{cor:torelli} applies to $A=\mathrm{Jac}(C_2)$ with $Y_1 = f(C_1)$, $Y_2 = C_2$. As a reference object to compare the Tannakian data corresponding to these two curves we take $\scrM = \delta_{W_{g-1}(C_2)}$, since by construction we have $W_{g-1}(C_2) = Y_i + \cdots + Y_i$ and hence 
\[
 \scrM \;\hookrightarrow\; \delta_{Y_i} * \cdots * \delta_{Y_i} \in \langle \delta_{Y_i} \rangle
\]
for $i=1,2$. The Tannaka groups of the two curves are isomorphic: By example~\ref{ex:jacobians}(a)  they only depend on whether the curve is hyperelliptic or not, and if one of the two curves is hyperelliptic, then so must be the other because this can be read off from the dimension of the singular locus of its theta divisor. Since the groups $G_u(\delta_{Y_i})$ are simply connected and the projections $p_i: G_u(\delta_{Y_i})\twoheadrightarrow G_u(\scrM)$ are isogenies, we can find an isomorphism
\[
 \varphi: \quad G_u(\delta_{Y_1}) \;\stackrel{\sim}{\longrightarrow}\; G_u(\delta_{Y_2})
 \quad \textnormal{with} \quad 
 p_1 = p_2 \circ \varphi.
\]
The discussion of the defining representations in example~\ref{ex:jacobians}(a) shows that possibly after replacing one of the two curves by its negative, we can assume that we have an isomorphism $\omega_u(\delta_{Y_1}) \simeq \varphi^*(\omega_u(\delta_{Y_2}))$. Then corollary~\ref{cor:torelli} says $Y_1 = Y_2 + a$ for some point $a\in A(\bbC)$, and it follows that the two curves are isomorphic.

\medskip

Note that in the classical Torelli theorem there is no preferred choice between the Abel-Jacobi image of a curve and its negative, the two are translates of each other only in the hyperelliptic case. This is reflected in the above representation theoretic argument by the fact that there is no preferred choice between the natural representation of $\Sl_{2g-2}(\bbC)$ and its dual representation.

\section{The Fourier-Mukai transform} \label{sec:fouriermukai}

In this section we discuss an interpretation of Tannaka groups as structure groups of certain principal bundles defined via the Fourier-Mukai transform. As an application we obtain that these groups are almost connected. 

\subsection{Algebraic principal bundles} \label{nori}

 Let $G$ be a linear algebraic group and $U$ an algebraic variety over the complex numbers. By a {\em principal bundle} with structure group $G$ on $U$ we mean an algebraic variety $\scrG$ that is equipped with a surjective flat affine morphism $\scrG\rightarrow U$ and with a right action~$m: \scrG \times G \rightarrow \scrG$ over $U$ such that the diagram
\[
\xymatrix@M=0.5em{
\scrG \times G \ar[r]^-m \ar[d]_-p & \scrG \ar[d] \\
\scrG \ar[r] & U
}
\]
is Cartesian, where $p$ denotes the projection onto the first factor. For any algebraic representation $E \in \Rep(G)$ we then get an associated vector bundle $\scrE = \scrG \times^G E$ by the usual contracted product. This gives a $\bbC$-linear exact tensor functor
\[
 \Phi \;=\; \Phi_\scrG: \quad \Rep(G) \; \longrightarrow \; \Mcoh(\scrO_U)
\]
to the category of coherent sheaves with the property that $\Phi(E)$ is a locally free sheaf of rank $\dim(E)$ for all $E\in \Rep(G)$. Conversely Nori has shown in~\cite[sect.~2]{NoriRepresentations} that any $\bbC$-linear exact tensor functor $\Phi: \Rep(G) \rightarrow \Mcoh(\scrO_U)$ with this latter property arises from a principal bundle in the above way and that up to isomorphism the bundle is uniquely determined by the functor. In the next section we will apply this to $G = G(\scrM)$ for $\scrM \in \Mhol(\scrD_A)$, with $\Phi$ defined via the Fourier-Mukai transform.

\medskip

We say that an algebraic vector bundle has a {\em reduction} to a principal bundle~$\scrG$ if it lies in the essential image of the functor $\Phi_\scrG$. Any vector bundle of rank $n$ has a reduction to a principal bundle with structure group~$\Gl_n(\bbC)$, but usually also many smaller ones. A reduction of a given vector bundle to a principal bundle~$\scrH$ is called a~{\em natural reduction} if \smallskip
\begin{itemize}
 \item the structure group $H$ of the principal bundle $\scrH$ is reductive, and \smallskip
 \item for any other reduction of the given vector bundle to a principal bundle~$\scrG$ with a reductive structure group $G$ there exists an embedding $\iota: H\hookrightarrow G$ such that the following diagram commutes:
\[
\xymatrix@M=0.5em{
 \Rep(G) \ar[rr]^-{\Phi_\scrG} \ar@{..>}[dr]_-{\iota^*} && \Mcoh(\scrO_U) \\
 & \Rep(H) \ar[ur]_-{\Phi_\scrH} &
}
\]
\end{itemize}
A result of Bogomolov~\cite[th.~2.1]{BogomolovStable} says that on varieties $U$ with $\dim H^0(U, \scrO_U) = 1$ any algebraic vector bundle has a unique such natural reduction, and for stable vector bundles on projective varieties the corresponding minimal reductive structure group is an algebraic analog of the holonomy of Chern connections~\cite{BKHolonomy}. We will be guided by this analogy, but replace stability by semisimplicity of our categories.

\subsection{The Fourier-Mukai transform} \label{subsec:fouriermukai}

Let $A^\natural$ be the moduli space of coherent line bundles with a flat connection on the abelian variety $A$ as in section~\ref{sec:fouriermukai-intro}, and denote by $\varphi: A^\natural \to \hat{A} = \Pic^\circ(A)$ the forgetful map. Let~$p_1$, $p_2$ be the two projections in the following diagram:
\[
\xymatrix@M=0.5em@C=2em@R=1em{
 & A\times A^\natural \ar[dl]_-{p_1} \ar[dr]^-{p_2} \\
 A && A^\natural
}
\]
The pullback $\scrP^\natural$ of the Poincar\'e bundle under the map $\id \times \varphi: A\times A^\natural \to A \times \hat{A}$ has a universal relative flat connection 
\[
  \nabla^\natural: \quad \scrP^\natural  \;\longrightarrow \; \Omega^1_{A\times A^\natural / A^\natural} \otimes_{\scrO_{A\times A^\natural}} \scrP^\natural
\]
and we consider the Fourier-Mukai transform 
\begin{align} \nonumber
 \FM_A: \quad \Dbcoh(\scrD_A) \; \stackrel{\sim}{\longrightarrow} \; & \; \Dbcoh(\scrO_{A^\natural}),\\ \nonumber
  \quad \scrM \; \mapsto \; & \; Rp_{2,*}\, \DR_{A\times A^\natural / A^\natural}(p_1^*(\scrM) \otimes (\scrP^\natural, \nabla^\natural) ),
\end{align}
as in~\cite{LaumonTransformation, RothsteinSheaves}. In what follows we use the shorter notation $\FM = \FM_A$ when there is no risk of confusion about the abelian variety under consideration. Recall that we use right $\scrD$-modules, so 
\[
 \DR_{A\times A^\natural / A^\natural}(\scrN) \,=\, \Bigl[ \, \cdots 
  \longrightarrow  \scrN \otimes_{\scrO_{A\times A^\natural}} \scrT_{A\times A^\natural / A^\natural}
  \longrightarrow  \scrN 
  \Bigr]
\]
where $\scrT_{A\times A^\natural / A}$ is the relative tangent sheaf. We place the complex in non-positive degrees to make it relatively perverse, with the final term $\scrN$ in degree zero. The reason for using right rather than left modules are the signs for the commutativity constraints:

\begin{propo} \label{prop:fmtensor}
The functor $\FM: \Dbhol(\scrD_A) \longrightarrow \Dbcoh(\scrO_{A^\natural})$ underlies a tensor functor with respect to the tensor structures given by the convolution product on the source and the usual derived tensor product on the target. 
\end{propo}

{\em Proof.} For any homomorphism $f: B\to A$ of abelian varieties, let $f^\natural: A^\natural \to B^\natural$ be the induced homomorphism between the corresponding moduli spaces of rank one flat connections. By~\cite[prop.~3.3.2]{LaumonTransformation} we have for any $\scrM \in \Dbhol(\scrD_B)$ a natural isomorphism 
\[
 \FM_A(f_\dag(\scrM)) \;\stackrel{\sim}{\longrightarrow}\; 
 Lf^{\natural *}(\FM_B(\scrM)).
\]
We apply this to the addition morphism $a: A\times A \to A$. Using the definition of the convolution product, the above isomorphism, the compatibility of the Fourier-Mukai transform with external tensor product and the fact that $a^\natural: A^\natural \to A^\natural \times A^\natural$ is the diagonal embedding, we get
\begin{eqnarray*} 
\FM_A(\scrM_1 * \scrM_2)
&\;=\;&
\FM_A(a_\dag(\scrM_1 \boxtimes \scrM_2)) \\[0.2em]
&\;\simeq\;&
La^{\natural*}(\FM_{A\times A}(\scrM_1 \boxtimes \scrM_2)) \\[0.2em]
&\;\simeq\;&
La^{\natural*}(\FM_A(\scrM_1) \boxtimes \FM_A(\scrM_2)) \\[0.2em]
&\;\simeq\;&
\FM_A(\scrM_1)\otimes^L_{\scrO_{A^\natural}} \FM_A(\scrM_2)
\end{eqnarray*}
for $\scrM_1, \scrM_2 \in \Dbhol(\scrD_A)$. We claim that $\FM_A$ is a tensor functor via the composite isomorphisms
\[
 \FM_A(\scrM_1 * \scrM_2) \; \stackrel{\sim}{\longrightarrow} \; \FM_A(\scrM_1) \otimes^L_{\scrO_{A^\natural}} \FM_A(\scrM_2),
\]
i.e.~that these isomorphisms satisfy the usual compatibilities with the associativity, commutativity and unit constraints. The crucial part is the commutativity constraint where we must check that the right signs appear. We claim that the following diagram commutes, where for a variety $X$ we write $\sigma_X: X\times X \to X\times X, (x,y)\mapsto (y,x)$ for the involution that interchanges the two factors and where all the arrows denote the natural isomorphisms:
\[
\xymatrix@R=1em@C=1em@M=0.5em{
\FM(\scrM_1) \otimes^L_{\scrO_{A^\natural}} \FM(\scrM_2) \ar[dd]
\ar[rr]
&& \FM(\scrM_2) \otimes^L_{\scrO_{A^\natural}} \FM(\scrM_1) \ar[dd]
\\ 
& \circled{1} &
\\ 
La^{\natural *}(\FM(\scrM_1) \boxtimes \FM(\scrM_2)) \ar[dddd]
\ar[rr] 
&&  La^{\natural *} ( \sigma_{A^\natural, *} (\FM(\scrM_2) \boxtimes \FM(\scrM_1))) \ar[dd]
\\ \\
& \circled{2} & La^{\natural *} (\sigma_{A^\natural, *} (\FM(\scrM_2 \boxtimes\scrM_1))) \ar[dd]
\\ \\
La^{\natural*}(\FM(\scrM_1 \boxtimes \scrM_2)) \ar[dd]
\ar[rr] 
&& La^{\natural*}(\FM(\sigma_{A,*}(\scrM_2 \boxtimes \scrM_1))) \ar[dd]
\\
& \circled{3} & 
\\ 
\FM(a_\dag(\scrM_1 \otimes \scrM_2)) 
\ar[rr] 
&& \FM(a_\dag(\sigma_{A,*}(\scrM_2 \boxtimes \scrM_1)))
}
\]
Here $\circled{1}$ uses the identification $a\circ \sigma = a$, while the commutativity of the square $\circled{3}$ follows from the naturality of the isomorphism $La^{\natural*} \circ \FM \simeq \FM\circ a_\dag$. The middle square $\circled{3}$ arises by applying $La^{\natural*}$ and the K\"unneth isomorphisms for the pushforward under 
\[
 (p_2, p_2): \quad A\times A^\natural \times A \times A^\natural \;\longrightarrow\; A^\natural \times A^\natural
\]
to the following square, where we write $\scrN_i := p_1^*(\scrM_i) \otimes (\scrP^\natural, \nabla^\natural)$:
\[
\xymatrix@C=2em@M=0.5em{
\DR_{A\times A^\natural/A^\natural}(\scrN_1) \boxtimes \DR_{A\times A^\natural/A^\natural}(\scrN_2)
\ar[r] \ar[dd]
& \sigma_* 
(
\DR_{A\times A^\natural/A^\natural}(\scrN_1) \boxtimes \DR_{A\times A^\natural/A^\natural}(\scrN_2)
)
\ar[d]
\\
& 
\sigma_{A\times A^\natural, *} (\DR_{(A\times A^\natural)^2/(A^{\natural})^2}(\scrN_1\boxtimes \scrN_2))
\ar[d]
\\
\DR_{(A\times A^\natural)^2/(A^{\natural})^2}(\scrN_1 \boxtimes \scrN_2) 
\ar[r] 
& \DR_{(A\times A^\natural)^2/(A^{\natural})^2}(\sigma_{A\times A^\natural, *}(\scrN_1\boxtimes \scrN_2))
}
\]
This amounts to the compatibility of the relative de Rham functor for {\em right} $\scrD$-modules with the action of the symmetric group on external tensor products and follows by the same sign computation as in~\cite[prop.~1.5]{MaximSaitoSchuermann}. \qed

\medskip

We now want to apply the generic vanishing property~\eqref{eq:gvt} from section~\ref{sec:fouriermukai-intro} to get a tensor functor with values in locally free sheaves on any finitely generated tensor subcategory of $\M(A)$. However, by definition $\M(A)$ is the quotient of $\Mhol(\scrD_A)$ by the Serre subcategory of {\em all} negligible modules, whereas for the construction of principal bundles we want to work on a {\em fixed} open subset $U\subseteq A^\natural$ and only take the quotient by those negligibles whose Fourier-Mukai transform vanishes on this given subset. So let $\Mhol(\scrD_A,U) \subset \Mhol(\scrD_A)$ be the full abelian subcategory of all holonomic modules with the property that the generic vanishing condition \eqref{eq:gvt} holds on $U$ for all their subquotients and their duals, and denote by $\M(A, U)$ the quotient of this full abelian subcategory by the Serre subcategory of all negligible modules contained in it. We have a commutative diagram
\[
\xymatrix@M=0.5em{
 \Mhol(\scrD_A, U) \ar@{^{(}->}[r] \ar[d] & \Mhol(\scrD_A) \ar[d] \\
 \M(A, U) \ar@{^{(}->}[r]^-{\exists! \, i} & \M(A)
}
\]
where $i$ comes from the universal property of the quotient $\M(A)$ and is a fully faithful embedding by~\cite[lemma 12.3]{KrWVanishing}. As in theorem~13.2 of loc.~cit.~its essential image is stable under convolution, making $\M(A, U)$ a rigid abelian tensor category in a way compatible with our previous constructions: For any~$\scrM \in \Mhol(A, U)$, the tensor categories that it generates in~$\M(A, U)$ and in~$\M(A)$ are equivalent. Using the tensor functor
\[
 \Phi: \quad \M(A, U) \;\longrightarrow\; \Mcoh(\scrO_U), \quad \scrM \;\mapsto \; \scrH^0(\FM(\scrM))|_U,
\]
we can rephrase theorem~\ref{thm:structure_group} as follows:

\begin{thm} \label{structure_group}
For any $\scrM\in \Mhol(\scrD_A, U)$, the vector bundle $\Phi(\scrM) \in \Mcoh(\scrO_U)$ is induced by a principal bundle 
\[ \scrG \;\longrightarrow \; U \] 
whose structure group is isomorphic to the Tannaka group $G=G(\scrM)$. If $\scrM$ is semisimple, then $\scrG$ is a natural reduction for this vector bundle.  
\end{thm}

{\em Proof.} By the generic vanishing property~\eqref{eq:gvt} the functor $\Phi$ is exact and takes values in locally free sheaves, and it is a tensor functor by proposition~\ref{prop:fmtensor}. So the composite functor
\[
\xymatrix@M=0.5em@C=2.5em{
 \Rep(G) \ar[r]^-{\sim} 
 & \langle \scrM \rangle \;\subset\; \M(A, U) \ar[r]^-{\Phi}
 & \Mcoh(\scrO_U)
}
\]
satisfies Nori's properties of section~\ref{nori}: It is a $\bbC$-linear exact tensor functor sending each representation to a locally free sheaf of the same rank. Thus $\Phi =\Phi_\scrG$ for some principal bundle~$\scrG$ as claimed in the first part of the theorem. For the second part notice that  
\[
 \dim H^0(U, \scrO_U) \;=\; \dim H^0(A^\natural, \scrO_{A^\natural}) \;=\; 1,
\]
where the first identity holds since we assumed that the complement of $U\subseteq A^\natural$ is of codimension at least two, and the second is shown in~\cite[th.~2.4.1]{LaumonTransformation}. So we know by the result of Bogomolov~\cite[th.~2.1]{BogomolovStable} that the vector bundle $\Phi(\scrM)$ has a natural reduction $\scrH$. Let~$H$ be the corresponding minimal reductive structure group. If~$\scrM$ is semisimple, then it corresponds to a semisimple representation of the Tannaka group $G=G(\scrM)$. This representation is moreover faithful by construction, hence it follows that the unipotent radical of the Tannaka group $G$ is trivial, i.e.~$G$ is a reductive group. Hence by the minimality of a natural reduction there exists an embedding $\iota: H\hookrightarrow G$ such that the diagram
\[
\xymatrix{
 \Rep(G) \ar[rr]^-{\Phi_\scrG} \ar@{..>}[dr]_-{\iota^*} && \Mcoh(\scrO_U) \\
 & \Rep(H) \ar[ur]_-{\Phi_\scrH} &
}
\]
commutes. Now any reductive subgroup of a reductive group is determined by its invariants in the tensor powers of a faithful self-dual representation~\cite[prop.~3.1(c)]{DeligneHodgeCycles}. It follows that if $\iota$ were not an isomorphism, we could find a non-trivial irreducible representation $\one \not\simeq W\in \Rep(G)$ whose restriction $\iota^*(W)\simeq \one \oplus W'$ contains a trivial representation as a direct summand. On the geometric side this would give us a simple object $\scrN \in \langle \scrM \rangle$ with $\scrN \not \simeq \delta_0$ such that on the open subset $U\subseteq A^\natural$ we have a splitting
\[
 \Phi(\scrN) \;=\; \scrH^0(\FM(\scrN))|_U \;\simeq \; \scrO_U \oplus \scrE
\]
as a direct sum of two vector bundles, one of which is trivial. Any simple object of the quotient category $\rmM(A, U)$ can be represented by a simple object in $\Mhol(\scrD_A, U)$, and the latter is determined uniquely up to isomorphism; so in what follows we view~$\scrN$ as a simple holonomic $\scrD_A$-module. Now the reconstruction result in~\cite[cor.~21.3]{SchnellHolonomic} says that for any simple non-negligible $\scrN \in \Mhol(A, U)$, its Fourier-Mukai transform can be reconstructed in a natural way from the coherent vector bundle $\scrF = \Phi(\scrN)$; more precisely
\begin{equation} \label{eq:reconstruction}
 \FM(\scrN) \;\simeq\; 
 (\tau_{\leq \ell - 1} \circ \Delta \circ \tau_{\leq \ell - 2} \circ \Delta \circ \cdots \circ \tau_{\leq 1} \circ \Delta)(j_*(\scrF))
\end{equation}
for any odd $\ell \geq \dim(A)$, where $\Delta = \RHom(-, \scrO_{A^\natural})$ and $j: U\hookrightarrow A^\natural$ denotes the open embedding. It follows by functoriality of the right hand side that the decomposition of the vector bundle $\scrF\simeq \scrO_U \oplus \scrE$ extends to a decomposition of $\FM(\scrN)$ as a direct sum of two coherent sheaf complexes. Note that we are working with algebraic coherent sheaves throughout, so the Fourier-Mukai transform has an inverse; so the previous decomposition comes from a decomposition $\scrN \simeq \delta_0 \oplus \scrN'$ for some $\scrN'\in \Mhol(\scrD_A)$, which contradicts our simplicity assumption.  \qed

\medskip

\subsection{Almost connectedness} \label{sec:almostconnected}

The above easily implies that the Tannaka groups are almost connected in the sense that their group of connected components is a finite abelian group Cartier dual to a finite group of points on the abelian variety, see theorem~\ref{thm:componentgroup}. For perverse sheaves this is an unpublished result by Weissauer~\cite{WeConnected}, who proved the corresponding statement also over finite fields; our proof is completely different, it works only in characteristic zero  but in contrast to loc.~cit.~it also applies to irregular holonomic $\scrD_A$-modules:

\begin{thm} \label{componentgroup}
For any $\scrM\in \Mhol(\scrD_A)$ with Tannaka group $G=G(\scrM)$ we have an epimorphism
\[
 \pi_1(\hat{A}, 0) \; \twoheadrightarrow \;  G/G^\circ
\]
from the fundamental group of the dual abelian variety onto the group of connected components of $G$. In particular $G/G^\circ$ is a finite abelian group.
\end{thm}

{\em Proof.} Let $\pi = G/G^\circ$ be the finite group of connected components and $W=\bbC[\pi]$ its regular representation. Then $W$ corresponds via the Tannakian formalism to some module $\scrN \in \Mhol(\scrD_A)$ with Tannaka group $G(\scrN)=\pi$ because $\pi$ acts faithfully on the regular representation. Replacing the original module $\scrM$ with $\scrN$ we may hence assume that the group $G=G(\scrM)$ is itself finite and in particular reductive, so $\scrM$ can be assumed to be semisimple. 

\medskip

Now consider as above the algebraic vector bundle $\scrE=\Phi(\scrM)$ over a Zariski open subset $U \subseteq A^\natural$ whose complement has codimension at least two. By theorem~\ref{structure_group} this vector bundle has a natural reduction to a principal bundle $\scrG$ with structure group $G$. By minimality the total space~$\scrG$ is connected: Otherwise the stabilizer of a connected component $\scrH \subset \scrG$ would be a proper subgroup $H \subset G$, and the composite morphism
\[
 \scrH \times^H G \;\hookrightarrow\;
 \scrG \times^H G \;\twoheadrightarrow\; \scrG \times^G G \;=\; \scrG 
\]
would be a morphism of principal bundles, hence an isomorphism. So $\scrH$ would reduce $\scrG$ to the smaller reductive structure group $H \subset G$, a contradiction. 

\medskip

As the structure group of our principal bundle is finite, we get that $\scrG \rightarrow U$ is a connected finite \'etale Galois cover with Galois group $G$. Fixing $u\in U(\bbC)$, any such cover is given by an epimorphism 
\[ \pi_1(U, u) \; \twoheadrightarrow \; G \;=\; \Aut(\scrG/U). \]
But
$
 \pi_1(U, u) = \pi_1(A^\natural, u) = \pi_1(\hat{A}, 0)
$
since the complement of the subset $U \subseteq A^\natural$ has codimension at least two and since~$A^\natural \twoheadrightarrow \hat{A}$ is a topologically trivial fibration with contractible fibers, so the result follows. \qed

\medskip

\begin{lem} \label{lem:character}
Conversely any character of the Tannaka group corresponds to a skyscraper sheaf: Any simple holonomic $\scrD_A$-module $\scrN \in \Mhol(\scrD_A, U)$ with $\dim (\omegaFM{u}(\scrN)) = 1$ has the form $\scrN \simeq \delta_a$ for some $a\in A(\bbC)$.
\end{lem}

{\em Proof.} By~\cite[prop.~21.1]{SchnellHolonomic} we know that the coherent sheaf $\scrF = \scrH^0(\FM(\scrN))$ is reflexive, and by construction its rank is $\dim(\omegaFM{u}(\scrN)) = 1$. Since the variety $A^\natural$ is smooth, any reflexive sheaf of rank one on it is a line bundle by~\cite[prop.~1.9]{HartshorneStable}. So it follows that $\scrF$ is a line bundle, and then the formula~\eqref{eq:reconstruction} for the Fourier-Mukai transform shows that the complex $\FM(\scrN)$ is quasi-isomorphic to the line bundle $\scrF$ placed in degree zero. Moreover, the first Chern class of this line bundle must vanish by~\cite[prop.~24.1]{SchnellHolonomic}, so altogether
$
\FM(\scrM) \in \Pic^\circ(A^\natural)\simeq\Pic^\circ(\hat{A})\simeq A.
$
If $a\in A(\bbC)$ denotes the corresponding point, it follows from the definitions that $\scrM \simeq \delta_a$. \qed 

\medskip 

For convenience of the reader we also include the interpretation of direct and inverse images under an isogeny $p: A\twoheadrightarrow B$ in terms of restriction and induction functors as stated in corollary~\ref{cor:isogeny}, which follows from theorem~\ref{componentgroup} as in~\cite{WeConnected}.

\medskip

{\em Proof of corollary~\ref{cor:isogeny}.} {\em (a)} The direct and inverse images for~$p$ descend to a pair of exact adjoint functors between $\M(A)$ and $\M(B)$. We want to restrict these to a pair of adjoint functors between suitable finitely generated tensor subcategories, starting from $\scrM \in \Mhol(\scrD_A)$. The direct image $\scrN = p_\dag(\scrM)$ has $p^*(\scrN) \simeq \bigoplus_{x\in \ker(p)} t_x^*(\scrM)$ for the translations $t_x: A\rightarrow A$. Applying the direct image functor again we arrive at the object \medskip 
\[
 p_\dag(p^*(\scrN)) \;\simeq \! \bigoplus_{x\in \ker(p)} p_\dag(t_x^*(\scrM)) 
 \;\simeq \! \bigoplus_{x\in \ker(p)} p_\dag(\scrM) 
 \;\simeq \; \scrN^{\oplus d} \bigskip
\]
for $d=\deg(p)$. Since this is again an object inside the category $\langle \scrN \rangle$, we obtain a pair 
\[
\xymatrix@M=1em{
 p^* \, : \, \langle \scrN \rangle \ar@<0.7ex>[r] &
 \langle p^*(\scrN) \rangle \, : \, p_\dag \ar@<0.7ex>[l]
}
\]
of exact adjoint functors. Here the right adjoint $p_\dag$ is a tensor functor corresponding to the restriction functor for an embedding $\iota: G(\scrN)\hookrightarrow G(p^*(\scrN))$, so the left adjoint~$p^*$ is the corresponding induction functor. Usually the latter exists only on the level of infinite dimensional algebraic representations. Since here it already exists in the finite dimensional setup, the image of $\iota$ is a subgroup of finite index; but any object of $\langle \scrN \rangle$ is a subquotient of the direct image of some object in $\langle \scrM \rangle$, so the composite homomorphism $G(\scrN)\hookrightarrow G(p^*(\scrN)) \twoheadrightarrow G(\scrM)$ is injective and the claim now easily follows. 

\medskip 

{\em (b)} We apply the previous result to $p: A\twoheadrightarrow B= A/K$ where $K\subset A(\bbC)$ denotes the finite subgroup which is Cartier dual to the group of connected components of the group~$G=G(\scrM)$. Note that if $\scrM$ is semisimple, then $p_\dag(\scrM)$ is semisimple. We have:
\begin{align*}
 \omegaFM{u}(\scrM)|_{G^\circ} \;\; \textnormal{is irreducible} &\;\;\Longleftrightarrow\;\; p_\dag(\scrM) \;\; \textnormal{is irreducible} \\
 &\;\;\Longleftrightarrow\;\; \dim \End_{\rmM(B)}(p_\dag(\scrM)) = 1 \\
 &\;\;\Longleftrightarrow\;\; \dim \Hom_{\rmM(A)}(p^*p_\dag(\scrM), \scrM) = 1
\end{align*}
So we are done because $p^*p_\dag(\scrM) \simeq \bigoplus_{x\in K} t_x^*(\scrM)$. \qed

\section{Tensor categories of germs} \label{sec:germs}

We now pass to the proof of our main theorem~\ref{thm:gausstorus}. In the next section we will construct an exact tensor functor on $\M(A)$ with values in germs of local systems on characteristic varieties. In this preliminary section we introduce possible target categories and relate them to subgroups of multiplicative type and their normalizers.

\subsection{Germs of local systems} \label{sec:localsystems}

The target categories will be endowed with a convolution product defined via the diagram
\[
\xymatrix@M=0.5em{
 &   A^2 \times V \ar@{=}[d] \ar[dl]_-{\id \times \delta} \ar[dr]^-{a\times \id} & \\
 A^2 \times V^2 \ar@{=}[d] & A^2 \times_A T^*A \ar[dr]^-{\varpi}  \ar[dl]_-{\rho} & A\times V \ar@{=}[d] \\
 T^*A^2 && T^*A  
}
\]
where $\rho$ and $\varpi$ are induced by the addition morphism $a: A^2\rightarrow A$ and its differential, the diagonal map $\delta: V\hookrightarrow V^2$ on cotangent spaces. For the moment the cotangent space $V = H^0(A, \Omega_A^1)$ can be replaced by any irreducible variety $U$ and we still denote by
\[
 \varpi \;=\; a\times \id: \; A^2 \times U \;\longrightarrow \;A \times U \quad \textnormal{and} \quad \rho \;=\; \id\times \delta: \; A^2\times U \;\longrightarrow\; A^2\times U^2
\]
the maps induced by the addition morphism and by the diagonal. We are mainly interested in the case where $U\subseteq V$ is an open subset of the cotangent space.

\medskip

We want to consider local systems on closed subvarieties $\Lambda \subset A\times U$ such that the Gauss map $\gamma: \Lambda \twoheadrightarrow U$ given by the projection to the second factor is a finite \'etale cover. It is natural to define the {\em convolution} of such subvarieties $\Lambda_1, \Lambda_2 \subset A\times U$ by
\[
 \Lambda_1 * \Lambda_2 \;=\; \varpi(\Lambda_1 \times_U \Lambda_2) \;\subset\; A\times U
\]
but usually this convolution has singularities and its Gauss map will no longer be a finite \'etale cover. Nevertheless, for any $u\in U$ one can find a small analytic open subset $U_0 \subseteq U$ containing $u$ such that $(\Lambda_1 * \Lambda_2)\cap (A\times U_0)$ is a union of irreducible analytic subvarieties whose Gauss maps are local isomorphisms. So consider the category $\LS(A, U)$ of all sheaves $\scrF$ of complex vector spaces on $A\times U$ that restrict over any sufficiently small classical open $U_0 \subseteq U$ to a finite sum
\[
 \scrF|_{A\times U_0} \; \simeq \; \bigoplus_\alpha \, i_{\alpha, *}(\scrF_\alpha),
\] 
where $\scrF_\alpha$ are finite rank local systems on analytic subvarieties $i_\alpha: \Lambda_\alpha \hookrightarrow A\times U_0$ that project isomorphically onto $U_0$. This category contains all the local systems that we are interested in, and by construction it is stable under the convolution product  
$$
 *: \quad \LS(A, U) \times \LS(A, U) \; \longrightarrow \;  \LS(A, U), \quad
 \scrF_1 * \scrF_2 \;=\; \varpi_* \rho^{-1}(\scrF_1\boxtimes \scrF_2).
$$
For any point $u\in U$ we introduce a category $\LS(A, u)$ of {\em germs of local systems} as follows: The objects of this category are the pairs $\alpha = (\scrF_\alpha, U_\alpha)$ consisting of a Zariski open neighborhood $U_\alpha \subseteq U$ of $u$ together with an object $\scrF_\alpha\in \LS(A, U_\alpha)$, and morphisms are defined by
\[
 \Hom_{\LS(A,u)}(\alpha, \beta) \;=\; \lim_{\longrightarrow} \; \Hom_{\LS(A, W)}(\scrF_\alpha|_{A\times W}, \scrF_\beta|_{A\times W})
\]
where the limit runs over all Zariski open neighborhoods $W\subseteq U_\alpha \cap U_\beta$ of $u$. This is an abelian category, and we equip it with the convolution product induced by the previous one in the obvious way.

\begin{lem}
The convolution product naturally endows $\LS(A, U)$ and $\LS(A, u)$ for~$u\in U$ with the structure of an abelian tensor category such that the passage to germs is an exact tensor functor $\LS(A, U) \longrightarrow \LS(A, u), \, \scrF \mapsto (\scrF, U)$.
\end{lem}

{\em Proof.}
The convolution product is an exact bifunctor because $\varpi$ is finite on the occuring supports. 
The unit object for the convolution product on $\LS(A, U)$ is the constant sheaf $\one = \bbC_{\Lambda_0}$ on $\Lambda_0 = \{0\} \times U$ and
$
\one \stackrel{\sim}{\longrightarrow} \one * \one
$  
is the obvious unit isomorphism. We define the commutativity and associativity constraints on local systems $\scrF_1, \scrF_2, \scrF_3 \in \LS(A, U)$ as follows.

\medskip

For the commutativity constraint, let $\sigma: A^2 \times U \stackrel{\sim}{\longrightarrow} A^2 \times U$ be the involution given by $\sigma(z_1, z_2, v)=(z_2, z_1, v)$. Then we define the commutativity constraint $\psi_{\scrF_1, \scrF_2}$ as the composite of the isomorphisms in the diagram 
\[
\xymatrix@C=2em@M=0.7em{
 \scrF_1 * \scrF_2 \ar@{..>}[d]^-{\psi_{\scrF_1, \scrF_2}} \ar@{=}[r]
 & \varpi_* \rho^{-1}(\scrF_1\boxtimes \scrF_2) \ar[r]^-{\circled{1}}
 &  \varpi_*\, \sigma_* \,\sigma^{-1}\, \rho^{-1} (\scrF_1\boxtimes \scrF_2)  \ar[d]^-{\circled{2}} \\
 \scrF_2 * \scrF_1
 &
 \varpi_* \rho^{-1}(\scrF_2\boxtimes \scrF_1) \ar@{=}[l]
 & \varpi_* \, \sigma_* \, \rho^{-1} (\scrF_2\boxtimes \scrF_1) \ar[l]_-{\circled{3}}  
}
\]
where 
\begin{enumerate} 
\item[\circled{1}] is induced by the adjunction morphism $\id \rightarrow \sigma_* \sigma^{-1}$ 
\item[\circled{2}] is induced by the natural isomorphism $\sigma^{-1} \rho^{-1}(\scrF_1 \boxtimes \scrF_2) \;\simeq\; \rho^{-1}(\scrF_2\boxtimes \scrF_1)$ 
\item[\circled{3}] is the functoriality of direct images for $\varpi = \varpi \circ \sigma$. 
\end{enumerate} 
For the associativity constraint $\varphi_{\scrF_1, \scrF_2, \scrF_3}$ we use the natural isomorphisms in the following diagram
\[
\xymatrix@C=1.8em@M=0.7em{
(\scrF_1 * \scrF_2)* \scrF_3 \ar@{..>}[d]^-{\varphi_{\scrF_1, \scrF_2, \scrF_3}} \ar@{=}[r] 
& \varpi_* \rho^{-1}( (\scrF_1* \scrF_2) \boxtimes \scrF_3 )\ar[r]
& \varpi_{3*} \rho_3^{-1}((\scrF_1 \boxtimes \scrF_2) \boxtimes \scrF_3)\ar@{=}[d] \\
\scrF_1 * (\scrF_2 * \scrF_3) 
& \varpi_* \rho^{-1}( \scrF_1 \boxtimes (\scrF_2* \scrF_3))\ar@{=}[l] 
& \varpi_{3*} \rho_3^{-1}(\scrF_1 \boxtimes (\scrF_2 \boxtimes \scrF_3))\ar[l]
}
\]
where $\rho_3: A^3 \times U \rightarrow A^3 \times U^3$ and $\varpi_3: A^3 \times U \rightarrow A\times U$ are the diagonal and the addition morphism. One immediately verifies that with these definitions $\LS(A, U)$ is an abelian tensor category. We endow $\LS(A, u)$ with the tensor structure induced by the previous one so that the passage to germs is a tensor functor.
\qed

\medskip

One may show that the tensor category $\LS(A, u)$ is rigid, though we will not need this in what follows. Let $\gamma: A\times U \twoheadrightarrow U$ denote the projection. For $\scrF\in \LS(A, U)$ let~$\gamma_*(\scrF)(u)$ be the fiber of the local system $\gamma_*(\scrF)$ at a point $u\in U(\bbC)$. For fixed $u$ the functor
\[
 \LS(A, U) \; \longrightarrow \; \Vect(\bbC), \quad \scrF \; \mapsto \; \gamma_*(\scrF)(u)
\]
is a fiber functor which factors over the tensor category of germs $\LS(A, u)$.
But our definition of~$\LS(A, u)$ makes sense for all scheme theoretic points, so we can also work at the generic point $\eta \in U$ instead of choosing a closed point $u\in U(\bbC)$. For finitely generated tensor subcategories one can always specialize from the generic point to a closed point, using that each~$\LS(A, u)$ is a full abelian tensor subcategory of $\LS(A, \eta)$:

\begin{lem} \label{lem:specialization}
Any finitely generated tensor subcategory of $\LS(A, \eta)$ is equivalent to a tensor subcategory of $\LS(A, u)$ for some $u\in U(\bbC)$.
\end{lem}

{\em Proof.} Since every object of $\LS(A, \eta)$ is of finite length, any finitely generated tensor subcategory of $\LS(A, \eta)$ has at most countably many isomorphism classes of objects. By the axiom of choice we can replace any small tensor subcategory by one which is~{\em skeletal}~\cite[rem.~2.8.7]{EtingofGelaki}, i.e.~has only one object in each isomorphism class. So we may assume the given subcategory has only countably many objects. Each of these is defined on some Zariski open dense subset, and we can take $u\in U(\bbC)$ to be any point in the intersection of these subsets.  
\qed

\subsection{Germs of vector bundles} \label{sec:vectorbundles}

For a weaker but more flexible framework one can replace local systems by coherent vector bundles as follows. Let $\VB(A, U)$ be the category of analytic coherent sheaves~$\scrF \in \Mcoh(\scrO_{A\times U}^\an)$ with algebraic support such that the projection $\gamma: A\times U \twoheadrightarrow U$ restricts on this support $\Supp(\scrF) \subset A\times U$ to a finite morphism
\[ \gamma|_{\Supp(\scrF)}: \quad \Supp(\scrF) \;\twoheadrightarrow \; U \]
and~$\gamma_*(\scrF)$ is a locally free sheaf on $U$ (in~\cite[def.~1.2.1]{KraemerMicrolocalII} we also required $\gamma|_{\Supp(\scrF)}$ to be a flat morphism, but as one of the referees pointed out, this condition might depend on the scheme structure chosen on $\Supp(\scrF)$ and in any case it is not needed for the following constructions). The category $\VB(A, U)$ is not abelian, but it is still an {\em exact category} in the sense of~\cite{BuehlerExact} if we define its short exact sequences to be those coming from short exact sequences in $\Mcoh(\scrO_{A\times U}^\an)$. We view $\VB(A, U)$ as a tensor category for the product
\[
 \scrF_1 * \scrF_2 \;=\; \varpi_* \rho^*(\scrF_1 \boxtimes \scrF_2).
\]
For~$u\in U$ we obtain an exact category~$\VB(A, u)$ of~{\em germs of vector bundles} whose objects are the pairs $\alpha=(\scrF_\alpha, U_\alpha)$ of a Zariski open neighborhood $U_\alpha \subseteq U$ of $u$ and~$\scrF_\alpha \in \VB(A, U_\alpha)$, with morphisms
\[
 \Hom_{\VB(A,u)}(\alpha, \beta) \;=\; \lim_{\longrightarrow} \, \Hom_{\VB(A, U)}(\scrF_\alpha|_{A\times W}, \scrF_\beta|_{A\times W})
\]
where the limit runs over all Zariski open neighborhoods $W \subseteq U_\alpha \cap U_\beta$ of $u$. We equip these categories of germs with the convolution product induced by the one on $\VB(A, U)$. For $u\in U$ we get a commutative diagram of exact tensor categories
\[
\xymatrix@M=0.5em{
 \LS(A, U) \ar[r] \ar[d] & \LS(A, u) \ar[d] \\
 \VB(A, U) \ar[r] & \VB(A, u)
}
\]
where all arrows are {\em exact} tensor functors in the sense that they send exact sequences to exact sequences, and except for the bottom one they are faithful.

\subsection{Subgroups of multiplicative type and normalizers}  \label{sec:multiplicativeI}

An algebraic group $T$ is said to be {\em of multiplicative type} if it embeds into $\bbG_m^r$ for some $r\in \bbN$. In this case the group
$X = \CartierDual{T}$ of algebraic characters is a finitely generated abelian group,
and the weight space decomposition gives an equivalence of abelian tensor categories
\[ \Rep(T) \;\simeq \; \Vect_X(\bbC) \]
where the right hand side denotes the category of finite dimensional $X$-graded complex vector spaces. The group of multiplicative type is recovered from its characters as the Cartier dual $T = \CartierDual{X}$, and this sets up an antiequivalence between the categories of groups of multiplicative type and finitely generated abelian groups. 

\medskip

Homomorphisms from groups of multiplicative type to a linear algebraic group $G$ correspond to fiber functors $\omega: \Rep(G) \longrightarrow \Vect_X(\bbC)$, where by a {\em fiber functor} we mean any faithful exact tensor functor. In looking for such functors we do not require the abelian group $X$ to be finitely generated, but for the definition of $\omega$ it can be replaced by the subgroup
\[
 X(\omega) \;=\; \bigl\{ 
 x\in X \mid
 \omega(\rho)_x \neq 0 \; \textnormal{for some} \; \rho \in \Rep(G)
 \bigr\}
\]
of occuring weights, and this subgroup is finitely generated: Indeed, by the Tannakian characterization of algebraic groups~\cite[prop.~2.20b]{DM} the rigid abelian tensor category $\Rep(G)$ admits a tensor generator, more precisely it is generated by any faithful representation $\delta\in \Rep(G)$; so $X(\omega)$ is generated by the finite set $\{ x\in X \mid \omega(\delta)_x \neq 0 \}$ of indices that occur in the grading of the given representation. In this setup we denote by
\[
 T(\omega) \;=\; \CartierDual{X(\omega)} \;\hookrightarrow \; G
\]
the subgroup of multiplicative type defined by the grading on the fiber functor; note that by construction it embeds in $G$. For any linear algebraic group $N$ let us denote by $\Rep_X(N)$ the abelian tensor category of algebraic representations of the group $N$ on $X$-graded vector spaces
\[
 V \;=\; \bigoplus_{x\in X} \; V_x
\]
for which the group action permutes the graded pieces in the sense that for each $g\in G$ there exists a bijection $\sigma_g: X\to X$ with $g\cdot V_x \subseteq V_{\sigma_g(x)}$ for all $x\in X$. Then any diagram of tensor functors
\[
\xymatrix@C=1.8em@M=0.5em{
 \Rep(G) \ar[rr]^-\omega \ar@{..>}[dr]_-\exists && \Vect_X(\bbC) \\
 & \Rep_X(N) \ar[ur] &
}
\]
gives rise to a homomorphism
$
 N \rightarrow  N_G(T(\omega))
$
to the normalizer of the subgroup of multiplicative type constructed above.

\subsection{Back to Gauss maps}

Subgroups of multiplicative type arise from tensor functors to the categories of germs in section~\ref{sec:localsystems} as follows. As before let $\eta \in U$ be the generic point of the cotangent space to the abelian variety. If a germ $\alpha=(\scrF_\alpha, U_\alpha)$ in $\VB(A, \eta)$ or $\LS(A, \eta)$ is defined over a neighborhood of $u\in U(\bbC)$, then we denote by
\[
 X(\alpha, u) \;=\; \bigl\langle a \in A(\bbC) \mid (a,u) \in \Supp(\scrF_\alpha) \bigr\rangle \;\subset\; A(\bbC)
\]
the subgroup which is generated by the finitely many points of $\Supp(\scrF_\alpha)\cap \gamma_\alpha^{-1}(u)$ where~$ \gamma_\alpha: A\times U_\alpha \twoheadrightarrow U_\alpha $ is the projection. Thus $X(\alpha, u)$ is a finitely generated abelian group, and we denote its Cartier dual by
\[
 T(\alpha, u) \;=\; \CartierDual{X(\alpha, u)}.
\]
For $\alpha = (\scrF_\alpha, U_\alpha)\in \LS(A, u)$ it follows from the definition of $\LS(A, u)$ that  $\gamma_{\alpha*}(\scrF_\alpha)$ is a local system. On its fiber $ F_\alpha = \gamma_{\alpha *}(\scrF_\alpha)(u)$ we have the monodromy representation of $\pi_1(U_\alpha, u)$. We define the {\em algebraic monodromy group} of the germ $\alpha$ as the Zariski closure
 \[
  \Gal(\alpha, u) \;\subseteq\; \Gl(F_\alpha)
 \] 
 of the image of the monodromy representation
$
 \pi_1(U_\alpha, u) \rightarrow \Gl(F_\alpha)$.

\begin{thm} \label{thm:gradedfiberfunctor}
Let $G$ be a linear algebraic group and $\delta\in \Rep(G)$ a faithful representation. Then the following properties hold: \medskip
\begin{enumerate}
\item For $u\in U(\bbC)$, any exact tensor functor $F: \Rep(G) \rightarrow \VB(A, u)$ defines an embedding \medskip
\[
 T(\alpha, u) \; \hookrightarrow \; G
 \quad \textnormal{\em for the germ} \quad \alpha \;=\; F(\delta) \;\in \; \VB(A, u). \medskip
\]
\item For any exact tensor functor $F: \Rep(G) \rightarrow \LS(A, \eta)$ and very general $u\in U(\bbC)$ we have \smallskip
\[
 \Gal(\alpha, u) \;\hookrightarrow \; N_G(T(\alpha, u))
 \quad \textnormal{\em for} \quad \alpha \,=\, F(\delta) \,\in\, \LS(A, u) \,\subset\, \LS(A, \eta).
\]
\end{enumerate}
\end{thm}

{\em Proof.} (a) Taking the fiber of germs of vector bundles at the point $u$ gives a tensor functor
\[
 \omega: \quad
 \VB(A, u) \;\longrightarrow \; \Vect(\bbC), 
 \quad
 \alpha \;=\; (\scrF_\alpha, U_\alpha) \;\mapsto \; \gamma_{\alpha*}(\scrF_\alpha)(u).
\]
This is an exact functor because the functor which associates to a coherent vector bundle its fiber at a given point is exact, although not faithful. We claim that~$\omega$ admits a natural grading by the group $A(\bbC)$.
Indeed, for $\alpha = (\scrF_\alpha, U_\alpha) \in \VB(A, u)$ the support $\Supp(\scrF_\alpha)\subset A\times U_\alpha$ meets the fiber $A\times \{u\}$ only in finitely many points, hence if we write
\[
 i: \quad A \;\hookrightarrow\; A\times U_\alpha, \quad a \;\mapsto\; (a, u)
\]
for the embedding of that fiber and 
\[
\scrF_\alpha(a, u) \;:=\; \Gamma_{\{(a, u)\}}(A, i^*\scrF_\alpha)
\;\subseteq\; 
\Gamma(A, i^*\scrF_\alpha) \;=\; \omega(\alpha)
\]
for the group of sections supported in a given point $(a,u)$, we get a decomposition
\medskip
\begin{equation} \label{eq:grading-by-fiber}
 \omega(\alpha) 
 \;\;=\; \bigoplus_{a\in A(\bbC)} \scrF_\alpha(a,u) 
\end{equation}
where the sum is over the finitely many $a\in A(\bbC)$ with $(a,u) \in \Supp(\scrF_\alpha)$. One easily checks that this decomposition is compatible with tensor products: More precisely, for $\alpha, \beta \in \VB(A, u)$ the tensor functoriality gives a natural isomorphism \medskip
\[
  c_{\alpha, \beta}: \quad \omega(\alpha) \otimes \omega(\beta) \; \stackrel{\sim}{\longrightarrow} \;\omega(\alpha*\beta),
\]
and a look at the supports of the respective sheaves shows that $c_{\alpha, \beta}$ decomposes as a sum of isomorphisms \medskip
\[
 \bigoplus_{a+b = c} \scrF_\alpha(a,u) \otimes \scrF_\beta(b,u)
  \; \stackrel{\sim}{\longrightarrow} \;
  \varpi_*(\rho^{-1}(\scrF_\alpha \boxtimes \scrF_\beta))(c,u)
\]
for $c\in A(\bbC)$. The compatibility of the grading with the tensor product implies that for the objects in any given finitely generated tensor subcategory~$\langle \alpha \rangle \subset \VB(A, u)$, the grading only involves points in $X(\alpha, u) \subset A(\bbC)$. Applying this remark to the image of~$F: \Rep(G) \to \VB(A, u)$, we see that the tensor functor $\omega \circ F: \Rep(G)\to \Vect(\bbC)$ factors over a tensor functor to a category of graded vector spaces as indicated in the following diagram:
%
\[
\xymatrix{
\Rep(G) \ar[rr]^-{\omega \circ F} \ar[dr]_-{\exists}&& \Vect(\bbC) \\
& \Vect_X(\bbC) \ar[ur]
}
\]
Here the grading of the vector spaces is by elements of the group $X=X(\alpha, u)$, where one can take $\alpha = F(\delta)$ for any faithful $\delta \in \Rep(G)$; indeed, any faithful representation of an algebraic group generates the representation category of that group. To conclude the proof of (a), note that $\omega \circ F$ is a fiber functor because  any exact tensor functor between rigid abelian tensor categories with $\End(\one) = \bbC$ is faithful~\cite[prop.~1.19]{DM}. 

\medskip 

(b) It now seems tempting to replace $\VB(A, u)$ by $\LS(A, u)$ and try to show that for any $\alpha = (\scrF_\alpha, U_\alpha) \in \LS(A, u)$ the action of the monodromy group $\Gal(\alpha, u)$ on the fiber permutes the summands in~\eqref{eq:grading-by-fiber}. However, as one of the referees observed, this is not always true: If two components of $\Supp(\scrF_\alpha)$ intersect at a point $(a, u)$, a loop in $\pi_1(U_\alpha, u)$ can be lifted to each of these two components; if the endpoints of the two lifts differ, the monodromy along the loop cannot preserve the grading.

\medskip 

To avoid this problem, we start with a tensor functor to the category $\LS(A, \eta)$ of germs of local systems over a neighborhood of the generic point $\eta \in U$ rather than over a given closed point. Since $\Rep(G)$ is generated as a tensor category by any faithful $\delta \in \Rep(G)$, lemma~\ref{lem:specialization} says that any tensor functor $F: \Rep(G) \longrightarrow \LS(A, \eta)$ factors over a tensor functor
\[
 F: \quad \Rep(G) \;\longrightarrow\; \LS(A, u)
\]
for very general $u\in U(\bbC)$, and we can moreover choose this closed point $u$ in such a way that on the reduced closed subscheme which underlies the support of each of the germs $\alpha=(\scrF_\alpha, U_\alpha) \in F(\Rep(G))$, the Gauss map
$\gamma_\alpha: \Supp(\scrF_\alpha, U_\alpha)^\mathrm{red} \longrightarrow U_\alpha$ 
is \'etale over the chosen point. Then the problem described above disappears, and for any $\alpha \in F(\Rep(G))$ the natural monodromy operation of $\Gal(\alpha, u)$ on $\gamma_{\alpha,*}(\scrF_\alpha)(u)$ will permute the summands in~\eqref{eq:grading-by-fiber}. As before, the grading on the occuring vector spaces will be by points of the finitely generated subgroup $X=X(\alpha, u)\subset A(\bbC)$ where one can take $\alpha = F(\delta)$ for any faithful $\delta \in \Rep(G)$, and similarly the monodromy groups of all occuring local systems will be quotients of the monodromy group $N = \Gal(\alpha, u)$ for the chosen tensor generator. We therefore obtain a commutative diagram
\[
\xymatrix{
	\Rep(G) \ar[rr]^-{\omega \circ F} \ar[dr]_-{\exists}&& \Vect(\bbC) \\
	& \Rep_X(N) \ar[ur]
}
\]
where $\Rep_X(N)$ denotes the category of finite dimensional algebraic representations of $N$ whose underlying vector space is graded by $X$ such that the graded pieces are permuted by the group action. Hence we are in the setting of section~\ref{sec:multiplicativeI} and the claim follows.
\qed

\subsection{Dependence on the base point}

In the above the base point $u\in U(\bbC)$ has been fixed, and it is natural to ask how the subgroup of multiplicative type in theorem~\ref{thm:gradedfiberfunctor} depends on it.
More generally, for any algebraic subvariety $\Lambda \subset A\times U$ such that the projection $\gamma: \Lambda \twoheadrightarrow U$ is a finite morphism, consider the finitely generated abelian groups
\[
 X_u \;=\; \bigl \langle a\in A(\bbC) \mid (a, u) \in \Lambda \rangle
 \quad \textnormal{for} \quad u\;\in\;U(\bbC).
\]
For instance, if $\Lambda$ is the complement of the zero section in the conormal bundle to a smooth curve of genus two in its Jacobian, then $X_u$ is generated by a point on the curve, and is hence isomorphic to $\bbZ$ except for finitely many points $u\in U(\bbC)$ by Raynaud's theorem~\cite{RaynaudTorsion}. Returning to the general case, the following observation will be enough for our purpose:

\begin{lem} \label{lem:verygeneral}
The isomorphism type of $X_u$ is constant for very general $u\in U(\bbC)$.
\end{lem}

{\em Proof.} Shrinking $U$ we may assume $\gamma$ is a finite \'etale cover. Then the number $n$ of points in the fiber $F_u = \{ a\in A(\bbC) \mid (a, u) \in \Lambda \}$ will not depend on $u$, and locally in the classical topology we may identify all nearby fibers with each other. With these local identifications the subgroup 
$ R_u = \{ (c_a)_{a\in F_u}  \mid  {\textstyle\sum\nolimits_{a\in F_u}} c_a \cdot a = 0 \} \subseteq  \bbZ^n$
of relations between points of the fiber is constant for all $u$ outside countably many proper closed subvarieties: For $c=(c_1, \dots, c_n) \in \bbZ^n$ the corresponding  relation defines a closed subvariety 
$$
 S_c \;=\; \Bigl\{
 (a_i)_{1\leq i\leq n}\in A^n \mid \textstyle\sum_{i=1}^n c_i \cdot a_i \;=\; 0
 \Bigr\}
 \;\subseteq\;
 A^n \;=\; A\times \cdots \times A.
$$
Let $S'_c$ be the preimage of this subvariety in the $n$-fold fibered product $\Lambda \times_U \cdots \times_U \Lambda $
and denote by $S_c'' \subseteq S_c'$ the union of all those irreducible components whose image in~$U$ is a proper closed subset. Then the claim of the lemma will hold for all $u$ outside these countably many proper closed subsets. \qed

\section{Microlocalization} \label{sec:microlocalization}

We now construct an exact tensor functor from $\M(A)$ to the category $\LS(A, \eta)$ of germs of local systems. This will prove our main theorem~\ref{thm:gausstorus} by the specialization lemma~\ref{lem:specialization} and theorem~\ref{thm:gradedfiberfunctor}. Our tensor functor will arise from microlocal analysis and factors over a tensor category $\MM(A)$ of microdifferential modules.

\subsection{A reminder on microdifferential modules} \label{sec:microdifferential}

For a complex manifold $X$ let us denote by $\pi: T^*X\rightarrow X$ the projection from the total space of its cotangent bundle. On this total space we have the sheaf $\scrE_X$ of holomorphic microdifferential operators~\cite{SKK, KashiwaraSystems, SchapiraMicrodifferential}; recall that this is a sheaf of rings containing $\pi^{-1}(\scrD_X)$ and that the corresponding categories of right modules are related by the faithful exact functor
\[
 \Mod(\scrD_X) \;\longrightarrow \; \Mod(\scrE_X), \quad \scrM \; \mapsto \; \pi^{-1}(\scrM) \otimes_{\scrD_X} \scrE_X
\]
with
\[ \Char(\scrM) \;=\; \Supp \bigl(\pi^{-1}(\scrM) \otimes_{\scrD_X} \scrE_X \bigr). \]
Let $\Mhol(\scrE_X) \subset \Mod(\scrE_X)$ be the full abelian subcategory of all microdifferential modules which are {\em holonomic} in the sense that they are coherent and supported on a conic Lagrangian subvariety. We then get a functor $\Mhol(\scrD_X) \longrightarrow \Mhol(\scrE_X)$ and similarly a functor $\Dbhol(\scrD_X) \longrightarrow \Dbhol(\scrE_X)$ between the corresponding bounded derived categories. The direct image under a morphism $f: Y\rightarrow X$ is defined to be the functor
\[
 f_\dag: \quad \Dbhol(\scrE_Y) \;\longrightarrow \; \Dbhol(\scrE_X), \quad f_\dag(\scrM) \;=\; R\varpi_{*}(\rho^{-1}\scrM \otimes^L_{\rho^{-1} \scrE_Y} \scrE_{Y\to X}),
\]
where $\varpi = \varpi_f$ and~$\rho = \rho_f$ are the morphisms induced by $f$ and by its codifferential in the diagram
\[
T^* Y 
\;\stackrel{\rho_f}{\longleftarrow} \;
Y \times_X T^* X 
\;\stackrel{\varpi_f}{\longrightarrow}\; 
T^* X
\]
and $\scrE_{Y\to X}$ is the~$(\rho^{-1} \scrE_Y, \varpi^{-1}\scrE_X)$-bimodule on $Y\times_X T^*X$ defined in~\cite[def.~I.4.3.1]{SchapiraMicrodifferential}.

\begin{ex} \label{ex:projection}
If $f: Y = W\times X \twoheadrightarrow X$ is a projection, then its codifferential is the closed embedding
\[
 \rho_f: \quad  W\times T^* X \; \hookrightarrow \; T^*W \times T^*X \;=\; T^* Y
\]
arising from the inclusion of the zero section in the cotangent bundle to the fiber. We then have
\[
 \scrE_{Y\to X} \;\simeq\; \scrO_{W} \boxtimes \scrE_X
\]
as a sheaf on $T^*Y$ with support on the closed submanifold $W\times T^*X \subset T^* Y$ and with the obvious bimodule structure: This follows directly from the definition and from the compatibility with products in~\cite[prop.~II.3.3.6]{SchapiraMicrodifferential}, since for a projection the graph $W \times \Delta \subset W \times X\times X$ is a product of the fiber and the diagonal $\Delta \subset X\times X$.
\end{ex}

 Direct images are compatible with composition in the sense that for any other morphism $g: Z\rightarrow Y$ the diagram 
\[
\xymatrix@M=0.5em@C=0.6em@R=1.2em{
 && Z\times_X T^*X \ar[dl]_{\rho_{f, Z}} \ar[dr]^-{\varpi_{g, X}}  && \\
 & Z\times_Y T^*Y \ar[dl]_-{\rho_g} \ar[dr]^-{\varpi_g} && Y\times_X T^*X \ar[dl]_{\rho_f} \ar[dr]^-{\varpi_f} & \\
 T^*Z && T^*Y && T^*X
}
\]
induces a natural morphism
\[
 \rho_{f, Z}^{-1} (\scrE_{Z\to Y}) \otimes_{\scrE_Y} \varpi_{g,X}^{-1}(\scrE_{Y\to X}) \longrightarrow \scrE_{Z\to X}.
\] 
Outside the zero section this need not be an isomorphism, the algebraic tensor product on the left can be too small: This happens for instance when $\dim Z > 0, \dim X > 0$ and $Y = \{x\}$ for a point $x\in X$~\cite[chapt.~II, rem.~3.3.4]{SchapiraMicrodifferential}; the problem is that outside the zero section the symbols of microdifferential operators are holomorphic functions rather than polynomials. But if $f$ is smooth or $g$ is a closed embedding, then the above natural morphism is an isomorphism, and the underived tensor product on the left coincides with the derived tensor product as in~\cite[chapt.~II, lemma~3.5.1]{SKK}:

\begin{lem} \label{lem:bimodule}
If $f$ is smooth or $g$ is a closed embedding, then the above natural morphism is an isomorphism
\[
  \rho_{f,Z}^{-1}(\scrE_{Z\to Y})  \otimes_{\scrE_Y} \varpi_{g, X}^{-1}(\scrE_{Y\to X})
  \;\stackrel{\sim}{\longrightarrow}\;
  \scrE_{Z\to X}
\]
and we have
\[
\scrTor_i^{\scrE_Y}\Bigl(\rho_{f,Z}^{-1}(\scrE_{Z\to Y}), \, \varpi_{g, X}^{-1}(\scrE_{Y\to X})\Bigr)
\;=\; 0
\quad \textnormal{\em for all $i\neq 0$}. 
\]
\end{lem}

{\em Proof.} If $g$ is a closed embedding, this is shown in~\cite[prop.~II.3.3.3]{SchapiraMicrodifferential}. The case when~$f$ is smooth can be easily reduced to the former case; for lack of a reference we briefly sketch a proof. Since the claim is local, we can assume $f: Y = W \times X \to X$ is a projection. If we factor~$g$ over its graph, we are interested in a composite of a closed embedding with two projections
\[
\xymatrix@M=0.8em{
 Z \ar@{^{(}->}[r]^-i & Z\times Y \ar@{->>}[r]^-p & Y \;=\; W\times X \ar@{->>}[r]^-f & X.
}
\]
Here $i(z):=(z, g(z))$ and $p(z,y):=y$. Example~\ref{ex:projection} shows that
\[
 \scrE_{Z\times Y\to Y} \simeq \scrO_{Z} \boxtimes \scrE_Y
 \quad \textnormal{and} \quad 
 \scrE_{Y\to X} \simeq \scrO_{W} \boxtimes \scrE_X.
\]
Here $\scrE_{Z\times Y \to Y}$ is flat over $\scrE_Y$. Suppressing the sheaf-theoretic pullback from the notation, we get 
\begin{eqnarray*}
\scrE_{Z\times Y \to Y} \otimes_{\scrE_Y}^L \scrE_{Y\to X} 
&\;\simeq\; &
\scrE_{Z\times Y\to Y} \otimes_{\scrE_Y} \scrE_{Y\to X} \\
&\;\simeq\; &
(\scrO_{Z} \boxtimes \scrE_{Y})\otimes_{\scrE_{Y}} (\scrO_{W}\boxtimes \scrE_X) \\
&\;\simeq\; &
\scrO_{Z} \boxtimes \scrO_{W} \boxtimes \scrE_X \\
&\;\simeq\; & 
\scrE_{Z\times Y\to X}
\end{eqnarray*}
on $Z \times W \times T^* X \subset T^* (Z\times Y)$. On the other hand, since for the precomposition of any morphism with a closed embedding we already know the claim of the lemma, we also have isomorphisms
\begin{eqnarray*}
\scrE_{Z\to Z\times Y} \otimes_{\scrE_{Z\times Y}}^L \scrE_{Z\times Y\to \circled{?}} 
\;\stackrel{\sim}{\longrightarrow}\; 
\scrE_{Z\to Z\times Y} \otimes_{\scrE_{Z\times Y}} \scrE_{Z\times Y \to \circled{?}} 	
\;\stackrel{\sim}{\longrightarrow}\; 
\scrE_{Z\to \circled{?}}
\end{eqnarray*} 
for the projection to both $\circled{?} = Y$ and $\circled{?}=X$. Combining all these isomorphisms we obtain that
\begin{eqnarray*}
\scrE_{Z\to Y} \otimes_{\scrE_Y} \scrE_{Y\to X} 
&\;\simeq\; &
(\scrE_{Z\to Z\times Y}\otimes_{\scrE_{Z\times Y}} \scrE_{Z\times Y \to Y}) \otimes_{\scrE_Y} \scrE_{Y\to X} \\
&\;\simeq\;& 
\scrE_{Z\to Z\times Y}\otimes_{\scrE_{Z\times Y}}( \scrE_{Z\times Y \to Y} \otimes_{\scrE_Y} \scrE_{Y\to X} ) \\
&\;\simeq\;& 
\scrE_{Z\to Z\times Y}\otimes_{\scrE_{Z\times Y}} \scrE_{Z\times Y \to X} \\
&\;\simeq\;&
\scrE_{Z\to X}
\end{eqnarray*} 
where each of the occuring tensor products of bimodules is quasi-isomorphic to their left derived tensor product.
\qed 

\medskip 

\begin{cor} \label{cor:composition-of-direct-images}
If in the situation of the above lemma the morphisms $f$ and $g$ are proper, we obtain for any $\scrM\in \Dbhol(\scrE_Z)$ natural isomorphisms
\begin{align} \nonumber
 f_\dag (g_\dag (\scrM))
 & \;=\;  R\varpi_{f,*} \bigl( \rho_f^{-1} \bigl(R\varpi_{g,*}(\rho_g^{-1} \scrM \otimes^L_{\scrE_Z} \scrE_{Z\to Y})\bigr) \otimes^L_{\scrE_Y} \scrE_{Y\to X} \bigr)\\ \nonumber
 & \;\simeq \; R\varpi_{fg, *} \bigl( \rho_{fg}^{-1}(\scrM) \otimes^L_{\scrE_Z}  
 		\bigl(\rho_{f, Z}^{-1} (\scrE_{Z\to Y}) \otimes_{\scrE_Y} \varpi_{g,X}^{-1}(\scrE_{Y\to X}) \bigr) \bigr) \\ \nonumber
 & \;\simeq \; R\varpi_{fg, *} \bigl( \rho_{fg}^{-1}(\scrM) \otimes^L_{\scrE_Z} \scrE_{Z\to X} \bigr) \\ \nonumber
 & \;=\; (f g)_\dag(\scrM) 
\end{align}
\end{cor} 

{\em Proof.} The first of the two isomorphisms involves only the projection formula and base change for proper direct images of sheaves, the second is induced by the natural isomorphism for transfer bimodules in lemma~\ref{lem:bimodule}. \qed 

\medskip 

The above isomorphisms will be used to define the commutativity and associativity constraints in proposition~\ref{prop:microtensor}, though for that purpose we only need the special case of isomorphisms and projections. We also remark that for $\scrM = \pi^{-1}(\scrN) \otimes_{\pi^{-1}(\scrD_Z)} \scrE_Z$ with $\scrN \in \Dbhol(\scrD_Z)$, the isomorphism in the above corollary could also be obtained from the analogous statement for the composition of direct images of $\scrD$-modules, using that microlocalization commutes with proper direct images~\cite[th.~7.5]{SSIndex}.

\subsection{Microlocal convolution} \label{sec:microconvolution}

We now specialize to the situation where $X=A$ is an abelian variety, in which case the cotangent bundle $T^*A$ is the trivial bundle with fiber~$V=H^0(A, \Omega^1)$. Let $U\subseteq V$ be a non-empty open subset. For $\scrM \in \Mhol(\scrE_A|_{A\times U})$ we have as in section~\ref{sec:localsystems} the Gauss map
\[
 \gamma: \quad \Supp(\scrM) \;\subset \; A\times U \; \twoheadrightarrow \; U
\]
which is generically finite. Let $d\in \bbN_0$ be its degree. If $d=0$, then we say that~$\scrM$ is~{\em negligible}. The negligible modules form a Serre subcategory of $\Mhol(\scrE_A|_{A\times U})$, and in what follows we denote by $\MM(A, U)$ the abelian quotient category; the notation~$\mathrm{MM}$ here stands for {\em m}icrodifferential {\em m}odules. Similarly, we denote by $\DM(A, U)$ the Verdier quotient of the triangulated category $\Dbhol(\scrE_A|_{A\times U})$ by the thick subcategory of complexes whose cohomology sheaves are negligible. In the case when~$U=V$ we simply write 
\[ \MM(A):=\MM(A, V) \quad \textnormal{and} \quad  \DM(A):=\DM(A, V).
\]
The various categories of differential and microdifferential modules from above fit in the following commutative diagram, where the dotted arrows come from the universal property of the quotient categories~$\M(A)$ and $\D(A)$ in section~\ref{sec:tannakiansetting}: \medskip
\[
\xymatrix@R=1em@C=0em@M=0.5em{
 & \Dbhol(\scrD_A) \ar[rr] \ar@{-}[d] && \D(A) \ar@{..>}[dd] \\
 \Mhol(\scrD_A)  \ar[dd] \ar[rr] \ar@{^{(}->}[ur] & \ar[d] & \M(A) \ar@{..>}[dd] \ar@{^{(}->}[ur] \\
 & \Dbhol(\scrE_A) \ar@{-}[r] \ar@{-}[d] & \ar[r] & \DM(A) \ar[dd] \\
 \Mhol(\scrE_A) \ar[rr] \ar[dd] \ar@{^{(}->}[ur] & \ar[d] & \MM(A) \ar@{^{(}->}[ur] \ar[dd] \\
  & \Dbhol(\scrE_A|_{A\times U}) \ar@{-}[r] & \ar[r] & \DM(A, U) \\
 \Mhol(\scrE_A|_{A\times U}) \ar[rr] \ar@{^{(}->}[ur] && \MM(A, U) \ar@{^{(}->}[ur]
} \medskip
\]
For $\scrM_1, \scrM_2 \in \Dbhol(\scrE_A)$ we define 
\[
 \scrM_1 * \scrM_2 \;:=\; a_\dag(M_1\boxtimes M_2). \]
More generally we can define the convolution of $\scrM_1, \scrM_2 \in \Dbhol(\scrE_A|_{A\times U})$ via the following diagram
\[
 \xymatrix@M=0.8em{
 & A^2 \times U \ar[dl]_-{\rho_U} \ar[dr]^-{\varpi_U} \ar@{^{(}->}[d] & \\
 A^2 \times U^2 \ar@{^{(}->}[d] & A^2 \times V \ar[dl]_-{\rho} \ar[dr]^-{\varpi} & A\times U \ar@{^{(}->}[d] \\
 A^2 \times V^2 && A\times V
}
\]
where $\rho=\rho_a$ and $\varpi=\varpi_a$ are induced by the group law on the abelian variety. We put
\[
 \scrM_1 * \scrM_2 \;:=\; R\varpi_{U,*}\bigl(\rho_U^{-1}(\scrM_1 \boxtimes \scrM_2 \bigr) \otimes^L_{\rho_U^{-1} (\scrE_{A^2}|_{A^2 \times U})} (\scrE_{A^2 \to A})|_{A\times U}).
\]
For $U=V$ this reduces to our previous global definition of the convolution product on the categories $\Dbhol(\scrE_A|_{A\times V})=\Dbhol(\scrE_A)$ and $\DM(A, V)=\DM(A)$, but since our later constructions will work over open dense subsets of the cotangent space, it seems natural to include this more general case here.

\begin{propo} \label{prop:microtensor}
The above convolution product endows $\Dbhol(\scrE_A|_{A\times U})$ and its Verdier quotient $\DM(A, U)$ with the structure of a triangulated tensor category in a natural way such that all the functors in the following commutative diagram are tensor functors:
%
\[
\xymatrix@M=0.8em{
\Dbhol(\scrD_A) \ar[r] \ar[d] & \Dbhol(\scrE_A) \ar[r] \ar[d] & \Dbhol(\scrE_A|_{A\times U}) \ar[d] \\
\D(A) \ar[r] & \DM(A) \ar[r] & \DM(A, U)
}
\]
\end{propo}

{\em Proof.} Since all our constructions will be compatible with the restriction to open subsets of the form $A\times U\subseteq A\times V$ in the cotangent bundle, we only discuss the global case.  The main point is to endow the triangulated category $\Dbhol(\scrE_A)$ with a tensor structure with respect to the convolution product. Clearly the convolution product is a $\bbC$-linear triangulated bifunctor. If $i: \{0\} \hookrightarrow A$ denotes the inclusion of the origin of the abelian variety, the unit object for the tensor structure will be~$\one = i_\dag(\bbC)$. Using the compatibility of direct images with composition, we equip it  with the natural isomorphism 
\[ 
 u: \quad \one \;=\; i_\dag(\bbC) \;\stackrel{\sim}{\longrightarrow} \; a_\dag (i\times i)_\dag (\bbC \boxtimes \bbC) \; = \; \one * \one
\]  
which is induced by the identitification $i = a\circ (i\times i)$. To define the commutativity and associativity constraints and to verify the required compatibilities between them, let~$\scrM_1, \dots, \scrM_4 \in \Dbhol(\scrE_A)$, and for arbitrary $n\in \bbN$ denote by $a: A^n \rightarrow A$ the  addition morphism. The commutativity constraint $\psi_{\scrM_1, \scrM_2}$ is given by the natural isomorphisms in the diagram
\[
\xymatrix@C=3em@M=0.5em{
 \scrM_1 * \scrM_2 \ar@{..>}[rr]^-{\psi_{\scrM_1, \scrM_2}} \ar@{=}[d]
 &&  \scrM_2 * \scrM_1  \ar@{=}[d] \\
 a_\dag(\scrM_1\boxtimes \scrM_2) \ar[r]^-\sim
 & a_\dag \sigma_\dag(\scrM_2\boxtimes \scrM_1) \ar[r]^-\sim
 & a_\dag (\scrM_2\boxtimes \scrM_1)  \\
}
\]
where $\sigma: A\times A \longrightarrow A\times A, (x,y)\mapsto (y,x)$ is the involution that interchanges the two factors. Similarly, the associativity constraint $\varphi_{\scrM_1, \scrM_2, \scrM_3}$ is given by the natural isomorphisms in the following diagram:
\[
\xymatrix@C=1.8em@M=0.5em{
(\scrM_1 * \scrM_2) * \scrM_3 \ar@{..>}[d]^-{\varphi_{\scrM_1, \scrM_2, \scrM_3}} \ar@{=}[r] 
& a_\dag( (\scrM_1 * \scrM_2) \boxtimes \scrM_3 )\ar[r]
& a_{\dag} ((\scrM_1 \boxtimes \scrM_2) \boxtimes \scrM_3)\ar@{=}[d] \\
\scrM_1 * (\scrM_2 * \scrM_3) 
& a_\dag( \scrM_1 \boxtimes (\scrM_2* \scrM_3)) \ar@{=}[l] 
& a_{\dag} (\scrM_1 \boxtimes (\scrM_2 \boxtimes \scrM_3)) \ar[l]
}
\]
It remains to check that these commutativity and associativity constraints satisfy the pentagon and hexagon axioms. Putting $\scrM = \scrM_1 \boxtimes \cdots \boxtimes \scrM_4$, the pentagon axiom boils down to the commutativity of the five triangles in the following diagram of natural isomorphisms:
\[
\footnotesize
\xymatrix@C=0em@M=0.5em@R=2em{
 && (\scrM_1 * (\scrM_2 * \scrM_3)) * \scrM_4 \ar[dll] & \\
 \hbox to 7em{$\scrM_1 * ((\scrM_2 * \scrM_3) * \scrM_4)$} \ar@<3ex>[dd]  &&& \\
 && a_{\dag}(\scrM) \ar@{..>}[uu] \ar@{..>}[ull] \ar@{..>}[r] \ar@{..>}[dd] \ar@{..>}[dll] & \hbox to 1em{$((\scrM_1 * \scrM_2) * \scrM_3) * \scrM_4$}\ar[uul] \ar[ddl] \\
 \hbox to 7em{$\scrM_1 * (\scrM_2 * (\scrM_3 * \scrM_4))$} &&& \\
 && (\scrM_1 * \scrM_2) * (\scrM_3 * \scrM_4) \ar[ull] &
}
\]
Similarly, putting $\scrM_{ijk} = \scrM_i \boxtimes \scrM_j \boxtimes \scrM_k$ for $i,j,k\in \{1,2,3\}$, the hexagon axiom boils down to the commutativity of the three trapezoids in the following diagram of natural isomorphisms:
\[
\footnotesize
\xymatrix@C=-0.4em@R=1em@M=0.4em{
 && (\scrM_2 * \scrM_3) * \scrM_1 \ar[ddll]  && \\
 &&&  a_\dag(\scrM_{231}) \ar[dddlll] \ar[dddddd] \ar@{..>}[ul] \ar@{..>}[dr] & \\
 \scrM_1 *(\scrM_2 * \scrM_3) &&&& \scrM_2 * (\scrM_3 * \scrM_1) \\
 &&&& \\
 a_\dag(\scrM_{123}) \ar[dddrrr] \ar@{..>}[uu] \ar@{..>}[dd] &&&& \\
 &&&& \\
 (\scrM_1 * \scrM_2) * \scrM_3 \ar[ddrr] &&&& \scrM_2 * (\scrM_1 * \scrM_3) \ar[uuuu] \\
 &&& a_\dag(\scrM_{213}) \ar@{..>}[ur] \ar@{..>}[dl] & \\
 && (\scrM_2 * \scrM_1) * \scrM_3  &&
}
\]
Going back to the construction of the isomorphism for the composition of direct images in corollary~\ref{cor:composition-of-direct-images}, one now deduces the pentagon and hexagon axioms from compatibilities for  natural maps between tensor products of transfer bimodules. It follows that $\Dbhol(\scrE_A)$ is a triangulated tensor category with respect to the convolution product. Microlocalization commutes with proper direct images~\cite[th.~7.5]{SSIndex}, so we have natural isomorphisms 
\[
 (\scrN_1 * \scrN_2) \otimes_{\scrD_A} \scrE_A \; \stackrel{\sim}{\longrightarrow} \; (\scrN_1 \otimes_{\scrD_A} \scrE_A) * (\scrN_2 \otimes_{\scrD_A} \scrE_A)
 \quad
 \textnormal{for}
 \quad 
 \scrN_1, \scrN_2 \in \Dbhol(\scrD_A),
\]
and these are compatible with our associativity and commutativity constraints so that they make $\Dbhol(\scrD_A) \longrightarrow \Dbhol(\scrE_A)$ a tensor functor. Finally, the convolution product descends to the quotient category $\DM(A)$, and we endow the latter with the induced tensor structure so that the quotient functor is a tensor functor. \qed

\begin{cor}
The convolution product endows $\MM(A)$ with a natural structure of an abelian tensor category such that the quotient functor $\M(A) \longrightarrow \MM(A)$ is a faithful exact tensor functor.
\end{cor}

{\em Proof.}
We claim that the essential image of $\MM(A)$ in $\DM(A)$ is stable under the convolution product. For this it suffices to show that for any $\scrM_1, \scrM_2 \in \Mhol(\scrE_A)$, the truncation morphisms
\begin{equation} \label{eq:truncation}
 \scrM_1 * \scrM_2 \;\longrightarrow \; \tau_{\geq 0}(\scrM_1 * \scrM_2) \;\longleftarrow \; \scrH^0(\scrM_1 * \scrM_2)
\end{equation}
restrict to isomorphisms over some Zariski open dense subset. This holds over any Zariski open dense subset $U \subset V$ over which both $\Lambda_i = \Supp(\scrM_i)$ are finite: Indeed, then
\[
 \varpi: \quad \rho^{-1}(\Supp(\scrM_1 \boxtimes \scrM_2)) \cap \varpi^{-1}(U) \;=\; \Lambda_1 \times_U \Lambda_2\;\longrightarrow \; U
\]
is a finite morphism, i.e.~the addition morphism $a:A\times A \longrightarrow A$ is  non-characteristic for $\scrM_1 \boxtimes \scrM_2$ over the subset $A\times U \subset T^*A$ in the sense of~\cite[def.~II.3.1.2(b)]{SchapiraMicrodifferential}; but in the non-characteristic case theorem~3.4.4 and remark~3.1.7 in loc.~cit.~say that the derived pushforward functor for microdifferential modules coincides with the naive pushforward, hence the morphisms in~\eqref{eq:truncation} restrict over $U$ to isomorphisms.
\qed

\medskip

So if we can find a faithful exact tensor functor $\omega_\eta: \MM(A) \longrightarrow \LS(A, \eta)$ to the  category of germs from section~\ref{sec:localsystems}, we will get a fiber functor on any finitely generated tensor subcategory of $\M(A)$. For this we now introduce a reference system of simple holonomic microdifferential modules on abelian varieties.

\subsection{A class of simple test modules} \label{sec:simplemodules}

Let $\Omega\subseteq T^*A$ be an open subset. For a sheaf $\scrI \trianglelefteq \scrE_A|_\Omega$ of right ideals, one defines the {\em symbol ideal} $\sigma(\scrI) \trianglelefteq \scrO_\Omega$ as the ideal generated by the principal symbols of all microdifferential operators in~$\scrI$. Let us say that $\scrI$ is {\em reduced} if its symbol ideal $\sigma(\scrI)$ is the full ideal of functions vanishing on some closed subvariety of $\Omega$. A module $\scrM \in \Mod(\scrE_A|_\Omega)$ will be called {\em simple} if it can be written as a quotient of $\scrE_A|_\Omega$ by a reduced right ideal. Working on the preimages $\Omega = A\times U$ of suitable open subsets $U\subseteq V$ of the cotangent space, we want to attach to every conic Lagrangian subvariety a simple module in a natural way.

\medskip

To achieve this, recall that for any closed conic Lagrangian subvariety $\Lambda \subset T^*A$ there exists a Zariski open dense subset $U\subseteq V$ over which $\Lambda$ is finite \'etale in the sense that the projection 
\[
  \Lambda |_U \;=\; \Lambda\times_V U \;\subset\; T^*A |_U \;=\; A\times U \;\twoheadrightarrow \; U
\]
is a finite \'etale cover (possibly empty). For any such open subset we define a sheaf of right ideals
$$ \scrI \;=\; \scrI_{\Lambda, U} \; \trianglelefteq \; \scrE_{A, U} \; = \; \scrE_A|_\Omega
\quad \textnormal{on} \quad \Omega \;=\; A\times U$$
as follows. Let us cover $\Omega=A\times U$ by coordinate charts of the form $\Omega_0=A_0 \times U_0$ where~$A_0 \subset A$ and $U_0\subset U$ are open in the classical topology. Choosing the open subsets to be sufficiently small, we may assume that for each coordinate chart $\Omega_0$ one of the following two cases occurs:

\smallskip

\begin{itemize}
\item Either $\Omega_0\cap \Lambda = \varnothing$. In this case we put $\scrI|_{\Omega_0} = \scrE_A|_{\Omega_0}$. \smallskip
\item Or $\Omega_0\cap \Lambda$ projects isomorphically onto $U_0$. In this case, let $z=(z_1, \dots, z_g)$ be local coordinates on $A_0 \subset A$ which pull back to affine linear coordinates on the universal cover of the abelian variety $A$, and let $f_i: U_0 \rightarrow \bbC$ be the holomorphic functions defined by
$$
\quad \Omega_0 \cap \Lambda \;=\; \bigl\{ (z, \xi) \in A_0\times U_0 \mid z_i = f_i(\xi) \; \textnormal{for} \; i=1,2,\dots, g \bigr\}.
$$
Then each $f_i$ is homogenous of degree zero since $\Lambda$ is conic. So in our chosen coordinate system we may view $z_i - f_i(\xi)$ as a section of $\scrE_A(0)|_{\Omega_0}$, and we define
$$
  \scrI|_{\Omega_0} \;=\; \sum_{i=1}^g \; (z_i - f_i(\xi))\cdot \scrE_A|_{\Omega_0}.
$$
\end{itemize} 

\noindent
The formula for the transformation of microdifferential operators under coordinate changes~\cite[sect.~7.2, equation (7.8)]{KashiwaraDModules}, applied to the special case of affine linear coordinate changes where no higher derivatives occur, shows that the definition does not depend on the choice of the affine linear coordinates. Hence by patching we obtain a sheaf of right ideals $\scrI = \scrI_{\Lambda, U} \trianglelefteq \scrE_{A,U}$.  

\begin{lem}
In the above setting, the symbol ideal $\sigma(\scrI_{\Lambda, U}) \trianglelefteq \scrO_{\Omega}$ coincides with the ideal of all functions vanishing on $\Omega \cap \Lambda$. Hence the quotient by it defines a simple module
\[
 \scrC_{\Lambda, U} \;=\; \scrI_{\Lambda, U} \backslash \scrE_{A, U}
 \;\in\; \Mhol(\scrE_A|_\Omega)
 \quad
 \textnormal{\em with}
 \quad
 \Supp(\scrC_{\Lambda, U}) \;=\; \Omega \cap \Lambda.
\]
\end{lem}

{\em Proof.} We first verify the claim about the symbol ideal. Working locally on some chart $\Omega_0 = A_0 \times U_0$ as above, we may assume that the right ideal $\scrI_{\Lambda, U}|_{\Omega_0}\trianglelefteq \scrE_{A}|_{\Omega_0}$ is generated by the homogenous functions $P_i(z, \xi) = z_i - f_i(\xi)$, where the latter are considered as sections of $\scrE_A(0)|_{\Omega_0}$. Since these generators depend only linearly on the variables $z_i$, the Leibniz rule for the product of microdifferential operators shows that the commutators $[P_i, P_k] \in \scrE_A(-1)|_{\Omega_0}$ coincide in our linear coordinates with their principal symbol, which is given by the Poisson bracket \medskip
\[
 \bigl[ P_i, P_k \bigr] 
 \;=\; \bigl\{ P_i, P_k \bigr\}
 \;=\; \sum_j \; 
 \left( \frac{\partial P_i}{\partial \xi_j}\frac{\partial P_k}{\partial z_j} - \frac{\partial P_k}{\partial \xi_j}\frac{\partial P_i}{\partial z_j}
 \right)
 \;=\; \frac{\partial f_k}{\partial \xi_i} - \frac{\partial f_i}{\partial \xi_k}. 
\]
The function on the right hand side depends only on the variables $\xi_j$ but not on the variables $z_j$. On the other hand, since $\Lambda_0 = \Omega_0 \cap \Lambda$ is an isotropic subvariety, the ideal of functions vanishing on it is stable under the Poisson bracket. Thus on the right hand side of the above equation we have a function on $\Omega_0 = A_0\times U_0$ which only depends on $U_0$ but nevertheless vanishes on $\Lambda_0$. This is possible only for the zero function since by assumption $\Lambda_0$ is finite \'etale over $U_0$. Thus the generators $P_1, \dots, P_g$ for the right ideal~$\scrI_{\Lambda, U}$ commute with each other; moreover, the differentials of their principal symbols  are clearly linearly independent. By~\cite[prop.~I.4.1.5]{SchapiraMicrodifferential} it then follows that the principal symbols $\sigma(P_1), \dots, \sigma(P_g)$ already generate the entire symbol ideal~$\sigma(\scrI_{\Lambda, U})\unlhd \scrO_\Omega$. In fact our chosen generators satisfy $\sigma(P_i)=P_i$ since the $f_i(\xi)$ are homogenous functions of degree zero; the upshot of this discussion is that the symbol ideal~$\sigma(\scrI_{\Lambda, U})$ is generated by the functions $z_i - f_i(\xi)$, hence it is the ideal of all functions vanishing on $\Omega\cap \Lambda$ as claimed. 

\medskip 

In particular, the quotient $\scrC_{\Lambda, U}  = \scrI_{\Lambda, U} \backslash \scrE_{A, U}\in \Mod(\scrE_{A, U})$ is a simple module in the sense defined above. For the support of this module we get $\Supp(\scrC_{\Lambda, U}) = \Omega \cap \Lambda$, which is a Lagrangian subvariety of $\Omega$. Hence $\scrC_{\Lambda, U}$ is holonomic.  \qed

\medskip

Note that the simple module $\scrC_{\Lambda, U}$ comes with a canonical section on $\Omega = A\times U$, the class of $1\in \scrE_{A,U}$. We will use these sections to establish a compatibility property for convolutions of our simple modules. Let $\Lambda_1, \Lambda_2 \subset T^*A = A\times V$ be conic Lagrangian subvarieties, and put
$$
 \Lambda_1 * \Lambda_2 \;=\; \varpi(\rho^{-1}(\Lambda_1 \times \Lambda_2)) \;\subset \; A\times V
$$
as in section~\ref{sec:localsystems}. Let $U\subseteq V$ be any open subset over which $\Lambda_1, \Lambda_2$ and $\Lambda_1*\Lambda_2$ are finite \'etale so that the corresponding simple modules exist over this subset. 

\begin{lem} \label{lem:cotrace}
We have a natural homomorphism 
\[
 \iota_{\Lambda_1, \Lambda_2}: \quad \scrC_{\Lambda_1*\Lambda_2, U} \; \longrightarrow \; 
\scrC_{\Lambda_1, U} * \scrC_{\Lambda_2, U}.
\]
\end{lem}

{\em Proof.} Put $\Lambda = \Lambda_1 * \Lambda_2$. The homomorphism will be defined by the commutative diagram
\[
\xymatrix@C=2em@R=1.5em@M=0.5em{
\scrE_A|_{A\times U} \ar[rr]^-{P\, \mapsto \, u\cdot P} \ar@{->>}[dr]_-{P\, \mapsto \, 1\cdot P \;\;\; } && 
\scrC_{\Lambda_1, U} * \scrC_{\Lambda_2, U} \\
& \scrC_{\Lambda, U} \ar@{..>}[ur]_-{\exists !}
}
\]
where
$$u \;=\; \varpi_*(\rho^{-1}(1\boxtimes 1) \otimes 1_{A^2\to A}) \;\in\; H^0(A\times U, \scrC_{\Lambda_1, U} * \scrC_{\Lambda_2, U}). $$ 
To see that such a factorization exists, we must show that the ideal sheaf $\scrI_{\Lambda, U}$ is contained in the annihilator of $u$. This is a local question, so it suffices to show it on~$\Omega_0 = A\times U_0$ for small open subsets $U_0 \subset U$ in the classical topology.
We may assume
\[
 \Omega_0 \cap \Lambda_1 \;=\; \bigsqcup_{\alpha} \; \Lambda_{1,\alpha}
 \quad \textnormal{and} \quad
 \Omega_0 \cap \Lambda_2 \;=\; \bigsqcup_{\beta} \; \Lambda_{2,\beta},
\]
where
\[
 \Lambda_{1,\alpha} \;=\; \bigl\{ (g_{1,\alpha}(\xi), \xi) \mid \xi \in U_0 \bigr\}
 \quad \textnormal{and} \quad
 \Lambda_{2,\beta} \;=\; \bigl\{ (g_{2,\beta}(\xi), \xi)  \mid \xi \in U_0 \bigr\} 
\]
are the graphs of holomorphic maps 
$g_{1,\alpha}: U_0  \rightarrow A$ and $g_{2,\beta}: U_0  \rightarrow A$.
Note that the decomposition of the supports into disjoint closed subsets gives decompositions as direct sums
\begin{align} \nonumber
 \scrC_{\Lambda_1, U_0} \;=\; \bigoplus_\alpha \; \scrC_{1,\alpha} \quad
 &\textnormal{with} \quad  \Supp(\scrC_{1,\alpha}) = \Lambda_{1,\alpha}, \\ \nonumber
  \scrC_{\Lambda_2, U_0} \;=\; \bigoplus_\alpha \; \scrC_{2,\beta} \quad
  &\textnormal{with} \quad  \Supp(\scrC_{2,\beta}) = \Lambda_{2,\beta}. 
\end{align}
Let $1_\alpha \in \scrC_{1,\alpha}$ and $1_\beta \in \scrC_{2, \beta}$ be the components of $1\in \scrC_{\Lambda_1, U_0}$ and $1\in \scrC_{\Lambda_2, U_0}$ in the respective direct summands. We need to verify that the right ideal~$\scrI_{\Lambda, U_0}$ annihilates the section 
\[
u_{\alpha, \beta} \;=\; \varpi_*(\rho^{-1}(1_\alpha \boxtimes 1_\beta) \otimes 1_{A^2 \to A} ) \;\in \; \scrC_{1,\alpha} * \scrC_{2,\beta}
\]
for all $\alpha,\beta$. To this end, note that locally near any point of $\varpi(\rho^{-1}(\Lambda_{1,\alpha} \times \Lambda_{2,\beta}))$, we can write $\Lambda$ as the graph of $f=g_{1,\alpha} + g_{2,\beta}: U_0 \rightarrow A$. Hence locally near any such point we have
\[
 \scrI_{\Lambda, U_0} \;=\; \sum_{i=1}^g \; (z_i - f_{i}(\xi))\cdot \scrE_A|_{\Omega_0}
\]
where $z_i$ is the $i$-th coordinate in a local coordinate system $(z,\xi)$ on $T^*A=A\times V$ that arises from a suitable affine linear coordinate system on the universal cover and where $f_{i}$ denotes the $i$-th component of $f$. For a suitable local affine linear coordinate system $(z_1, z_2, \xi_1, \xi_2)$ on the product $T^*A^2 = A^2 \times V^2$, we have in $\scrE_{A^2 \to A}$ the relations
\begin{eqnarray} \nonumber
 1_{A^2 \to A}\cdot z_i
 &=& 
 (z_{1,i} + z_{2,i}) \cdot 1_{A^2 \to A},  \\ \nonumber
 1_{A^2 \to A}\cdot h(\xi)  
 &=& 
 h(\xi_1) \cdot 1_{A^2 \to A} \;\;=\;\; h(\xi_2) \cdot 1_{A^2 \to A}
\end{eqnarray}
for any homogenous holomorphic function $h(\xi)$, considered as a microdifferential operator. Thus locally
\[
1_{A^2 \to A} \cdot (z_i - f_i(\xi)) \;=\;
(z_{1,i} - g_{1,\alpha,i}(\xi_1)  + z_{2,i} - g_{2,\beta,i}(\xi_2)) \cdot 1_{A^2 \to A},
\]
and it follows that $u_{\alpha, \beta} \cdot (z_i - f_i(\xi)) = 0$ as required. \qed

\subsection{Second microlocalization} \label{sec:secondmicro}

We now pass from holonomic microdifferential modules to local systems on the smooth locus of their support, which will yield a variant of the second microlocalization compatible with convolution. To define the required tensor functor we apply a change of coefficients $\scrM \mapsto \scrM^\bbR = \scrM\otimes_{\scrE_A} \scrE_A^\bbR$ via the ring extension $\scrE_A \subset \scrE_A^\bbR$. For the definition of the latter we refer the reader to~\cite[sect.~1.4]{KashiwaraSystems} and only quote the following result:

\begin{thm} \label{thm:microclassification}
Let $\Omega \subset T^*A$ be an open subset, and let $\scrC_\Lambda \in \Mod(\scrE_A|_\Omega)$ be simple with support a smooth Lagrangian subvariety $\Lambda \subset \Omega$. Then any coherent module $\scrM \in \Mod(\scrE_A|_\Omega)$ supported on $\Lambda$ satisfies $\scrExt_{\scrE_A}^i(\scrC_\Lambda, \scrM^\bbR)=0$ for all $i\neq 0$, and
\[ \scrHom_{\scrE_A}(\scrC_\Lambda, \scrM^\bbR) \]
is a local system on $\Lambda$ whose rank is equal to the multiplicity of $\scrM$ along $\Lambda$. 
\end{thm}

{\em Proof.} This is shown in \cite[th.~3.2.1]{KashiwaraSystems}. \qed 

\medskip

We will apply this result for the simple modules $\scrC_\Lambda = \scrC_{\Lambda, U}$ from the previous section, omitting the open subset $U \subseteq V$ from the notation when there is no risk of confusion. By the same abuse of notation we will also omit the open subset for objects of the category $\LS(A, \eta)$ of germs from section~\ref{sec:localsystems}. We are interested in the functor
\[
 F: \quad  \MM(A) \; \longrightarrow \; \LS(A, \eta),
 \quad
 \scrM \;\mapsto \;  \scrHom_{\scrE_A}(\scrC_{\Lambda}, \scrM^\bbR)|_{U}
\]
where on the right hand side $U \subseteq V$ is taken to be the maximal Zariski open over which $\Lambda = \Supp(\scrM)$ is finite \'etale. The functor $F$ is faithful and exact, and our choice of the simple modules in the previous section implies:

\begin{thm} \label{thm:microlocal-fiberfunctor}
The functor $F: \MM(A) \longrightarrow \LS(A, \eta)$ underlies a tensor functor. 
\end{thm}

{\em Proof.} 
We must find isomorphisms 
\[
 c_{\scrM_1, \scrM_2}: \quad
 F(\scrM_1) * F(\scrM_2) \; \stackrel{\sim}{\longrightarrow} \; F(\scrM_1 * \scrM_2)
 \quad 
 \textnormal{for}
 \quad
 \scrM_1, \scrM_2 \;\in\; \MM(A)
\]
which are functorial and compatible with the associativity, commutativity and unit constraints. Put $\Lambda_i = \Supp(\scrM_i)$, and let $U\subset V$ be the maximal open subset where both $\scrM_i$ are defined and over which the Gauss maps for $\Lambda_1, \Lambda_2$ and $\Lambda = \Lambda_1 * \Lambda_2$ are finite \'etale. Writing 
\[
 \scrM_{12}^\bbR \;=\; \scrM_1^\bbR \boxtimes \scrM_2^\bbR
 \quad \textnormal{and} \quad
 \Lambda_{12} \;=\; \Lambda_1 \times \Lambda_2,
\]
we define $c_{\scrM_1, \scrM_2}$ by the following diagram: 
\[
\xymatrix@R=3.4em@C=0.25em@M=0.8em{
F(\scrM_1) * F(\scrM_2) \ar[r]^-\sim_-{\circled{1}} \ar@{..>}[dddd]_-{c_{\scrM_1, \scrM_2}}
& \varpi_* \rho^{-1}  \scrHom_{\scrE_{A^2}}(\scrC_{\Lambda_{12}}, \scrM_{12}^\bbR)|_U \ar[d]_-{\circled{2}} \\
& \varpi_* \scrHom_{\rho^{-1}\scrE_{A^2}}(\rho^{-1}\scrC_{\Lambda_{12}}, \rho^{-1}\scrM_{12}^\bbR)|_U \ar[d]_-{\circled{3}} \\ 
& \varpi_* \scrHom_{\varpi^{-1}\scrE_A}(\rho^{-1}\scrC_{\Lambda_{12}} \otimes \scrE_{A^2\to A}, \rho^{-1}\scrM_{12}^\bbR \otimes \scrE_{A^2 \to A}^\bbR )|_U \ar[d]_-{\circled{4}} \\ 
& \varpi_* \scrHom_{\varpi^{-1}\scrE_A}(\varpi^{-1}\scrC_{\Lambda}, \rho^{-1}\scrM_{12}^\bbR \otimes \scrE_{A^2 \to A}^\bbR )|_U \ar[d]_-{\circled{5}}
\\ 
F(\scrM_1 * \scrM_2) \ar[r]_-{\circled{6}}^-\sim
& \scrHom_{\scrE_A}(\scrC_{\Lambda}, \varpi_*(\rho^{-1}\scrM_{12}^\bbR \otimes \scrE_{A^2 \to A}^\bbR) )|_U 
}
\]
Here \circled{1} is induced by the K\"unneth morphism
\[
\scrHom_{\scrE_A}(\scrC_{\Lambda_1}, \scrM_1^\bbR) \boxtimes \scrHom_{\scrE_A}(\scrC_{\Lambda_2}, \scrM_2^\bbR)
\;\longrightarrow \;
\scrHom_{\scrE_{A\times A}}(\scrC_{\Lambda_1} \boxtimes \scrC_{\Lambda_2}, \scrM_1^\bbR \boxtimes \scrM_2^\bbR),
\]
while
\begin{itemize}
	\item[\circled{2}] is the natural restriction map for $\rho$,
	\item[\circled{3}] is the tensor product with the inclusion $\scrE_{A^2 \to A} \hookrightarrow \scrE_{A^2 \to A}^\bbR$,
	\item[\circled{4}] is given by the morphism $\varpi^{-1}(\scrC_\Lambda) \rightarrow \rho^{-1}(\scrC_{\Lambda_{12}}) \otimes \scrE_{A^2\to A}$ from lemma~\ref{lem:cotrace},
	\item[\circled{5}] is the adjunction isomorphism,
	\item[\circled{6}] is given by the identity $(\scrM_1 \boxtimes \scrM_2)^\bbR=\scrM_1^\bbR \boxtimes \scrM_2^\bbR$ and by the compatibility of the functor $(-)^\bbR$ with direct images on the locus where the morphism~$\varpi$ is finite~\cite[chapt.~II, th.~3.5.5]{SKK}.
\end{itemize}
Note that $c_{\scrM_1, \scrM_2}$ is a morphism between local systems of the same rank. To see that it is an isomorphism, it suffices to show that the morphisms \circled{2} -- \circled{5} are injective. This can be checked locally in the classical topology, so in what follows we will work on a classical open subset of the form $\Omega_0 = A\times U_0$ as in the proof of lemma~\ref{lem:cotrace}. Here $U_0$ can be arbitrarily small, so by the local classification of holonomic microdifferential modules in theorem~\ref{thm:microclassification} we may assume each $\scrM_i^\bbR|_{\Omega_0}$ is isomorphic to a finite direct sum of copies of the module $\scrC_{\Lambda_i}^\bbR|_{\Omega_0}$. Since \circled{2} -- \circled{5} are functorial under $\scrE_A^\bbR$-module homomorphisms and compatible with direct sums, it only remains to show injectivity when $\scrM_i = \scrC_{\Lambda_i}$. In this special case theorem~\ref{thm:microclassification} says that\medskip
\[
 \scrHom_{\scrE_{A^2}}(\scrC_{\Lambda_{12}}, \scrM_{12})|_{U\times U}
 \;=\;
 \scrHom_{\scrE_{A^2}}(\scrC_{\Lambda_{12}}, \scrM_{12}^\bbR)|_{U\times U}
 \;=\; \bbC_{\Lambda_{12}\cap (A\times U)} \medskip
\] 
is the constant local system, generated by the global section $\id: \scrC_{\Lambda_{12}}\to \scrC_{\Lambda_{12}}$. By construction the composite of the morphisms \circled{2} -- \circled{5} sends this global section to the morphism 
\[ 
  \scrC_\Lambda 
  \;\longrightarrow \; \scrC_{\Lambda_1} * \scrC_{\Lambda_2}
  \;=\;
  \varpi_*(\rho^{-1} \scrM_{12} \otimes \scrE_{A^2 \to A})
  \;\subset\; \varpi_*(\rho^{-1} \scrM_{12}^\bbR \otimes \scrE_{A^2 \to A}^\bbR)
\]
from lemma~\ref{lem:cotrace}. One easily sees from this description that the composite of \circled{2} -- \circled{5} is injective and hence an isomorphism.

\medskip

It remains to check that the $c_{\scrM_1, \scrM_2}$ satisfy the usual compatibilities with the unit, commutativity and associativity constraints of $\MM(A)$ and $\LS(A, \eta)$. Since all the morphisms in the relevant diagrams are defined globally, the commutativity of the diagrams can be checked locally in the classical topology.~For any section $f$ of~$F(\scrM_1)* F(\scrM_2)$ over a classical open neighborhood $U_0\subset V$ of a given point, we may after further shrinking the neighborhood assume that $f=\varpi_*\rho^{-1}(f_1\boxtimes f_2)$ with $f_i: \scrC_{\Lambda_i} \rightarrow \scrM_i^\bbR$ locally defined over $U_0$. But in this case the definitions imply that the local section $c_{\scrM_1, \scrM_2}(f)$ of $F(\scrM_1 * \scrM_2) = \scrHom_{\scrE_A}(\scrC_\Lambda, (\scrM_1*\scrM_2)^\bbR)|_{U_0}$ is the composite map
\[
c_{\scrM_1, \scrM_2}(f): \quad 
\xymatrix@C=3em@M=0.5em{
 \scrC_\Lambda \ar[r]^-{\iota_{\Lambda_1, \Lambda_2}}
 & \scrC_{\Lambda_1} * \scrC_{\Lambda_2} \ar[r]^-{f_1 * f_2}
 & \scrM_1^\bbR * \scrM_2^\bbR \ar[r]^-\sim
 & (\scrM_1 * \scrM_2)^\bbR
}
\]
The compatibility with the unit, commutativity and associativity constraints then reduces to simple properties of the maps $\iota$ from lemma~\ref{lem:cotrace}. For instance, the compatibility with the associativity constraints means that for any conic Lagrangian subvarieties $\Lambda_1, \Lambda_2, \Lambda_3 \subset A\times V$ that are finite \'etale over a Zariski open $U\subset V$ which is small enough so that the convolutions $\Lambda_i * \Lambda_j$ and $\Lambda_i*\Lambda_j*\Lambda_k$ are still finite \'etale over it, the following diagram commutes:
\[
\xymatrix@C=5em@M=0.6em@R=2.5em{
\scrC_{\Lambda_1*\Lambda_2*\Lambda_3, U} \ar@{=}[d] \ar[r]^-{\iota_{\Lambda_1, \Lambda_2 * \Lambda_3}} 
& \scrC_{\Lambda_1,U} * \scrC_{\Lambda_{2} * \Lambda_{3},U} \ar[r]^-{\id*\iota_{\Lambda_2,\Lambda_3}}
& \scrC_{\Lambda_1,U} * ( \scrC_{\Lambda_2,U}*\scrC_{\Lambda_3,U}) \\
\scrC_{\Lambda_1*\Lambda_2*\Lambda_3, U} \ar[r]^-{\; \iota_{\Lambda_1 * \Lambda_2, \Lambda_3}}
& \scrC_{\Lambda_{1}*\Lambda_{2},U} * \scrC_{\Lambda_3,U} \ar[r]^-{\iota_{\Lambda_1,\Lambda_2} * \id}
& (\scrC_{\Lambda_1,U} * \scrC_{\Lambda_2,U}) * \scrC_{\Lambda_3,U}  \ar[u]^-\varphi
}
\]
But this follows easily by looking at the images of the unit section. \qed

\medskip

\begin{rem} \label{rem:no-extension}
Even if $\Lambda \subset T^*A$ is the conormal variety to a smooth subvariety of~$A$, we have defined the simple modules $\scrC_\Lambda = \scrC_{\Lambda, U}$ only over the open subset $U\subseteq V$ where the Gauss map $\gamma: \Lambda \rightarrow V$ restricts to a finite \'etale cover, and in general they cannot be extended over the branch locus: If such an extension of $\scrC_\Lambda$ exists, then for any $\scrM \in \Mhol(\scrD_A)$ with $\CC(\scrM) = \Lambda$ the germ 
\[ 
 F\bigl(\pi^{-1}(\scrM)\otimes_{\scrD_A} \scrE_A \bigr) \;\in\; \LS(A, \eta) 
\] 
extends to a local system $\scrF$ on the complement of the zero section in $\Lambda$, hence its monodromy along small loops $\alpha: [0,1] \rightarrow \Lambda$ around any branch point $p\in \Lambda$ of the Gauss map must be trivial. Suppose for example that~$\gamma$ is a non-trivial cover of degree $\deg(\gamma) = 2$. Then the triviality of the local monodromy near $p$ would imply that in a suitable basis the local monodromy of $\gamma_*(\scrF)$ near the image point $\gamma(p)$ contains the matrix
$
(
\begin{smallmatrix}
0 & 1 \\
1 & 0
\end{smallmatrix}
)
$
which interchanges the stalks at the two points of the fiber $\gamma^{-1}(\gamma(p))$. Since~$F$ is a tensor functor, it would follow that this permutation matrix lies in the Tannaka group $G(\scrM)$. But after a translation we can always assume that $\scrM$ has trivial determinant, and hence $G(\scrM) \subseteq \Sl_2(\bbC)$. This shows that an extension over the branch locus cannot exist.
\end{rem}

\section{Appendix: A twistor variant} \label{sec:twistor}

In this appendix we sketch a twistor variant of the Fourier-Mukai transform that gives another source of subgroups of multiplicative type. We do not know how to see Weyl groups in this picture: The Fourier-Mukai transform for Higgs sheaves will only give vector bundles, not local systems, so we lack the monodromy operation and can only apply part (a) of theorem~\ref{thm:gradedfiberfunctor}. However, for simple holonomic $\scrD_A$-modules that do not arise from Hodge modules, the twistor construction usually yields bigger subgroups of multiplicative type than microlocalization (see example~\ref{ex:nonconic}).

\subsection{Twistor modules}

The framework of pure twistor modules~\cite{Sabbah_Polarizable, MochizukiAsymptotic, MochizukiWild} 
gives a natural way to pass from semisimple holonomic $\scrD_X$-modules on a smooth projective variety~$X$ to coherent sheaves on the cotangent bundle $T^*X$. While the resulting coherent sheaves are algebraic, the constructions pass through the analytic category. Let $\bbA^1$ be the affine line with coordinate $z$. Inside the sheaf of holomorphic relative differential operators we have the subalgebra
\[
 \scrR_X^{\,an} \;=\; 
 \bigl\langle \scrO_{X\times \bbA^1}^{\, an} \oplus z\cdot \scrT_{X\times \bbA^1 / \bbA^1}^{an} \bigr\rangle 
 \;\subset\; \scrD_{X\times \bbA^1 / \bbA^1}^{\, an}
\]
generated by the holomorphic functions and by the relative vector fields vanishing at~$z=0$. This {\em analytic Rees algebra} is a deformation from holomorphic differential operators to holomorphic fiberwise polynomial functions on the cotangent bundle in the sense that
\[
 \scrR^{\, an}_X/(z-1) \;\simeq\; \scrD_X^{\, an} \quad \textnormal{and} \quad \scrR^{\,an}_X/(z) \;\simeq\; \Sym^\bullet(\scrT_X^{an}). \]
Note that since we are working on a projective variety, the categories of holonomic respectively coherent modules for these two sheaves of rings are equivalent to their algebraic counterparts
\[
 \Mhol(\scrD_X^{\, an}) \;\simeq\; \Mhol(\scrD_X) \quad \textnormal{and} \quad
 \Coh(\Sym^\bullet(\scrT_X^{an})) \;\simeq\; \Coh(\scrO_{T^*X}).
\] 
A first interpolation between these categories is the abelian category $\Mhol(\scrR_X^{\, an})$ of holonomic right $\scrR_X^{\,an}$-modules~\cite[def.~1.2.4]{Sabbah_Polarizable}, but it is too large to have good functorial properties. Motivated by the framework of Hodge modules, Sabbah and Mochizuki introduced as a much better interpolation the abelian category of wild pure polarizable {\em twistor modules} of weight~$w\in \bbZ$. We denote this category by $\MT(X, w)^\p$ as in chapter 4 of loc.~cit., to which we refer the reader for details. We have a faithful exact functor $\MT(X, w)^\p \rightarrow \Mhol(\scrR_X^{\, an})$, denoted $\scrM \mapsto \scrM''$ in what follows,
that gives a diagram
\[
\xymatrix@M=0.5em@R=2em{
&& \Mhol(\scrD_X) && \{1\} \ar@{_{(}->}[d]^-{i_1} \\
 \MT(X,w)^\p  \ar@/^1.3pc/[urr]^-{\Xi_\DR} \ar@/_1.3pc/[drr]_-{\Xi_\Dol} \ar[r] & \Mhol(\scrR_X^{\, an}) \ar[ur]_-{i_1^*} \ar[dr]^-{i_0^*} & & \textnormal{over} & \bbA^1\\
&& \Mcoh(\scrO_{T^*X}) && \{0\} \ar@{^{(}->}[u]_{i_0}
}
\]
Here
$
 \Xi_\DR:  \MT(X, w)^\p \stackrel{\sim}{\longrightarrow} \Mhol(\scrD_X)^\semisimple \subset  \Mhol(\scrD_X)
$
is an equivalence between the abelian category of wild polarizable pure twistor modules of weight $w$ and the one of semisimple holonomic $\scrD_X$-modules~\cite[th.~19.4.1]{MochizukiWild}.

\subsection{Some functorial properties} \label{sec:laumonformula}

Even though we are ultimately interested in the pure case, it seems convenient to discuss functorial properties in the wider setup of the abelian categories $\MTM(X)$ of algebraic {\em mixed twistor modules}, for whose construction and functorial properties we refer the reader to~\cite[def.~7.2.1 and sect.~13]{MochizukiMixed}. Inside the bounded derived category of the category of mixed twistor modules we consider the full subcategory  
$\rmD_X \subset \D^b(\MTM(X))$ 
of all finite sums $\bigoplus_{w\in \bbZ} M_w[-w]$ of pure twistor modules~$M_w \in \MT(X, w)^\p$.
This subcategory contains all the pure twistor modules of weight zero and comes with the functors we want: We have an external product 
\[
 \boxtimes: \quad \rmD_X \times \rmD_X \; \longrightarrow \; \rmD_{X\times X}
\]
that is compatible with the usual external product via the functors $\Xi_\DR$ and $\Xi_\Dol$, and any morphism $f: Y\rightarrow X$ of smooth projective varieties induces a direct image functor
\[
 f_\dag: \quad \rmD_Y \; \longrightarrow \; \rmD_X
\]
that is compatible with the usual direct image functor for $\scrD$-modules via $\Xi_\DR$.
The direct image functor is also compatible with $\Xi_\Dol$ in the sense that we have an analog of Laumon's formula in~\cite[th.~2.4]{PS_GVT} and~\cite[construction 2.3.2]{Laumon_TrafoCanoniques}. For simplicity we here only verify this formula for the direct image under a projection, which will be enough for our purpose (the general case can be deduced via the graph factorization of an arbitrary morphism into a closed immersion followed by a projection):

\begin{propo} \label{prop:laumonformula}
Consider a projection $f: Y=X\times Z \to X$, and let $\varpi=\varpi_f$ and $\rho=\rho_f$ be the morphisms induced by $f$ and by its differential as in section~\ref{sec:microdifferential}, then for any $\scrM\in \rmD_Y$ we have a natural isomorphism
\[
 \Xi_\Dol \bigl( f_\dag(\scrM) \bigr) \;\stackrel{\sim}{\longrightarrow} \;
  R\varpi_{*}  L\rho^* \bigl( \, \Xi_\Dol(\scrM) \bigr)
  \quad \textnormal{\em in} \quad
  \Dbcoh(\scrO_{T^*X}).
\]
\end{propo}

{\em Proof.} Passing to direct summands we may assume that $\scrM \in \MT(Y, w)^\p$ is a single twistor module. Let $\scrM'' \in \Mhol(\scrR_Y^{\, an})$ be the underlying Rees module.
%
Then by~\cite[ex.~1.4.1]{Sabbah_Polarizable} the relative Spencer complex gives
\[
 f_\dag(\scrM'') \;\simeq \; Rf_* \bigl(\scrM'' \otimes_{\scrO_{Y\times \bbA^1}^{\, an}} \wedge^{-\bullet}(z\cdot \scrT_{Y\times \bbA^1/X\times \bbA^1}^{an})\bigr),
\]
where the derived direct image functor on the right hand side can be computed via a Dolbeault resolution as in loc.~cit. So far this holds for any $\scrM''\in \Mhol(\scrR_Y^{\, an})$. In our situation $\scrM''$ underlies a polarizable pure twistor module, so by \cite[th.~18.1.1]{MochizukiWild} all the cohomology sheaves $\scrH^i(f_\dag(\scrM''))$ underlie polarizable pure twistor modules as well. Then in particular all these cohomology sheaves are {\em strict} in the sense that the multiplication by $z$ is injective on them. The individual terms of the Dolbeault resolution which computes the direct image are also strict. So lemma~\ref{lem:strict} below gives a natural isomorphism
\[
 Li_0^*(f_\dag(\scrM'')) \;\simeq\; 
 Rf_* \bigl( i_0^*(\scrM'') \otimes_{\scrO_{Y}^{\, an}} \wedge^{-\bullet} ( \scrT_{Y/X}^{an}) \bigr)
\]
in $\Dbcoh(\scrO_{\Sym^\bullet(\scrT_X^{an})})$. The rest of the proof works in the same way as in~\cite[th.~2.4]{PS_GVT} with the appropriate translation from left to right modules. \qed

\medskip

For completeness we include the following fact, putting $\scrR = \scrR_X^{\, an}$ for brevity:

\begin{lem} \label{lem:strict}
Let $K = [\cdots \to K^\nu \to K^{\nu+1} \to \cdots ]$ be a bounded complex of strict modules in $\Mod(\scrR)$, and assume that the cohomology sheaves $\scrH^\nu(K)$ are strict for all $\nu\in \bbZ$. Then in the derived category $\rmD^b(\Mod(\scrR/(z)))$ we have a natural isomorphism
\[ K \otimes_\scrR^L \scrR/(z) \;\;\simeq\;\; \bigl[\cdots \to K^\nu/zK^\nu \to K^{\nu+1}/zK^{\nu+1} \to \cdots \bigr].
\]

\end{lem}

{\em Proof.} Let us denote by $K/zK$ the complex on the right hand side. The quotient map $q: K \rightarrow K/zK$ induces on cohomology sheaves a map that factors uniquely over maps~$r_\nu$ as indicated in the following diagram:

\[
\xymatrix@M=0.5em@C=0em@R=3em{
\scrH^\nu(K) \ar[dr] \ar[rr]^-{\scrH^\nu(q)} && \scrH^\nu(K/zK) \\
& \scrH^\nu(K) \otimes_\scrR \scrR/(z) \ar@{..>}[ur]_-{\; \exists! \, r_\nu}
}
\]
A simple diagram chase, using that the cohomology sheaf $\scrH^{\nu+1}(K)$ as well as the terms $K^\nu$ and $K^{\nu+1}$ are strict, shows that $r_\nu$ is an isomorphism. It then follows in particular that for any quasi-isomorphism $F\rightarrow K$ between complexes of strict modules in $\Mod(\scrR)$ with strict cohomology sheaves, the reduction $F/zF \rightarrow K/zK$ is still a quasi-isomorphism. Taking $F\rightarrow K$ to be a flat resolution of our given complex, we get the claimed isomorphism in the derived category. \qed

\subsection{The Fourier-Mukai transform for twistor modules}

Now let $X=A$ be an abelian variety. We want to interpret the Fourier-Mukai transform of section~\ref{subsec:fouriermukai} as the restriction to $z=1$ of a transform for twistor modules. In doing so we need to replace the moduli space $A^\natural$ by a twistor deformation $E(A)$, the moduli space of generalized flat rank one connections. Here by a {\em generalized connection} or more specifically by a {\em $\lambda$-connection} on a coherent vector bundle $\scrE\in \Mcoh(\scrO_A)$ we mean a $\bbC$-linear morphism
$
 \nabla: \scrE \rightarrow \Omega_A^1 \otimes_{\scrO_A} \scrE
$
satisfying 
\[
\nabla (f\cdot s) \;=\; df\otimes \lambda s + f\cdot \nabla{s}
\] 
for all local sections $s$ of $\scrE$ and $f$ of $\scrO_A$. We say $\nabla$ is {\em flat} if its square vanishes; this notion of a generalized flat connection makes sense for any $\lambda \in H^0(A, \scrO_A)=\bbC$. For~$\lambda = 1$ we recover a usual flat connection while for $\lambda = 0$ we get a Higgs field.

\medskip

The target space of our Fourier-Mukai transform for twistor modules will be the moduli space $E(A)$ of pairs consisting of a  line bundle $\scrL \in \hat{A} = \Pic^\circ(A)$ and a generalized flat connection~$\nabla: \scrL \to \Omega^1_A \otimes_{\scrO_A} \scrL$. For more about this moduli space we refer the reader to~\cite[sect.~10]{SchnellHolonomic}. It is a torsor over the dual abelian variety~$\hat{A}$ via the forgetful morphism 
$$\pi: \quad E(A) \; \longrightarrow \; \hat{A} \;=\; \Pic^\circ(A), \quad (\scrL, \nabla) \; \mapsto \; \scrL$$ 
and comes with a natural morphism $\lambda: E(A) \rightarrow \bbA^1$ so that the moduli spaces of flat connections in the usual sense and of Higgs fields on line bundles arise as the fibers 
\[
 \lambda^{-1}(1) \;=\; A^\natural \quad \textnormal{and} \quad \lambda^{-1}(0) \;=\; \hat{A}\times V
 \quad \textnormal{for} \quad V\;=\;H^0(A, \Omega_A^1).
\]
Let $\scrP^\lambda$ be the pullback of the Poincar\'e bundle under $\id\times \pi: A\times E(A) \to A\times \hat{A}$. By~\cite[lemma~10.5]{SchnellHolonomic} this pullback comes with a universal relative generalized flat connection
\[
 \nabla^\lambda: \quad \scrP^\lambda \;\longrightarrow \; \Omega_{A\times E(A) / E(A)}^1 \otimes_{\scrO_{A\times E(A)}} \scrP^\lambda.
\]
Taking the relative de Rham complex for the projection $p: A\times E(A) \rightarrow E(A)$, we define the Fourier-Mukai transform on the derived category of Rees modules as in loc.~cit.~by
\begin{align} \nonumber
 \FM: \quad
 \D^b(\scrR_A^{\, an}) &\;\longrightarrow \; \D^b(\scrO_{E(A)}^{\, an}), \\ \nonumber
 \scrM & \;\;\;\mapsto\;
 Rp_{*} \, \DR_{A\times E(A) / E(A)}\bigl((\id \times \lambda)^* \scrM \otimes_{\scrO_{A\times E(A)}^{\, an}} (\scrP^\lambda, \nabla^\lambda) \bigr). 
\end{align}
By definition any holonomic module $\scrM \in \Mhol(\scrR_A^{\, an})$ is {\em good} in the sense of~\cite[\S 1.1.c]{Sabbah_Polarizable} and hence admits good filtrations over the preimage of each sufficiently small open subset of the affine line. It follows that the complex $\FM(\scrM)$ has coherent cohomology sheaves. So we obtain a functor
\[
 \FM: \quad \rmD_A \;\longrightarrow \; \Dbcoh(\scrO_{E(A)}^{\, an})
\]
on the category $\rmD_A$ of pure twistor complexes of weight zero from section~\ref{sec:laumonformula}. 

\medskip

The generic vanishing property~\eqref{eq:gvt} for the Fourier-Mukai transform on $\Mhol(\scrD_A)$ carries over to the case of twistor modules as follows. Let us say that a subvariety of $E(A)$ is a {\em linear subvariety} if it has the form $E(B)\hookrightarrow E(A)$ where $A\twoheadrightarrow B$ is an epimorphism of abelian varieties. Then for every twistor module $\scrM \in \MT(A, w)^\p$ there is an open subset $U\subseteq E(A)$ whose complement is a finite union of translates of proper linear subvarieties such that
\begin{equation} \label{eq:gvt-twistor} 
 \scrH^i(\FM(\scrM))|_U
 \;\; \textnormal{is} \;\;
 \begin{cases}
 \; \textnormal{locally free} & \textnormal{for $i=0$}, \\
 \; \textnormal{zero} & \textnormal{for $i\neq 0$}.
 \end{cases} 
\end{equation}
The idea of the proof is to transfer the codimension estimates in~\cite[th.~18.1]{SchnellHolonomic} from the moduli space $A^\natural$ to all of $E(A)$ by showing that the restriction of the Fourier-Mukai transform to each twistor line has locally free cohomology sheaves.

\medskip

\subsection{Another fiber functor}

By the above twistor reformulation we may apply the constructions from section~\ref{subsec:fouriermukai} not only at $z=1$ but also at $z=0$. 
In this case our Fourier-Mukai transform is just the relative Fourier-Mukai transform for coherent sheaves on the abelian scheme $A\times V \rightarrow V$ (see the proof of theorem~\ref{thm:torus-twistor} below).

\medskip

To see what we have gained from this, let $U\subseteq E(A)$ be the complement of a finite union of translates of proper linear subvarieties, and consider  $U_z = U\cap \lambda^{-1}(z)$ for~$z\in \bbA^1(\bbC)$. Let 
\[ \Mhol(A, U_1)^\semisimple \; \subset \; \Mhol(\scrD_A, U_1) \]
be the full subcategory of all semisimple holonomic modules $\scrN \simeq \Xi_\DR(\scrM)$ for which the corresponding twistor module $\scrM \in \MT(A, 0)^\p$ satisfies the vanishing condition~\eqref{eq:gvt-twistor} and the Gauss map
\[
\xymatrix@M=0.7em{
 \gamma: \quad \Supp(\Xi_\Dol(\scrM)|_{A\times U_0}) \;\subset \; A\times U_0  \ar@{->>}[r]^-p & U_0
}
\]
is a finite morphism. We view $\Mhol(\scrD_A, U_1)^\semisimple$ as a tensor category with respect to the 
tensor product obtained from the convolution product by discarding all negligible direct summands. Then the above discussion leads to the following result, where for objects $\scrE \in \VB(A, U_0)$ in the tensor category of section~\ref{sec:vectorbundles} and for $u\in U_0(\bbC)$ we denote by
\[
 X_u(\scrE) \;=\; \langle a\in A(\bbC) \mid (a, u) \in \Supp(\scrE) \rangle \;\subset\; A(\bbC)
\]
the group generated by the points in the fiber of the corresponding Gauss map.

\begin{thm} \label{thm:torus-twistor}
For $U\subseteq E(A)$ as above, we have a commutative diagram of tensor functors
\[
\xymatrix@M=0.5em@C=2em{
\Mhol(\scrD_A, U_1)^\semisimple \ar@{..>}[dr]_-{\exists \, \Phi} \ar[rr]^-{\scrH^0(\FM(-))|_{U_0}} && \Mcoh(\scrO_{U_0}^{\, an}) \\
& \VB(A, U_0) \ar[ur]_-{p_*}  &
}
\]
Taking the fiber at $u\in U_0(\bbC)$, we obtain for any $\scrM \in \Mhol(\scrD_A, U_1)^\semisimple$ an embedding
\[ \CartierDual{X_u} \hookrightarrow G(\scrM)
\quad \textnormal{\em where} \quad
X_u \;=\; X_u(\Phi(\scrM)).
\]
\end{thm}

{\em Proof.}
Via the functor $\Xi_\DR$ we consider $\Mhol(\scrD_A, U_1)^\semisimple$ as a full subcategory of~$\MT(A, 0)^\p$. Taking the inverse image under the projection
$q:  A\times U_0  \rightarrow A\times V$
followed by the involution $\iota = \id_A \times (-\id_V)$ and tensoring with the restriction of the Poincar\'e line bundle
\[
 \scrP_{U_0} \;:=\; \scrP^\lambda|_{A\times U_0},
\]
we define
\[
 \Phi \;=\; q^* \iota^* (\Xi_\Dol(-)) \otimes \scrP_{U_0}: \quad
 \Mhol(\scrD_A, U_1)^\semisimple \;\longrightarrow \; \VB(A, U_0). \medskip
\]
For $\scrM_1, \scrM_2 \in \Mhol(\scrD_A, U_1)^\semisimple $ we get a natural isomorphism\medskip 
\[
 \Phi(\scrM_1 * \scrM_2) \; \simeq \; \Phi(\scrM_1) * \Phi(\scrM_2)
\]
by applying proposition~\ref{prop:laumonformula} to the direct image of $\scrM=\scrM_1 \boxtimes \scrM_2$ under the sum morphism $f: A\times A \to A$ (which can be considered as a projection onto a direct factor by precomposing it with the automorphism $(x,y)\mapsto (x,x+y)$) and then using the projection formula and the fact that for the Poincar\'e bundle we have a natural isomorphism
\[ \varpi_f^*(\scrP_{U_0}) \; \simeq \; \rho_f^*(\scrP_{U_0} \boxtimes \scrP_{U_0}). \]
One may check that these isomorphisms $\Phi(\scrM_1 * \scrM_2) \simeq \Phi(\scrM_1) * \Phi(\scrM_2)$ satisfy the usual compatibilities with the associativity, commutativity and unit constraints, so that~$\Phi$ becomes a tensor functor. We also note that 
\[
 Rp_*(q^*(-)\otimes \scrP_{U_0}):
 \quad
 \Dbcoh(\scrO_{A\times V}) \;\longrightarrow \; \Dbcoh(\scrO_{\hat{A} \times V})
\] 
is the Fourier-Mukai transform between the derived categories of coherent sheaves on the abelian scheme $A\times V \rightarrow V$ and its dual. The twist by the involution $\iota$ makes the diagram of tensor functors commutative by the same base change arguments as in~\cite[proof of prop.~12.1(iii)]{SchnellHolonomic}, and the embedding $\CartierDual{X_u} \hookrightarrow G(\scrM)$ then comes from the first part of theorem~\ref{thm:gradedfiberfunctor}.  
\qed

\medskip

\begin{ex} \label{ex:nonconic}
If $\scrM \in \Mhol(\scrD_A)$ underlies a pure Hodge module, then theorem~\ref{thm:torus-twistor} yields the same subgroup of multiplicative type in $G(\scrM)$ as our previous microlocal construction in theorem~\ref{thm:microlocal-fiberfunctor}: Applying the functor from Hodge modules to twistor modules in~\cite[\S 13.5]{MochizukiMixed}, one sees that here the twistor module associated to $\scrM$ has as its $\scrR_A^{\, an}$-module the Rees module \medskip 
\[
 R_F(\scrM) \;:=\; \bigoplus_{k\in \bbZ} F_k(\scrM)\cdot z^k \;\subseteq\; \scrM \otimes \bbC[z, z^{-1}]
\]
for the Hodge filtration. Then $\Xi_\mathrm{Dol}(\scrM)\simeq gr_F(\scrM)$ is the associated graded for the Hodge filtration. The latter is a good filtration, hence $\Supp(\Xi_\mathrm{Dol})(\scrM)=\Char(\scrM)$ and the subgroup of multiplicative type in theorem~\ref{thm:torus-twistor} arises from the characteristic cycle as in our previoius microlocal construction.

\medskip 

However, theorem~\ref{thm:torus-twistor} gives more for semisimple holonomic $\scrD_A$-modules that do not come from pure Hodge modules, in particular for modules with irregular singularities like those in irregular Hodge theory~\cite{Sabbah_Irregular}. For such $\scrD_A$-modules $\scrM \in \Mhol(\scrD_A)$ one usually has
\[
 \Supp(\Phi(\scrM)) \;\neq\; \Char(\scrM),
\]
and the subgroup of multiplicative type in theorem~\ref{thm:torus-twistor} can be strictly bigger than the one in theorem~\ref{thm:gausstorus}. Suppose for example that $A$ is an elliptic curve, and take a holomorphic function $f: A^*=A\setminus \{0\} \rightarrow \bbC$ with a pole at the origin. Consider the trivial smooth bundle $\scrC^\infty_{A^*}$ endowed with the flat connection $\nabla = d + \partial f+\overline{\partial f}$ and the constant Hermitian metric $h(1,1) = 1$. This is a harmonic bundle of rank one in the sense of~\cite[sect.~1.2.1.1]{MochizukiWild}; the corresponding Higgs bundle is the trivial line bundle $\scrO_{A^*}^{\, an}$ equipped with the Higgs field $\theta = \partial f$. By loc.~cit.~our chosen harmonic bundle is the twistor deformation of a unique simple module $\scrM^* \in \Mhol(\scrD_{A^*})$. The minimal extension $\scrM \in \Mhol(\scrD_A)$ of the latter has its characteristic variety contained in the zero section plus the fiber over the origin. On the other hand the support $\Supp(\Phi(\scrM))\subset A\times V$ coincides over the open subset $A^*=A\setminus \{0\} \subset A$ with the graph of the differential $\partial f: A^* \to V = H^0(A, \Omega^1_A)$. It follows that for the fibers of the projection $\gamma: A\times V \twoheadrightarrow V$ over any non-zero point $u\in V(\bbC)$ we have
\begin{eqnarray*}
\Char(\scrM) \cap \gamma^{-1}(u) 
&=& 
\{ (0, u)\} \\ 
\Supp(\Phi(\scrM))\cap \gamma^{-1}(u)
&\supseteq & \{(a,u)) \mid f(a) = u\}.
\end{eqnarray*} 
Hence we see that while the subgroup of multiplicative type from theorem~\ref{thm:gausstorus} is trivial in this case, the one from theorem~\ref{thm:torus-twistor} is not: More precisely, it is the Cartier dual of the group generated by the finitely many points $a\in A^*(\bbC)$ with $f(a)=u$.
\end{ex}

\medskip

\subsection*{Acknowledgements.} I would like to thank Claude Sabbah for many inspiring discussions and his kind hospitality at \'Ecole Polytechnique, Christian Schnell for explaining to me his conjectures on the Fourier-Mukai transform that motivated the paper, Takuro Mochizuki and Will Sawin for various exchanges, and the referees for their many useful comments and suggestions. 
This project was funded by the DFG research grant {\em Holonome $\scrD$-Moduln auf abelschen Variet\"aten}.

\bigskip

\backmatter

\bibliographystyle{smfplain}
\bibliography{Bibliography}

\providecommand{\bysame}{\leavevmode ---\ }
\providecommand{\og}{``}
\providecommand{\fg}{''}
\providecommand{\smfandname}{\&}
\providecommand{\smfedsname}{\'eds.}
\providecommand{\smfedname}{\'ed.}
\providecommand{\smfmastersthesisname}{M\'emoire}
\providecommand{\smfphdthesisname}{Th\`ese}
\begin{thebibliography}{10}

\bibitem{AndreottiTorelli}
{\scshape {Andreotti, A.}} -- {\og {On a theorem of Torelli}\fg}, \emph{{Amer.
  J. Math.}} \textbf{{80}} ({1958}), p.~{801--828}.

\bibitem{ABPerverse}
{\scshape {Arinkin, D. and Bezrukavnikov, R.}} -- {\og {Perverse coherent
  sheaves}\fg}, \emph{{Mosc. Math. J.}} \textbf{{10}} ({2010}), p.~{3--29}.

\bibitem{BKHolonomy}
{\scshape {Balaji, V. and Koll\'ar, J.}} -- {\og {Holonomy groups of stable
  vector bundles}\fg}, \emph{{Publ. RIMS, Kyoto Univ.}} \textbf{{44}} ({2008}),
  p.~{183--211}.

\bibitem{BeilinsonHolonomic}
{\scshape {Beilinson, A.}} -- {\og {Constructible sheaves are holonomic}\fg},
  \emph{{Selecta Math.}} \textbf{{22}} ({2016}), p.~{1797–1819}.

\bibitem{BogomolovStable}
{\scshape {Bogomolov, F. A.}} -- {\og {Stable vector bundles on projective
  surfaces}\fg}, \emph{{Russian Acad. Sci. Sb. Math.}} \textbf{{81}} ({1995}),
  p.~{397--419}.

\bibitem{BuehlerExact}
{\scshape {B\"uhler, T.}} -- {\og {Exact categories}\fg}, \emph{{Expo. Math.}}
  \textbf{{28}} ({2010}), p.~{1--69}.

\bibitem{CarterFinite}
{\scshape {Carter, R. W.}} -- \emph{{Finite groups of Lie type --- Conjugacy
  classes and complex characters}}, {Wiley}, {1985}.

\bibitem{CGIntermediateJacobian}
{\scshape {Clemens, C. H. and Griffiths, P. A.}} -- {\og {The intermediate
  Jacobian of the cubic threefold}\fg}, \emph{{Annals of Math.}} \textbf{95}
  (1972), p.~281--356.

\bibitem{CommelinMTC}
{\scshape {Commelin, J.}} -- {\og {The Mumford-Tate conjecture for products of
  abelian varieties}\fg}, \emph{{Algebraic Geometry}} \textbf{{6}} ({2019}),
  p.~{650--677}.

\bibitem{ATLAS}
{\scshape {Conway, J. H. et al.}} -- \emph{{ATLAS of finite groups}}, Clarendon
  Press, 1985.

\bibitem{CoulembierTannakian}
{\scshape {Coulembier, K.}} -- {\og {Tannakian categories in positive
  characteristic}\fg}, \emph{{Duke Math. J.}} \textbf{{169}} ({2020}).

\bibitem{DeligneHodgeCycles}
{\scshape {Deligne, P.}} -- {\og {Hodge cycles on abelian varieties}\fg}, in
  \emph{{Hodge cycles, Motives and Shimura varieties}}, Lecture Notes in Math.,
  vol. 900, Springer Verlag, 1982, p.~9--100.

\bibitem{DM}
{\scshape {Deligne, P. and Milne, J. S.}} -- {\og {Tannakian categories}\fg},
  in \emph{Hodge Cycles, Motives, and Shimura varieties}, Lecture Notes in
  Math., vol. 900, Springer Verlag, 1982, p.~\mbox{101--228}.

\bibitem{EtingofGelaki}
{\scshape {Etingof, P., Gelaki, S., Nikshych, D. and Ostrik, V.}} --
  \emph{{Tensor categories}}, {Math. Surveys Monogr.}, vol. {205}, {Amer. Math.
  Soc.}, {2015}.

\bibitem{FKGauss}
{\scshape {Franecki, J. and Kapranov, M.}} -- {\og {The Gauss map and a
  noncompact Riemann-Roch \mbox{formula} for constructible sheaves on
  semiabelian varieties}\fg}, \emph{{Duke Math. J.}} \textbf{104} (2000),
  p.~171--180.

\bibitem{GaL}
{\scshape {Gabber, O. and Loeser, F.}} -- {\og {Faisceaux pervers
  $\ell$-adiques sur un tore}\fg}, \emph{Duke Math. J.} \textbf{83} (1996),
  p.~501--606.

\bibitem{GinzburgCharacteristic}
{\scshape {Ginzburg, V.}} -- {\og {Characteristic varieties and vanishing
  cycles}\fg}, \emph{Invent. Math.} \textbf{84} (1986), p.~327--402.

\bibitem{GreenMinuscule}
{\scshape {Green, R. M.}} -- \emph{{Combinatorics of minuscule
  representations}}, {Cambridge Tracts in Math.}, vol. {199}, {Cambridge Univ.
  Press}, {2013}.

\bibitem{HartshorneStable}
{\scshape {Hartshorne, R.}} -- {\og {Stable reflexive sheaves}\fg},
  \emph{{Math.~Annalen}} \textbf{{254}} ({1980}), p.~{121--176}.

\bibitem{HTT}
{\scshape {Hotta, R., Takeuchi, K. and Tanisaki, T.}} -- \emph{{D-modules,
  perverse sheaves and representation theory}}, Progress in Math., vol. 236,
  {Birkh{\"a}user}, 2008.

\bibitem{KashiwaraSystems}
{\scshape {Kashiwara, M.}} -- \emph{{Systems of microdifferential equations}},
  {Progress in Math.}, vol.~{34}, {Birkh\"auser}, {1983}.

\bibitem{KashiwaraRepresentation}
\bysame , {\og {Representation theory and $\scrD$-modules on flag
  varieties}\fg}, in \emph{{Orbites unipotentes et réprésentations III}},
  vol. {173-174}, {1989}, p.~{55--109}.

\bibitem{KashiwaraDModules}
\bysame , \emph{{$\mathscr{D}$-modules and microlocal calculus}}, {Iwanami
  Series in Modern Mathematics, Transl. Math. Monographs}, vol. {217},
  {American Mathematical Society}, {2003}.

\bibitem{KashiwaraTStructures}
\bysame , {\og {t-structures on the derived categories of holonomic
  $\mathscr{D}$-modules and coherent $\scrO$-modules}\fg}, \emph{{Mosc. Math.
  J.}} \textbf{{4}} ({2004}), p.~{847–868}.

\bibitem{KatzSatoTate}
{\scshape {Katz, N. M.}} -- \emph{{Convolution and equidistribution: Sato-Tate
  theorems for finite-field Mellin transforms}}, {Ann. of Math. Stud.}, vol.
  {180}, {Princeton Univ. Press}, {2012}.

\bibitem{KraemerSemiabelian}
{\scshape {Kr{\"a}mer, T.}} -- {\og {Perverse sheaves on semiabelian
  varieties}\fg}, \emph{{Rend. Semin. Mat. Univ. Padova}} \textbf{{132}}
  ({2014}), p.~{83--102}.

\bibitem{KraemerThreefolds}
\bysame , {\og {Cubic threefolds, Fano surfaces and the monodromy of the Gauss
  map}\fg}, \emph{{Manuscripta Math.}} \textbf{{149}} ({2016}), p.~{303--314}.

\bibitem{KraemerMicrolocalII}
\bysame , {\og {Characteristic cycles and the microlocal geometry of the Gauss
  map II}\fg}, \emph{{J. Reine Angew. Math.}} \textbf{{774}} ({2021}),
  p.~{53--92}.

\bibitem{KrWSmall}
{\scshape {Kr{\"a}mer, T. and Weissauer, R.}} -- {\og {Semisimple super
  Tannakian categories with a small tensor generator}\fg}, \emph{{Pacific J.
  Math.}} \textbf{{276}} ({2015}), p.~{229--248}.

\bibitem{KrWSchottky}
\bysame , {\og {The Tannaka group of the theta divisor on a generic principally
  polarized abelian variety}\fg}, \emph{{Math. Z.}} \textbf{{281}} ({2015}),
  p.~{723--745}.

\bibitem{KrWVanishing}
\bysame , {\og {Vanishing theorems for constructible sheaves on abelian
  varieties}\fg}, \emph{{J. Alg. Geom.}} \textbf{{24}} ({2015}), p.~{531--568}.

\bibitem{LaumonTransformation}
{\scshape {Laumon, G.}} -- {\og {Transformation de Fourier
  g\'en\'eralis\'ee}\fg}, {\url{arXiv:alg-geom/9603004}}.

\bibitem{Laumon_TrafoCanoniques}
\bysame , {\og {Transformations canoniques et sp\'ecialisation pour les D
  -modules filtr\'es}\fg}, \emph{{Differential Systems and Singularities
  (Luminy, 1983), Ast\'erisque}} \textbf{{130}} ({1985}), p.~{56–129}.

\bibitem{MaximSaitoSchuermann}
{\scshape {Maxim, L., Saito, M. and Sch\"urmann, J.}} -- {\og {Symmetric
  products of mixed Hodge modules}\fg}, \emph{{J. Math. Pures Appl.}}
  \textbf{{96}} ({2011}), p.~{462--483}.

\bibitem{MazurMessing}
{\scshape {Mazur, B. and Messing, W.}} -- in \emph{{Universal Extensions and
  One Dimensional Crystalline Cohomology}}, {Lecture Notes in Math.}, vol.
  {370}, {Springer}, {1974}.

\bibitem{MochizukiAsymptotic}
{\scshape {Mochizuki, T.}} -- {\og {Asymptotic behaviour of tame harmonic
  bundles and an application to pure twistor $D$-modules}\fg}, \emph{Mem. Amer.
  Math. Soc.} \textbf{185} (2007).

\bibitem{MochizukiWild}
\bysame , \emph{{Wild Harmonic Bundles and Wild Pure Twistor
  $\mathscr{D}$-modules}}, {Ast{\'e}risque}, vol. {340}, {Soci{\'e}t{\'e}
  Math{\'e}matique de France}, {2011}.

\bibitem{MochizukiMixed}
\bysame , \emph{{Mixed twistor $\scrD$-modules}}, {Lecture Notes in Math.},
  vol. {2125}, {Springer Verlag}, {2015}.

\bibitem{NoriRepresentations}
{\scshape {Nori, M. V.}} -- {\og {On the representations of the fundamental
  group}\fg}, \emph{{Compos. Math.}} \textbf{33} (1976), p.~29--41.

\bibitem{PS_GVT}
{\scshape {Popa, M. and Schnell, C.}} -- {\og {Generic vanishing theory via
  mixed Hodge modules}\fg}, \emph{{Forum of Math. Sigma}} \textbf{{1}}
  ({2013}), p.~{1--60}.

\bibitem{RaynaudTorsion}
{\scshape {Raynaud, M.}} -- {\og {Courbes sur une vari\'et\'e ab\'elienne et
  points de torsion}\fg}, \emph{{Invent. Math.}} \textbf{{71}} ({1983}),
  p.~{207--233}.

\bibitem{RothsteinSheaves}
{\scshape {Rothstein, M.}} -- {\og {Sheaves with connection on abelian
  varieties}\fg}, \emph{{Duke Math. J.}} \textbf{{84}} ({1996}), p.~{565--598}.

\bibitem{Sabbah_Polarizable}
{\scshape {Sabbah, C.}} -- {\og {Polarizable Twistor
  $\mathscr{D}$-modules}\fg}, \emph{{Ast{\'e}risque}} \textbf{300} (2005).

\bibitem{Sabbah_Irregular}
\bysame , \emph{{Irregular Hodge theory}}, {Memoires de la SMF}, vol. {156},
  {2018}.

\bibitem{SaitoCC}
{\scshape {Saito, T.}} -- {\og {The characteristic cycle and the singular
  support of a constructible sheaf}\fg}, \emph{{Invent. Math.}} \textbf{{207}}
  ({2017}), p.~{597--695}.

\bibitem{SKK}
{\scshape {Sato, M., Kashiwara, M. and Kawai, T.}} -- {\og {Hyperfunctions and
  Pseudo-differential Equations II}\fg}, Lecture Notes in Math., vol. {287},
  {Springer Verlag}, {1973}, p.~{265--529}.

\bibitem{SchapiraMicrodifferential}
{\scshape {Schapira, P.}} -- \emph{{Microdifferential Systems in the Complex
  Domain}}, {Grundlehren Math. Wiss.}, vol. {269}, {Springer Verlag}, {1985}.

\bibitem{SSIndex}
{\scshape {Schapira, P. and Schneiders, J.-P.}} -- {\og {Elliptic Pairs I.
  Relative Finiteness and Duality}\fg}, \emph{{Ast\'erisque}} \textbf{{224}}
  ({1994}), p.~{5--60}.

\bibitem{SchnellHolonomic}
{\scshape {Schnell, C.}} -- {\og {Holonomic $\mathscr{D}$-modules on abelian
  varieties}\fg}, \emph{{Publ. Math. Inst. Hautes {\'E}tudes Sci.}}
  \textbf{{121}} ({2015}), p.~{1--55}.

\bibitem{SchreiederTheta}
{\scshape {Schreieder, S.}} -- {\og {Theta divisors with curve summands and the
  Schottky problem}\fg}, \emph{{Math. Annalen}} \textbf{{365}} ({2016}),
  p.~{1017--1039}.

\bibitem{SpringerVeldkamp}
{\scshape {Springer, T. A. and Veldkamp, F. D.}} -- \emph{{Octonions, Jordan
  Algebras and Exceptional Groups}}, {Springer Verlag}, {2000}.

\bibitem{WeBN}
{\scshape {Weissauer, R.}} -- {\og {Brill-Noether sheaves}\fg},
  {\url{arXiv:math/0610923}}.

\bibitem{WeissauerDegenerate}
\bysame , {\og {On subvarieties of abelian varieties with degenerate Gauss
  mapping}\fg}, {\url{arXiv:1110.0095}}.

\bibitem{WeConnected}
\bysame , {\og {Why certain Tannaka groups attached to abelian varieties are
  almost connected}\fg}, {\url{arXiv:1207.4039}}.

\bibitem{WeissauerVanishing}
\bysame , {\og {Vanishing theorems for constructible sheaves on abelian
  varieties over finite fields}\fg}, \emph{{Math. Annalen}} \textbf{{365}}
  ({2016}), p.~{559--578}.

\end{thebibliography}

\end{document}